\documentclass[11pt]{article}

\usepackage[utf8]{inputenc}
\usepackage{amssymb,amsmath,amsthm,mathtools}
\usepackage{array}
\usepackage{graphicx}
\usepackage{xcolor}

\usepackage{tikz}
\usepackage[numbers]{natbib}
\usepackage[bookmarks=true,pdfstartview=FitH,colorlinks,
            linkcolor=black,filecolor=black,citecolor=black,
            urlcolor=black]{hyperref}
\usepackage{aliascnt}

\usepackage{float}     
\usepackage{cleveref}  

\usetikzlibrary{arrows.meta}

\allowdisplaybreaks[2]

\newtheorem{theorem}{Theorem}[section]

\newaliascnt{lemma}{theorem}
\newtheorem{lemma}[lemma]{Lemma}
\aliascntresetthe{lemma}

\newaliascnt{example}{theorem}

\aliascntresetthe{example}

\newaliascnt{claim}{theorem}
\newtheorem{claim}[claim]{Claim}
\aliascntresetthe{claim}

\newaliascnt{corollary}{theorem}
\newtheorem{corollary}[corollary]{Corollary}
\aliascntresetthe{corollary}

\newaliascnt{construction}{theorem}

\aliascntresetthe{construction}

\newaliascnt{question}{theorem}

\aliascntresetthe{question}

\newaliascnt{fact}{theorem}

\aliascntresetthe{fact}

\newaliascnt{proposition}{theorem}
\newtheorem{proposition}[proposition]{Proposition}
\aliascntresetthe{proposition}

\newaliascnt{conjecture}{theorem}

\aliascntresetthe{conjecture}

\newaliascnt{definition}{theorem}
\newtheorem{definition}[definition]{Definition}
\aliascntresetthe{definition}

\newaliascnt{remark}{theorem}

\aliascntresetthe{remark}

\newaliascnt{observation}{theorem}
\newtheorem{observation}[observation]{Observation}
\aliascntresetthe{observation}

\usepackage{cleveref}
\crefname{theorem}{Theorem}{Theorems}
\crefname{lemma}{Lemma}{Lemmas}
\crefname{example}{Example}{Examples}
\crefname{claim}{Claim}{Claims}
\crefname{corollary}{Corollary}{Corollaries}
\crefname{construction}{Construction}{Constructions}
\crefname{question}{Question}{Questions}
\crefname{fact}{Fact}{Facts}
\crefname{proposition}{Proposition}{Propositions}
\crefname{conjecture}{Conjecture}{Conjectures}
\crefname{definition}{Definition}{Definitions}
\crefname{remark}{Remark}{Remarks}
\crefname{observation}{Observation}{Observations}
\crefname{subsection}{Subsection}{Subsections}

\newcommand{\PP}{\mathbb{P}}

\newcommand{\Lip}{\operatorname{Lip}}
\newcommand{\dist}{\operatorname{dist}}

\title{Matching Random Point Sets with Minimum Distortion}
\author{
  Itai Benjamini\thanks{
    Faculty of Mathematics and Computer Science,
    Weizmann Institute of Science, Rehovot, Israel.
  }
  \and
  Mark Braverman\thanks{
    Department of Computer Science,
    Princeton University, Princeton, NJ, USA.
  }
  \and
  Elad Tzalik\thanks{
    Faculty of Mathematics and Computer Science,
    Weizmann Institute of Science, Rehovot, Israel. Supported by the Adams Fellowship Program of the Israel Academy of Sciences and Humanities.}
  \and
  Georgii Zakharov\thanks{
    Mathematical Institute,
    University of Oxford, Oxford, UK.
  }
}

\usepackage{caption}
\date{}

\begin{document}

\maketitle

\begin{abstract}
How well can two independent random samples from a metric measure space be matched while preserving pairwise distances? Concretely, given two independent copies, each of $n$ points sampled from the space, what is their typical metric distortion as finite $n$-point metric spaces? We prove that the space $[0,1]$ equipped with the uniform measure already exhibits a rich behavior. We show that with high probability the distortion of the two independent samples is $n^{2/3}$ up to polylogarithmic factors.

The lower bound reflects two competing obstructions. A bijection that maps points close to their original positions or their reflections in $[0,1]$ must expand the smallest gaps, in both directions. On the other hand, we show that large displacement from both the identity and reflection itself forces large distortion. Balancing these obstructions yields the lower-bound exponent $2/3$. For the upper bound, we construct a partial matching of the points that maps each gap between consecutive points of length $d\lesssim n^{-4/3}$ to a pair at distance $\widetilde{O}(dn^{1/3})$, with analogous control in the reverse direction. We then use a delicate multiscale oscillation construction to extend this partial matching to a low-distortion bijection.
\end{abstract}

\noindent
\begin{minipage}{\textwidth}
  \centering
  \includegraphics[width=\linewidth]{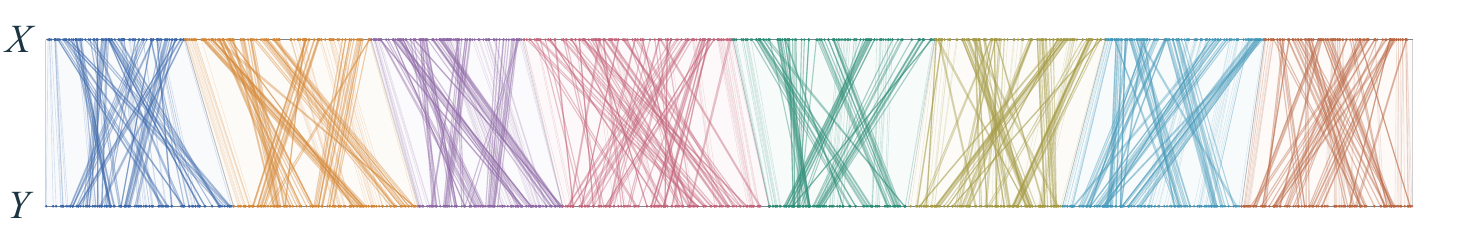}
  \captionsetup{font=small,skip=3pt,hypcap=false}
  \captionof{figure}{Illustration of two independent samples of $1{,}000$ points
  in $[0,1]$, matched using a near-optimal
  mapping. Each color identifies a pair of corresponding intervals
  whose points are matched bijectively.}
  \label{fig:local-matching}
\end{minipage}
\par\smallskip
\noindent
\begin{minipage}{\textwidth}
  \centering
  \includegraphics[width=\linewidth]{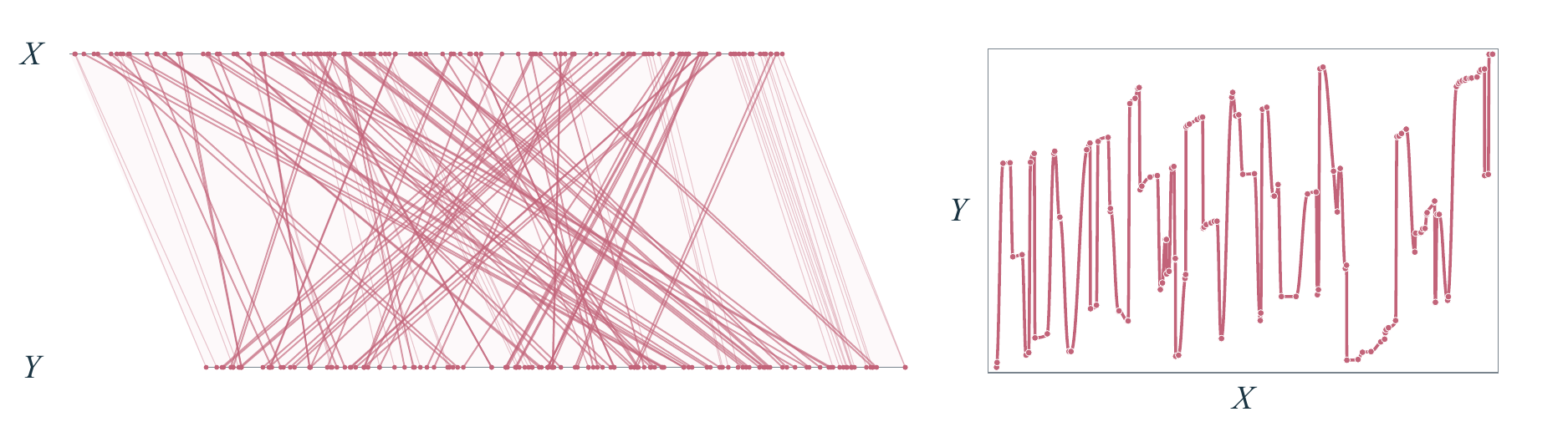}
  \captionsetup{font=small,skip=3pt,hypcap=false}
  \captionof{figure}{An enlarged view of the pink interval.
  Left: the restricted matching.
  Right: the matched pairs $(x,f(x))$, joined in increasing order of $x$. Dots represent the actual matched pairs.}
  \label{fig:local-oscillations}
\end{minipage}
\par

\clearpage

\section{Introduction}
\label{sec:introduction}

How well can two independent random samples $X,Y$ of $n$ points from a
metric measure space be matched while preserving pairwise distances?
We study this question through the minimum distortion between $X$ and
$Y$, viewed as finite metric spaces with the induced metric. Our aim is
to understand the typical growth of this distortion with the sample
size.

For finite metric spaces $X,Y$ of the same cardinality with at least
two points, and a bijection $f:X\to Y$, define
\[
 \begin{gathered}
  \Lip(f)=\max_{x\ne x'}
  \frac{d_Y(f(x),f(x'))}{d_X(x,x')},
  \qquad
  \dist(f)=\Lip(f)\Lip(f^{-1}),\\
  c(X,Y)=\min\{\dist(f):f:X\to Y\text{ is bijective}\}.
 \end{gathered}
\]
The distortion of a bijection is the product of its maximum expansion
and contraction factors. In particular, it is invariant under
rescaling either metric. Thus $c(X,Y)$ measures how closely the two
metric spaces can be matched.

The minimum distortion between two independent samples can be thought
of as a notion of \emph{metric variance}. It is well known that, for
a square-integrable real random variable $A$ and an independent copy
$A'$, we have
$\operatorname{Var}(A)=\frac12\mathbb E[(A-A')^2]$.
This identity expresses variance through the discrepancy between
independent copies, without referring to their mean. Similarly,
comparing $X$ and $Y$ through $c(X,Y)$ provides a multiplicative measure
of the variability of a random finite metric space, without requiring
a definition of its expected metric space.

\subsection{Main result}

In this work, we study the minimum distortion between two copies, $X$ and $Y$, each consisting of $n$ points drawn uniformly and independently from $[0,1]$. We show that even this one-dimensional setting is highly non-trivial and exhibits rich geometric structure.
Our main result is the following theorem:

\begin{theorem}
\label{thm:main}
There exists an absolute constant $C>0$ such that, with high probability,
\[
 \frac{n^{2/3}}{(\log n)^C}
 \le c(X,Y)\le
 n^{2/3}(\log n)^C.
\]
\end{theorem}

\paragraph{Why typical distortion?}
The expected distortion is a natural first statistic to consider. In our setting, however, it fails
to capture the typical behavior, as 
$\mathbb E[c(X,Y)]=\infty$ for sufficiently large $n$ (see Appendix \ref{sec:inf-exp}).
The reason is simple: an exceptionally close pair in one sample, with
no comparably close pair in the other, forces large distortion under
every bijection. Such rare configurations produce a sufficiently
heavy tail to make the expectation infinite. We therefore study the
order of magnitude distortion encountered with high probability when comparing
two random configurations.

\paragraph{How (not) to match points.}
The simplest bijection matches points in order. Writing
$X=\{x_1<\cdots<x_n\}$ and $Y=\{y_1<\cdots<y_n\}$, it sends
$x_i$ to $y_i$. The difficulty with this matching is that the smallest
gaps in the two samples need not occur at corresponding positions.
With high probability, the smallest gap in $X$ has size
$n^{-2+o(1)}$, while the gap at the same position in $Y$ has size
$n^{-1+o(1)}$. This pair alone requires expansion by a factor
$n^{1+o(1)}$. The smallest gap in $Y$ causes the same problem for
the inverse map. Thus both Lipschitz constants are
$n^{1+o(1)}$, giving distortion $n^{2+o(1)}$ with high probability.

This suggests matching the closest pair in $X$ to the closest pair in
$Y$, with both gaps of length $n^{-2+o(1)}$ whp, and extending this partial matching to a
bijection with low distortion. This controls the distance ratio for
the selected pair, but leaves a geometric constraint on all remaining
points. The two pairs can lie in very different parts of the interval:
for example, one may lie near an endpoint and the other near the middle.
They then divide $[0,1]$ into parts containing very different numbers of
points, which need to be matched. We prove that, whp, every bijection matching the minimum gap pair in $X$ to the minimum gap pair in $Y$ has large distortion of order at least $n^{1-o(1)}$.

\paragraph{The main difficulty.}
A more flexible approach is to match each exceptionally short gap to
a suitable gap within a limited distance, without insisting that its
partner be the shortest gap in the other sample. The difficulty is
then to choose these matches so that they can be completed to a
bijection with low distortion. Although exceptionally short gaps involve
only a small fraction of the points, matching them with low distortion
requires a carefully coordinated mapping of the remaining points, as
illustrated in the left panel of \Cref{fig:local-oscillations}.
In particular, many neighboring pairs at the typical distance of order
$1/n$, which initially posed no obstacle, are mapped much farther apart
to accommodate the displacement required by the short-gap matches.
The oscillations distribute this displacement across the remaining
points while controlling expansion in both directions.
These oscillations must be coordinated in both samples: nearby points
of $X$ must have nearby images in $Y$, and nearby points of $Y$ must
have nearby preimages in $X$, within the allowed expansion factor.
Searching farther for partners makes short gaps easier to match, but
also increases the displacement that these oscillations must
accommodate. Balancing these two requirements yields the exponent $2/3$.
This balance identifies the relevant scales, but realizing it requires
matching nearby short gaps together, resolving conflicts between
prescribed matches, and completing the resulting partial matching
while controlling expansion in both directions. All these components
of the mapping are delicate and have nontrivial interactions. We
survey them in the technical overview, \Cref{sec:proof-overview}.
\subsection{Related work}

\paragraph{Geometric similarity of random point sets in $\mathbb R$.}
A closely related line of work asks whether two independent
realizations of the same random process have the same geometry
at large scales. 
Peled~\cite{Peled2010} studied quantitative versions of this question
for random subsets of the line. Basu and Sly~\cite{BasuSly2014}
subsequently developed a multiscale construction for embeddings of
random sequences, resolving the question for Poisson processes on
the line: two independent realizations are almost surely roughly
isometric. Rough isometries allow both multiplicative and additive
errors in distances. We highlight that additive error allows exceptionally short
gaps to be mapped to substantially longer ones, whereas our
distortion objective requires relative control even for these gaps.
More recently, Li, Yu, and Zheng~\cite{LiYuZheng2026} established
a contrasting result for independent Bernoulli subsets of a product
of two regular trees, which are almost surely not quasi-isometric.
These results show how the ambient geometry affects the existence
of geometrically controlled correspondences between independent
random configurations.

\paragraph{Minimum-distortion bijections in $\mathbb R$.}
The computational problem of finding a minimum-distortion bijection
between two prescribed finite metric spaces was studied by Kenyon,
Rabani, and Sinclair~\cite{KRS2009}. For point sets on the line,
they gave a polynomial-time algorithm that finds an \emph{optimal}
bijection when the minimum distortion is below a fixed constant, with a subsequent work by
Chandran et al.~\cite{ChandranEtAl2008} improving this threshold.
Nayyeri and Raichel~\cite{NR2015} subsequently developed exact and
approximation algorithms with running-time bounds in terms of the
optimum distortion and the \emph{spread} of the input, namely its
diameter after rescaling the minimum distance to be one.
On the hardness side, Hall and Papadimitriou~\cite{HP2005} showed
that, even for point sets on the line, approximating the minimum
distortion within factors growing polynomially with the number
of points is NP-hard in the regime of polynomially large optimum
distortion. We prove that when the two point sets are random, finding a mapping that achieves a near-optimal distortion can be done whp.

\paragraph{Random geometric matching.}
A related line of work studies matchings between random point sets
that keep matched points close to one another. Goel, Rai, and
Krishnamachari~\cite{GRK2005} studied bottleneck matching between
independent uniform samples, minimizing the largest distance
between matched points, and used it to establish sharp thresholds
for random geometric graphs. Holroyd, Pemantle, Peres, and
Schramm~\cite{HPPS2009} studied invariant perfect matchings between
independent Poisson processes, focusing on the tail of the distance
from a typical point to its partner. More recently, Elboim, Spinka,
and Yakir~\cite{ESY2026} studied invariant perfect matchings between
randomly perturbed lattices and the lattice, obtaining optimal
tail bounds under suitable assumptions on the perturbations.

\paragraph{Mappings to grids and lattices.}
The classical work of Burago and Kleiner~\cite{BuragoKleiner1998}
and McMullen~\cite{McMullen1998} shows that, in dimensions at least
two, even uniformly separated point sets with bounded holes need
not admit a bi-Lipschitz bijection to a lattice.
Therefore, obstructions to bounded distortion can arise even in the
absence of exceptionally close pairs.
For finite lattice subsets, Dymond, Kalu\v{z}a, and
Kopeck\'a~\cite{DKK2018} showed that the minimum Lipschitz constant
of a bijection onto a regular grid can be arbitrarily large.
Dymond~\cite{DymondGrid} subsequently obtained polylogarithmic
upper bounds for suitable random lattice subsets.
These finite grid-mapping problems concern one-sided Lipschitz
control and a prescribed regular target, while our problem
requires control in both directions between two random samples.

\subsection{Organization of the paper}
We begin with the probabilistic estimates and metric reductions used
throughout the paper in \Cref{sec:prelim}.
In \Cref{sec:proof-overview}, we explain the origin of the $2/3$ exponent
and outline the proofs of the lower and upper bounds.
Next, in \Cref{sec:permutation-model}, we relate displacement to distortion
in the equally spaced model (where we study permutations of $[n]$ and their distortion), yielding a short proof of the lower bound
in \Cref{sec:lower-bound}.
We present the upper-bound construction in \Cref{sec:upper-bound}.

\paragraph{Use of AI.}
The proof strategy of the main theorem was developed by the authors, yet we incorporated a simplification found by AI. ChatGPT-6 suggested using the exponential representation in \Cref{prop:exponential-model} to replace our previous balls and bins type argument which allowed to use the independence of the gaps in the exponential model, and to simplify the approach behind the proof of \Cref{prop:matching-patterns}. After this simplification, we realized that our previous $n^{2/3 +o(1)}$ upper bound, should yield an $n^{2/3}(\log n)^{O(1)}$ upper bound, when choosing the parameters more delicately. We subsequently used ChatGPT-6 to choose the parameters for this refinement. We used ChatGPT-6 extensively to assist with writing the paper and generating the figures. 

After completing the main results we conjectured that controlling
displacement through oscillations, as we did in the upper bound, is necessary. We asked the model to prove the existence of oscillating permutation minors in any permutation with distortion and displacement nearly equal. ChatGPT-6 supplied the proof (and defined woven permutations) appearing in Appendix \ref{sec:woven}. One component of the AI's proof, \Cref{lem:windows}, provides a shorter proof to \Cref{thm:shift-obstruction}. We decided to retain our original proof in the main body of the paper.

\paragraph{Acknowledgments}
We thank Uri Grupel for useful discussions.

\section{Preliminaries}
\label{sec:prelim}

Throughout, let $n\ge 2$. We denote an $n$-point set $Z\subseteq\mathbb{R}$, by
$Z=\{z_1<\cdots<z_n\}$.  Let $\mathcal{O}$ denote the
distribution of the random set $Z=\{U_1,\ldots,U_n\}$, where $U_1,\ldots,U_n$ are iid
uniform random variables on $[0,1]$.  These variables are distinct almost surely, so
$Z$ has size $n$ almost surely.
For $1\le i<n$, define $\Delta_i^Z=z_{i+1}-z_i$, and put
$\Delta_{\min}^Z=\min_{1\le i<n}\Delta_i^Z$.

The following proposition is a well-known exponential representation of uniform spacings (see
\cite[Sec.~4.1]{Pyke1965} or \cite[Ch.~V, Thm.~2.2]{Devroye1986}).

\begin{proposition}
\label{prop:exponential-spacings}
Let $E_0,\ldots,E_n \sim \textnormal{Exp}(1)$ be iid and let $S_{n+1}=\sum_{i=0}^nE_i$.  Then,
\[
 \bigl(z_1,\Delta_1^Z,\ldots,\Delta_{n-1}^Z,1-z_n\bigr)
 \stackrel{d}{=}
 \frac{1}{S_{n+1}}(E_0,E_1,\ldots,E_n).
\]
\end{proposition}

The following observation follows trivially from the definition of Lipschitz constant of a function.
\begin{observation}
\label{obs:similarity-invariance}
Let $M,N$ be finite metric spaces and let $f:M\to N$ be a bijection. If the
metrics of $M$ and $N$ are multiplied by constants $\lambda>0$ and $\mu>0$, respectively, then $\Lip(f)$ is multiplied by $\mu/\lambda$ and
$\Lip(f^{-1})$ by $\lambda/\mu$.  Consequently, 
\begin{itemize}
    \item $\dist(f)$ and $c(M,N)$
remain the same and
    \item for $M, N \subseteq \mathbb{R}$, translating the sets does not change $c(M, N)$ as well.
\end{itemize}
\end{observation}

This observation combined with \Cref{prop:exponential-spacings} allows us to reformulate our problem.

\begin{proposition}\label{prop:exponential-model}
    Let $X, Y\sim \mathcal{O}$ be iid. Let $E_1^X, \ldots, E_{n-1}^X$ and $E_1^Y, \ldots, E_{n-1}^Y$ be iid exponential random variables with mean 1. Let $\widetilde X=\{ \widetilde x_1<\cdots<\widetilde x_n\}$ be defined as follows: $\tilde x_1 = 0$ and, for all $i\in[n-1]$, let $\tilde x_{i+1} =\tilde x_i +E_i^X$. Similarly, define $\widetilde Y=\{ \widetilde y_1<\cdots<\widetilde y_n\}$. Then,
    $$c\bigl(X, Y\bigr) \stackrel{d}{=} c\bigl(\widetilde X, \widetilde Y\bigr).$$
\end{proposition}

In order to bound $c(X, Y)$ from above we need to present some bijection $f:X \to Y$ with small Lipschitz constants. The following claim shows that we only need bound the quantities $\frac{|f(x_2) - f(x_1)|}{x_2 - x_1}$ and $\frac{|f^{-1}(y_2) - f^{-1}(y_1)|}{y_2- y_1}$ for the nearby pairs of points $x_1, x_2 \in X$ and  $y_1, y_2 \in Y$.

\begin{claim}
\label{cl:adjacent-gap}
Let $X=\{x_1<\cdots<x_n\} \subseteq \mathbb{R}$ and $Y=\{y_1<\cdots<y_n\} \subseteq \mathbb{R}$ and let $f:X \to Y$ be a bijection. Then, 
\[ 
 \Lip(f)  =\max_{1\le i<n}
   \frac{|f(x_{i+1})-f(x_i)|}{x_{i+1}-x_i},\quad\quad
 \Lip(f^{-1}) =\max_{1\le i<n}
   \frac{|f^{-1}(y_{i+1})-f^{-1}(y_{i})|}{y_{i+1}-y_i}.
    \]
\end{claim}

\begin{proof}
Let
\[
 L=\max_{i<n}
  \frac{|f(x_{i+1})-f(x_i)|}{x_{i+1}-x_i}.
\]
For $1 \leq j<k \leq n$, by the triangle inequality we get
\begin{align*}
|f(x_k)-f(x_j)|
&\leq \sum_{i=j}^{k-1}|f(x_{i+1})-f(x_i)| \\
&\leq L\sum_{i=j}^{k-1}(x_{i+1}-x_i) = L(x_k-x_j).
\end{align*}
So, $\Lip(f) =  L$. The formula for $f^{-1}$ is proved similarly.
\end{proof}

\begin{corollary}
\label{cr:adjacent-gap}
    
    Let $X=\{x_1<\cdots<x_n\}$ and $Y=\{y_1<\cdots<y_n\}$. Let $f:X \to Y$ be a bijection and let
$\pi:[n]\to [n]$ be the permutation satisfying $f(x_i)=y_{\pi(i)}$. Then,
\[ 
 \Lip(f)  =\max_{1\le i<n}
   \frac{|y_{\pi(i+1)}-y_{\pi(i)}|}{x_{i+1}-x_i},\quad\quad
 \Lip(f^{-1}) =\max_{1\le i<n}
   \frac{|x_{\pi^{-1}(i+1)}-x_{\pi^{-1}(i)}|}{y_{i+1}-y_i}.
    \]
\end{corollary}

\subsection{Spacing estimates}
\label{subsec:rank-spacing}
\label{subsec:window-minima}

First, we will need the following simple bound that can be found in~\cite{Janson2018}.

\begin{claim}
    \label{cl:exp_sum_tail_bounds}
    Let $S_n$ denote the sum of $n$ iid exponential variables with mean 1. There exists a universal constant $c>0$ such that, for all large enough $n$, we have
    $$\mathbb{P}\biggl(n/2 \leq S_{n+1} \leq 2n\biggr) \geq 1 - \exp(-cn).$$
    
\end{claim}
The following claim bounds $x_j - x_i$ for all $1 \leq i < j \leq n$ satisfying $j \geq i+ C\log n$.

\begin{claim}
\label{cl:global-rank-spacing}
There exists $C > 0$ such that for $X \sim \mathcal{O}$ whp the following holds simultaneously for all $1\le i<j\le n$ with $j-i\ge C\log n$,
\[\frac{j-i}{Cn}\le x_j-x_i\le\frac{C(j-i)}{n}.\]
\end{claim}

\begin{proof}
Let $E_0,\ldots,E_n \sim \textnormal{Exp}(1)$ be iid and let $S_{n+1}=\sum_{i=0}^nE_i$. By \Cref{prop:exponential-spacings},
\[
 \bigl(x_1,\Delta_1^X,\ldots,\Delta_{n-1}^X,1-x_n\bigr)
 \stackrel{d}{=}
 \frac{1}{S_{n+1}}(E_0,E_1,\ldots,E_n).
\]
By \Cref{cl:exp_sum_tail_bounds}, we whp  have $n/2 \leq S_{n+1} \leq 2n$. Also, for arbitrary indices $1 \leq i < j\leq n$ satisfying $j \geq i+ C\log n$ it holds that
\begin{equation}
    x_{j} - x_i \stackrel{d}{=} \frac{1}{S_{n+1}}\sum_{k = i}^{j-1} E_k.
    \label{eq:x_j-x_i_exponential_formula}
\end{equation}
We first claim that uniformly over all pairs $1 \leq i < j\leq n$ satisfying $j \geq i+ C\log n$, we have 
$$\mathbb{P}\biggl(\frac{j-i}{C/2}\le \sum_{k = i}^{j-1} E_k\le(C/2)(j-i)\biggr) = 1 - o\left(n^{-2}\right).$$
Indeed, by using the constant $c>0$ of \Cref{cl:exp_sum_tail_bounds}, and setting $C > \max(4, 10c^{-1})$ we have for large enough $n$:
\begin{align*}
    \mathbb{P}\biggl(\frac{j-i}{C/2}\le \sum_{k = i}^{j-1} E_k\le(C/2)(j-i)\biggr) &\geq 1 - \exp\bigl(-c \cdot (j-i)\bigr) \\
    &\geq  1 - \exp\bigl(-c \cdot C \cdot \log n\bigr) = 1 - o(n^{-2}).
\end{align*}
bounding the complement of the above event, over the $O(n^2)$ pairs of indices $(i,j)$ and combining that with $n/2 \leq S_{n+1} \leq 2n$, which holds whp, and with~\eqref{eq:x_j-x_i_exponential_formula} we conclude the proof.
\end{proof}

Notice that for all $i \in [n-1]$ the quantity $x_{i+1}-x_i$ does not exceed the quantity $x_k - x_j$ for some $j \leq i < k$ satisfying $k - j = \lceil C \log n\rceil$. Therefore, \Cref{cl:global-rank-spacing} gives the following corollary.

\begin{corollary}
\label{cr:maximum-gap}
Let $C > 0$ be large. For $X \sim \mathcal{O}$ and for all $i\in[n-1]$, whp
\[\Delta_i^X\leq C\frac{\log n}{n}.\]
\end{corollary}

\begin{claim}
\label{cl:window-minimum-gap}
Let $I\subseteq\{1,\ldots,n-1\}$ be an interval of size $m \in [n-1]$, and let $X \sim \mathcal{O}$. Then
\[\mathbb{P}\left(\frac{1}{2mn\log n}\le\min_{i\in I}\Delta_i^X\le\frac{2\log n}{mn}\right)=1-O\left(\frac{1}{\log n}\right).
\]  
\end{claim}

\begin{proof}
Let $E_0,\ldots,E_n \sim \textnormal{Exp}(1)$ be iid and let $S_{n+1}=\sum_{i=0}^nE_i$. By \Cref{prop:exponential-spacings},
\[
 \bigl(x_1,\Delta_1^X,\ldots,\Delta_{n-1}^X,1-x_n\bigr)
 \stackrel{d}{=}
 \frac{1}{S_{n+1}}(E_0,E_1,\ldots,E_n).
\]

Recall that a minimum of $m$ iid exponential variables with mean 1 is an exponential variable with mean $1/m$ and thus $\min_{i\in I}E_i\sim\operatorname{Exp}(m)$. Hence,
\[
 \PP\left(\min_{i\in I}E_i<\frac{1}{m\log n}\right)
 \le\frac{1}{\log n},
 \qquad
 \PP\left(\min_{i\in I}E_i>\frac{\log n}{m}\right)=\frac1n.
\]
Since $\min_{i\in I}\Delta_i^X\stackrel{d}{=}\min_{i\in I}E_i/S_{n+1}$, combining the bounds above with \Cref{cl:exp_sum_tail_bounds} we conclude the desired estimate.
\end{proof}

\begin{corollary}
\label{cr:minimum-gap}
For $X \sim \mathcal{O}$, whp 
\[\frac{1}{3n^2\log n}\le\Delta_{\min}^X\le\frac{3\log n}{n^2}.
\]
\end{corollary}

\section{Technical Overview}
\label{sec:proof-overview}

Let $X\sim\mathcal O$ and $Y\sim\mathcal O$ be independent samples, and
write $X=\{x_1<\cdots<x_n\}$ and $Y=\{y_1<\cdots<y_n\}$.
Observe that bijections $f:X\to Y$ are in one-to-one correspondence with permutations
$\pi:[n]\to[n]$, by setting $f(x_i)=y_{\pi(i)}$.
We use this correspondence to analyze how well $f$ maps points of $X$
to points of $Y$.
Throughout this overview, $\widetilde O$ and $\widetilde\Omega$
suppress polylogarithmic factors.

\subsection{The lower bound}
\label{subsec:overview-lower-bound}

For the lower bound, we classify the bijections $f:X\to Y$ into two types: local or nonlocal, according to
how the permutation $\pi$ that induces $f$ maps elements, and a parameter
$T\in\{2,\ldots,n\}$ to be chosen later. Informally, a $T$-local bijection maps
every index to a $T$-neighborhood of itself, that is $|\pi(i)-i|\leq T$ for the permutation $\pi$ associated with $f$.  As the interval $[0,1]$ has reflection as an isometry, we also need to take reflections into account in the definitions.

More formally, we measure the locality of a permutation by its $\ell_\infty$ distance
from a global isometry of $[n]$. The only such isometries are the identity
and the reversal $\rho(i)=n+1-i$. We therefore define\footnote{Displacement
is usually measured only with respect to the identity.}
\[
 \operatorname{disp}(\pi)=\min_{\tau\in\{\mathrm{id},\rho\}}
              \max_{i\in[n]}|\pi(i)-\tau(i)|.
\]

\paragraph[Permutations with displacement at most T]{Permutations with $\operatorname{disp}(\pi)\le T$.}
We show that, with high probability, \emph{all} permutations $\pi$ with
$\operatorname{disp}(\pi)\le T$ induce mappings $f$ satisfying
$\operatorname{Lip}(f)=\widetilde\Omega(n/T)$.
To see this, choose $x_i,x_{i+1}$ achieving the minimum gap in $X$.
By \Cref{cr:minimum-gap}, this gap is $\widetilde O(n^{-2})$ with high
probability.

Fix a choice of $X$, and let $i$ be the index of its minimum gap. Since $\operatorname{disp}(\pi)\le T$, the points $x_i,x_{i+1}$ are mapped to $y_j,y_k$, where
either $ j,k \in [i-2T,i+2T]\cap[n]$, or $j,k \in [\rho(i)-2T,\rho(i)+2T]\cap[n]$.
Since $Y$ is chosen independently of $X$, we can apply \Cref{cl:window-minimum-gap}
to the windows $[i-2T,i+2T]\cap[n]$, and $[\rho(i)-2T,\rho(i)+2T]\cap[n]$. As each window contains $O(T)$ consecutive points, whp the minimum distance in each window is
$\widetilde\Omega(1/(nT))$.
Consequently, whp  $|f(x_{i})-f(x_{i+1})|= \widetilde\Omega(1/(nT))$\footnote{Observe there is no need to union bound over all permutations satisfying $\operatorname{disp}(\pi)\le T$, as the argument uses the same two windows for all of them.}. Comparing with the original gap gives
\[
 \operatorname{Lip}(f)\ge
 \frac{|f(x_{i+1})-f(x_i)|}{x_{i+1}-x_i}
 =\widetilde\Omega(n/T).
\]

Since $\operatorname{disp}(\pi^{-1})=\operatorname{disp}(\pi)$, applying the same argument with
$X$ and $Y$ interchanged gives $\operatorname{Lip}(f^{-1})=\widetilde\Omega(n/T)$.
Thus, whp,
\emph{all} permutations $\pi$ satisfying $\operatorname{disp}(\pi)\le T$ give rise to
mappings $f$ with $\operatorname{dist}(f)=\widetilde\Omega((n/T)^2)$.

\paragraph[Permutations with displacement at least T]{Permutations with $\operatorname{disp}(\pi)\ge T$.}
For permutations of large displacement, we use a different obstruction. Our main ingredient is the following theorem, showing that displacement itself causes distortion. Concretely, let $\operatorname{dist}(\pi)$ be the distortion of a permutation $\pi$, when $[n]$ is considered with the induced metric from $\mathbb{R}$. We prove:

\begin{theorem}
\label{thm:shift-obstruction-overview}
There exists a universal constant $c>0$ such that, for every $n\ge2$ and every
permutation $\pi:[n]\to[n]$,
\[
 \operatorname{dist}(\pi)\ge c\,\operatorname{disp}(\pi).
\]
\end{theorem}

The above theorem can be seen as a finite (and quantitative) analog of a theorem by Benjamini and Shamov \cite{Benjamini-shamov} which says that bi-Lipschitz bijections of $\mathbb{Z}$ are at finite displacement from an isometry. Using the above theorem, the lower bound for the random setting of bijections $f:X\to Y$ follows naturally from this theorem. Concretely, in \Cref{cl:lipschitz-bound} we prove that, with high probability,
$\operatorname{Lip}(f)\ge c'\operatorname{Lip}(\pi)/\log n$ and
$\operatorname{Lip}(f^{-1})\ge c'\operatorname{Lip}(\pi^{-1})/\log n$ simultaneously for every
bijection $f$, where $c'>0$ is an absolute constant.
Together with \Cref{thm:shift-obstruction-overview}, this gives that whenever $\operatorname{disp}(\pi)\ge T$,
\[
 \operatorname{dist}(f)=\widetilde\Omega(\operatorname{dist}(\pi))
 =\widetilde\Omega(\operatorname{disp}(\pi))
 =\widetilde\Omega(T).
\]

\paragraph[The 2/3 exponent]{The $2/3$ exponent.}
As explained, the two lower bounds above hold whp for all permutations $\pi$. On this event, every bijection $f:X\to Y$ satisfies
\[
 \operatorname{dist}(f)=\widetilde\Omega\!\left(
   \min\left\{\frac{n^2}{T^2},T\right\}\right).
\]
The exponent $2/3$ comes from balancing the terms, which happens for $T \approx n^{2/3}$. This proves that, with high probability,
$c(X,Y)=\widetilde\Omega(n^{2/3})$ as required.

\subsection{The upper bound}
\label{subsec:overview-upper-bound}

We construct a bijection $f:X\to Y$ satisfying
$\max\{\operatorname{Lip}(f),\operatorname{Lip}(f^{-1})\}=\widetilde O(n^{1/3})$.
Suppose first that every pair of consecutive points in each sample
were at distance at least $n^{-4/3}$. Since, with high probability,
every consecutive gap in both $X$ and $Y$ is $\widetilde O(1/n)$, the increasing matching
$f(x_i)=y_i$ and its inverse would have expansion
$\widetilde O(n^{1/3})$, as required, by \Cref{cl:adjacent-gap}.
Thus the need to depart from the increasing bijection comes from
exceptionally close pairs of consecutive points below the
$n^{-4/3}$ threshold. We call these pairs \emph{tiny}; the precise
cutoff is slightly below $n^{-4/3}$ by a polylogarithmic factor.

\subsubsection{Handling the tiny pairs}
\label{subsubsec:overview-tiny-pairs}

We say that a pair $y_j,y_{j+1}$ \emph{covers} a pair $x_i,x_{i+1}$ if
\[
 \frac{y_{j+1}-y_j}{x_{i+1}-x_i},\quad
 \frac{x_{i+1}-x_i}{y_{j+1}-y_j}\in\widetilde O(n^{1/3}).
\]
As explained above, when all gaps are at least $n^{-4/3}$, the
increasing matching sends every consecutive pair to a pair covering it.

The lower-bound argument shows that a
bijection of distortion $\widetilde O(n^{2/3})$ must have displacement
$\widetilde O(n^{2/3})$, measured relative to the identity or the
reversal. We seek a mapping close to the identity and try to map
each tiny pair $x_i,x_{i+1}$ to a covering pair $y_j,y_{j+1}$ with
$|j-i|=\widetilde O(n^{2/3})$. Similarly, we require the mapping to
control the tiny pairs in $Y$.

Observe that this approach does not work as is. Suppose that
$x_i,x_{i+1}$ and $x_{i+d},x_{i+d+1}$ are tiny, where
$d\ll n^{1/3}$. Choosing their covers independently within
$n^{2/3}$ index difference (for $i$ and $i+d$) may send them to pairs whose indices differ
by $\Theta(n^{2/3})$. The two pairs of images $f(x_i),f(x_{i+1})$ and
$f(x_{i+d}),f(x_{i+d+1})$ will then be at distance of order
$n^{-1/3}$, while $x_i$ and $x_{i+d}$ are at distance at most
$\widetilde O(d/n)$. Consequently, up to polylogs, the
expansion can be as large as
\begin{equation}
 \frac{n^{-1/3}}{d/n}=\frac{n^{2/3}}{d}.
 \label{eq:overview-nearby-tiny-expansion}
\end{equation}
This exceeds our desired bounds on
$\operatorname{Lip}(f)$ and $\operatorname{Lip}(f^{-1})$.

Therefore, tiny pairs whose indices are within
$\widetilde O(n^{1/3})$ of one another require extra care and need to be
handled together. Concretely, we will force nearby tiny pairs (from
either copy) to be mapped together using \emph{the same shift}.
Notice that a priori it is not clear that all these choices, mapping
tiny pairs to covers, especially with the constraint that nearby tiny pairs
are mapped together, can be made consistently. We show, however, that
such a choice is possible, which is the main probabilistic argument
in the upper bound (the rest of the argument is geometric in nature,
and uses only basic properties of $\mathcal{O}$).

To handle both samples simultaneously, we call an index $i$ tiny if
either $x_i,x_{i+1}$ or $y_i,y_{i+1}$ is a tiny pair. Our main probabilistic argument shows that we can choose the partners
consistently and obtain the following interval matches.

\begin{proposition}[\Cref{prop:matching-patterns}, informal]
\label{prop:overview-matching}
For suitable $H=n^{2/3}\log^{O(1)}n$ and $D=n^{1/3}\log^{O(1)}n$, with high probability there
exist a partition $[n]=J_1\sqcup\cdots\sqcup J_\ell$
into consecutive parts, a family
\[
 \mathcal I=\{I_1,\bar I_1,\ldots,I_m,\bar I_m\}
\]
of pairwise disjoint intervals of indices, and a bijection $g:X_U\to Y_U$,
where $U=\bigcup_{I\in\mathcal I}I$ and $Z_A=\{z_i:i\in A\}$ for
$Z\in\{X,Y\}$, satisfying:
\begin{enumerate}
\item Each part satisfies $H\le |J_j|\le 3H$, and
each pair satisfies $|I_k|=|\bar I_k|=O(D\log n)$.
\item For every tiny index $i$, there is an interval
$I\in\mathcal I$ containing both $i$ and $i+1$.
\item For every $k$, the map $g$ sends $X_{I_k}$
increasingly onto $Y_{\bar I_k}$ and $X_{\bar I_k}$ increasingly onto
$Y_{I_k}$. In each mapping, the first and last points are mapped to
the first and last points, respectively.
\item For every $I\in\mathcal I$, $\max\{\operatorname{Lip}(g|_{X_I}),\operatorname{Lip}(g^{-1}|_{Y_I})\}=\widetilde O(n^{1/3})$.
\item Each pair $I_k,\bar I_k$
is contained in one part $J_j$, and both intervals are at distance
at least $D/2$ from its endpoints. Any two distinct intervals in
$\mathcal I$ are separated by at least $D/2$.
\end{enumerate}
\end{proposition}

We omit from the overview a sketch of the proof and highlight what the proposition guarantees.
By (2), every tiny pair in either sample has both endpoint indices in
one of the intervals of $\mathcal I$. Property (3) pairs these intervals
and matches their points increasingly, in both directions between $X$
and $Y$. Property (4) gives the desired bound $\widetilde O(n^{1/3})$
on the expansion in both directions of $g$ when restricted to an interval.

Properties (1) and (5) impose additional size and separation requirements
whose merit has not yet been motivated. They organize the construction
at two scales. The parts $J_1,\ldots,J_\ell$ have size $\Theta(H)$,
where $H=n^{2/3}\log^{O(1)}n$, and each matched pair of intervals lies
within a single part. The intervals in $\mathcal I$ are separated from
one another and from the endpoints of their containing parts by at
least $D/2$ indices. These requirements
provide the additional geometric conditions needed for the rest of
the proof, as we explain below.

\subsubsection{An interlude: the role of oscillations}
\label{subsubsec:overview-oscillation-interlude}
The partial matching $g$ from the previous proposition was chosen with
the aim of extending it to a bijection from $X$ to $Y$ with the desired
Lipschitz bounds. A priori, however, it is unclear whether this can be
done, and if so, how. We first illustrate this challenge through a toy
example involving a permutation of equally spaced points. The example
also shows that the bound in \Cref{thm:shift-obstruction} is sharp up to
constant factors in some cases.

Imagine completing the mapping from $X$ to $Y$ from left to right, mapping the points of $X$ whose indices lie in $J_j$ (whose length is $\approx n^{2/3}$) to the points of $Y$ with indices in the same part. For each $J_j$ we will handle the displacements prescribed by tiny pairs inside it as we encounter them.

For our toy model, let $m=k^2$, where $k\ge3$ is odd. We consider permutations of
$[m]$, with $m\approx n^{2/3}$ and $k\approx n^{1/3}$. Ignoring
rounding and logarithmic factors, $[m]$ represents the points of some $J_j$. We also require $\pi(1)=1$ to model that the left boundary of $J_j$ is mapped onto itself (together with fixing the right boundary of $J_{j-1}$, this controls the expansion across the boundary between $J_{j-1}$ and $J_j$).

Now suppose that we also require $\pi(k)=m-k+1$. Think of this as
modelling the requirement to map a tiny pair at index $k$ to a covering
pair at index $m-k+1$. This forces $k$ to move by
$m-2k+1=(k-1)^2\approx k^2$ indices. Can we find a permutation
$\pi:[m]\to[m]$ satisfying both $\pi(1)=1$ and $\pi(k)=m-k+1$,
while keeping $\operatorname{Lip}(\pi)$ and
$\operatorname{Lip}(\pi^{-1})$ bounded by $O(k)$?

We now describe such a permutation. Moreover, the
Lipschitz bounds in both directions strongly constrain how we can
accommodate this displacement.

\begin{figure}[H]
\centering
\begin{tikzpicture}[x=.60cm,y=1cm,font=\small]
\draw[gray!25] (0,6.3) -- (24,6.3);
\draw[gray!25] (0,0) -- (24,0);
\node[anchor=east] at (-.6,6.3) {$r$};
\node[anchor=east] at (-.6,0) {$\pi(r)$};
\draw[teal!70!black,line width=.6pt,opacity=.55,-{Latex[length=2.2mm,width=1.4mm]},shorten <=3.5pt,shorten >=3.5pt] (10,6.3) -- (2,0);
\draw[teal!70!black,line width=.6pt,opacity=.55,-{Latex[length=2.2mm,width=1.4mm]},shorten <=3.5pt,shorten >=3.5pt] (11,6.3) -- (7,0);
\draw[teal!70!black,line width=.6pt,opacity=.55,-{Latex[length=2.2mm,width=1.4mm]},shorten <=3.5pt,shorten >=3.5pt] (12,6.3) -- (12,0);
\draw[teal!70!black,line width=.6pt,opacity=.55,-{Latex[length=2.2mm,width=1.4mm]},shorten <=3.5pt,shorten >=3.5pt] (13,6.3) -- (17,0);
\draw[teal!70!black,line width=.6pt,opacity=.55,-{Latex[length=2.2mm,width=1.4mm]},shorten <=3.5pt,shorten >=3.5pt] (14,6.3) -- (22,0);
\draw[teal!70!black,line width=.6pt,opacity=.55,-{Latex[length=2.2mm,width=1.4mm]},shorten <=3.5pt,shorten >=3.5pt] (15,6.3) -- (23,0);
\draw[teal!70!black,line width=.6pt,opacity=.55,-{Latex[length=2.2mm,width=1.4mm]},shorten <=3.5pt,shorten >=3.5pt] (16,6.3) -- (16,0);
\draw[teal!70!black,line width=.6pt,opacity=.55,-{Latex[length=2.2mm,width=1.4mm]},shorten <=3.5pt,shorten >=3.5pt] (17,6.3) -- (13,0);
\draw[teal!70!black,line width=.6pt,opacity=.55,-{Latex[length=2.2mm,width=1.4mm]},shorten <=3.5pt,shorten >=3.5pt] (18,6.3) -- (6,0);
\draw[teal!70!black,line width=.6pt,opacity=.55,-{Latex[length=2.2mm,width=1.4mm]},shorten <=3.5pt,shorten >=3.5pt] (19,6.3) -- (3,0);
\draw[teal!70!black,line width=.6pt,opacity=.55,-{Latex[length=2.2mm,width=1.4mm]},shorten <=3.5pt,shorten >=3.5pt] (20,6.3) -- (4,0);
\draw[teal!70!black,line width=.6pt,opacity=.55,-{Latex[length=2.2mm,width=1.4mm]},shorten <=3.5pt,shorten >=3.5pt] (21,6.3) -- (5,0);
\draw[teal!70!black,line width=.6pt,opacity=.55,-{Latex[length=2.2mm,width=1.4mm]},shorten <=3.5pt,shorten >=3.5pt] (22,6.3) -- (14,0);
\draw[teal!70!black,line width=.6pt,opacity=.55,-{Latex[length=2.2mm,width=1.4mm]},shorten <=3.5pt,shorten >=3.5pt] (23,6.3) -- (15,0);
\draw[teal!70!black,line width=.6pt,opacity=.55,-{Latex[length=2.2mm,width=1.4mm]},shorten <=3.5pt,shorten >=3.5pt] (24,6.3) -- (24,0);
\draw[orange!90!black,line width=.8pt,opacity=.95,-{Latex[length=2.2mm,width=1.4mm]},shorten <=3.5pt,shorten >=3.5pt] (5,6.3) -- (21,0);
\draw[orange!90!black,line width=.8pt,opacity=.95,-{Latex[length=2.2mm,width=1.4mm]},shorten <=3.5pt,shorten >=3.5pt] (6,6.3) -- (18,0);
\draw[orange!90!black,line width=.8pt,opacity=.95,-{Latex[length=2.2mm,width=1.4mm]},shorten <=3.5pt,shorten >=3.5pt] (7,6.3) -- (11,0);
\draw[orange!90!black,line width=.8pt,opacity=.95,-{Latex[length=2.2mm,width=1.4mm]},shorten <=3.5pt,shorten >=3.5pt] (8,6.3) -- (8,0);
\draw[orange!90!black,line width=.8pt,opacity=.95,-{Latex[length=2.2mm,width=1.4mm]},shorten <=3.5pt,shorten >=3.5pt] (9,6.3) -- (1,0);
\draw[blue!75!black,line width=.9pt,opacity=1,-{Latex[length=2.2mm,width=1.4mm]},shorten <=3.5pt,shorten >=3.5pt] (1,6.3) -- (9,0);
\draw[blue!75!black,line width=.9pt,opacity=1,-{Latex[length=2.2mm,width=1.4mm]},shorten <=3.5pt,shorten >=3.5pt] (2,6.3) -- (10,0);
\draw[blue!75!black,line width=.9pt,opacity=1,-{Latex[length=2.2mm,width=1.4mm]},shorten <=3.5pt,shorten >=3.5pt] (3,6.3) -- (19,0);
\draw[black,line width=1.1pt,opacity=1,-{Latex[length=2.2mm,width=1.4mm]},shorten <=3.5pt,shorten >=3.5pt] (0,6.3) -- (0,0);
\draw[black,line width=1.1pt,opacity=1,-{Latex[length=2.2mm,width=1.4mm]},shorten <=3.5pt,shorten >=3.5pt] (4,6.3) -- (20,0);
\fill[black] (0,0) circle (1.8pt);
\fill[black] (0,6.3) circle (1.8pt);
\node[above=4pt] at (0,6.3) {1};
\node[below=4pt] at (0,0) {1};
\fill[black] (1,0) circle (1.8pt);
\fill[black] (1,6.3) circle (1.8pt);
\node[above=4pt] at (1,6.3) {2};
\node[below=4pt] at (1,0) {2};
\fill[black] (2,0) circle (1.8pt);
\fill[black] (2,6.3) circle (1.8pt);
\node[above=4pt] at (2,6.3) {3};
\node[below=4pt] at (2,0) {3};
\fill[black] (3,0) circle (1.8pt);
\fill[black] (3,6.3) circle (1.8pt);
\node[above=4pt] at (3,6.3) {4};
\node[below=4pt] at (3,0) {4};
\fill[black] (4,0) circle (1.8pt);
\fill[black] (4,6.3) circle (1.8pt);
\node[above=4pt] at (4,6.3) {5};
\node[below=4pt] at (4,0) {5};
\fill[black] (5,0) circle (1.8pt);
\fill[black] (5,6.3) circle (1.8pt);
\node[above=4pt] at (5,6.3) {6};
\node[below=4pt] at (5,0) {6};
\fill[black] (6,0) circle (1.8pt);
\fill[black] (6,6.3) circle (1.8pt);
\node[above=4pt] at (6,6.3) {7};
\node[below=4pt] at (6,0) {7};
\fill[black] (7,0) circle (1.8pt);
\fill[black] (7,6.3) circle (1.8pt);
\node[above=4pt] at (7,6.3) {8};
\node[below=4pt] at (7,0) {8};
\fill[black] (8,0) circle (1.8pt);
\fill[black] (8,6.3) circle (1.8pt);
\node[above=4pt] at (8,6.3) {9};
\node[below=4pt] at (8,0) {9};
\fill[black] (9,0) circle (1.8pt);
\fill[black] (9,6.3) circle (1.8pt);
\node[above=4pt] at (9,6.3) {10};
\node[below=4pt] at (9,0) {10};
\fill[black] (10,0) circle (1.8pt);
\fill[black] (10,6.3) circle (1.8pt);
\node[above=4pt] at (10,6.3) {11};
\node[below=4pt] at (10,0) {11};
\fill[black] (11,0) circle (1.8pt);
\fill[black] (11,6.3) circle (1.8pt);
\node[above=4pt] at (11,6.3) {12};
\node[below=4pt] at (11,0) {12};
\fill[black] (12,0) circle (1.8pt);
\fill[black] (12,6.3) circle (1.8pt);
\node[above=4pt] at (12,6.3) {13};
\node[below=4pt] at (12,0) {13};
\fill[black] (13,0) circle (1.8pt);
\fill[black] (13,6.3) circle (1.8pt);
\node[above=4pt] at (13,6.3) {14};
\node[below=4pt] at (13,0) {14};
\fill[black] (14,0) circle (1.8pt);
\fill[black] (14,6.3) circle (1.8pt);
\node[above=4pt] at (14,6.3) {15};
\node[below=4pt] at (14,0) {15};
\fill[black] (15,0) circle (1.8pt);
\fill[black] (15,6.3) circle (1.8pt);
\node[above=4pt] at (15,6.3) {16};
\node[below=4pt] at (15,0) {16};
\fill[black] (16,0) circle (1.8pt);
\fill[black] (16,6.3) circle (1.8pt);
\node[above=4pt] at (16,6.3) {17};
\node[below=4pt] at (16,0) {17};
\fill[black] (17,0) circle (1.8pt);
\fill[black] (17,6.3) circle (1.8pt);
\node[above=4pt] at (17,6.3) {18};
\node[below=4pt] at (17,0) {18};
\fill[black] (18,0) circle (1.8pt);
\fill[black] (18,6.3) circle (1.8pt);
\node[above=4pt] at (18,6.3) {19};
\node[below=4pt] at (18,0) {19};
\fill[black] (19,0) circle (1.8pt);
\fill[black] (19,6.3) circle (1.8pt);
\node[above=4pt] at (19,6.3) {20};
\node[below=4pt] at (19,0) {20};
\fill[black] (20,0) circle (1.8pt);
\fill[black] (20,6.3) circle (1.8pt);
\node[above=4pt] at (20,6.3) {21};
\node[below=4pt] at (20,0) {21};
\fill[black] (21,0) circle (1.8pt);
\fill[black] (21,6.3) circle (1.8pt);
\node[above=4pt] at (21,6.3) {22};
\node[below=4pt] at (21,0) {22};
\fill[black] (22,0) circle (1.8pt);
\fill[black] (22,6.3) circle (1.8pt);
\node[above=4pt] at (22,6.3) {23};
\node[below=4pt] at (22,0) {23};
\fill[black] (23,0) circle (1.8pt);
\fill[black] (23,6.3) circle (1.8pt);
\node[above=4pt] at (23,6.3) {24};
\node[below=4pt] at (23,0) {24};
\fill[black] (24,0) circle (1.8pt);
\fill[black] (24,6.3) circle (1.8pt);
\node[above=4pt] at (24,6.3) {25};
\node[below=4pt] at (24,0) {25};
\end{tikzpicture}
\caption{An oscillating permutation for $k=5$. (1) The black arrows represent the initial assumption on $\pi$, $\pi(1)=1$, and $\pi(k)=m-k+1$. (2) The blue arrows correspond to a walk from $1$ to $\pi(5)=21$ with low expansion. (3) The orange arrows are neighbours of already assigned images whose preimages under $\pi$ must lie near $[1,5]$. (4) The remaining arrows completing the permutation. }
\label{fig:overview-neighbor-progression}
\end{figure}
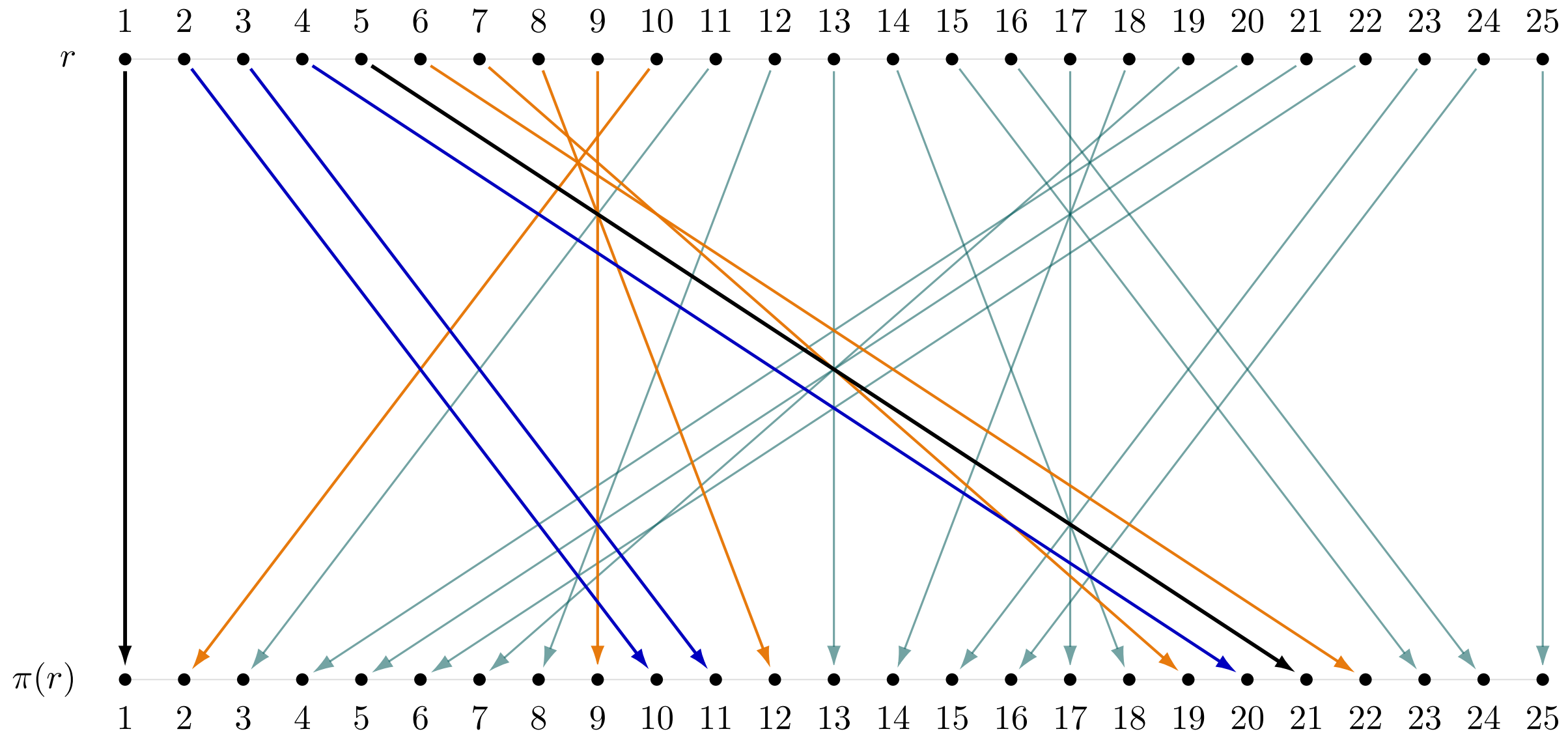

We will construct $\pi$ with
$\operatorname{Lip}(\pi),\operatorname{Lip}(\pi^{-1})\le 2k$,
giving distortion at most $4k^2=4m$. Since $\pi(1)=1$ and
$\pi(k)=m-k+1$, the first $k$ images must travel a distance of order
$k^2$. We choose these images so that they generally move to the right,
with consecutive images at most $2k$ indices apart, as shown by the
black and blue arrows in \Cref{fig:overview-neighbor-progression}.

Now consider the unused neighbors of these images. By the inverse
Lipschitz bound, their preimages must lie in $[k+1,3k]$. Thus the next
$2k$ images must revisit points spread across the interval we have
just traversed. For example, suppose we choose $\pi(2k)=2$. Then
$\pi(k)=m-k+1$ lies near the right endpoint, while $\pi(2k)=2$ lies
near the left endpoint. Thus, as indices increase from $k$ to $2k$,
the images must travel back across the interval, with consecutive
images at most $2k$ indices apart. We choose them so that they generally move
to the left, visiting the unused neighbors along the way, as shown
by the orange arrows.

The example illustrates how oscillations reconcile large displacement with small expansion in both directions. Reaching the prescribed image
requires moving across the interval, while controlling the inverse expansion requires returning to visit unused points near images already
assigned. Thus, a single prescribed displacement constrains the mapping
throughout the interval. The permutation described in this toy example
is shown in \Cref{fig:overview-neighbor-progression}. It achieves
expansion $O(k)$ in both directions and displacement $\Theta(k^2)$,
so its distortion is comparable to its displacement. The full
construction and its verification for general odd $k$ appear in
Appendix \ref{sec:overview-oscilating-verification}.

In the upper-bound proof, we face two additional difficulties not captured by the toy model:
\begin{enumerate}
\item \label{item:overview-varying-gaps}
The points left unmatched by the partial mapping $g$ are not equally spaced and have varying gap lengths, so we need to be able to handle variable gap lengths.
\item \label{item:overview-multiple-displacements}
The partial mapping $g$ may prescribe many displacements (intervals of $\mathcal I$), so we need
to allow for more flexible oscillations that can handle more than
a single displaced pair, in contrast to the toy example.
\end{enumerate}

In the next subsection, we describe how we overcome these difficulties
and conclude the overview.

\subsubsection{Implementing the oscillations}
\label{subsubsec:overview-implementing-oscillations}

We now explain how to adapt the above idea in the random setting.
Our goal is to extend the partial mapping $g$ on the points whose
indices lie in the intervals of $\mathcal I$ to a bijection from $X$
to $Y$ with small distortion. Let $p=C\lceil\log n\rceil$, where $C$
is a sufficiently large constant, and assume for simplicity that $p$
divides $n$. Partition $X$ and $Y$ into $M=n/p$ blocks of $p$
consecutive points each. We think of each block as a single object
and seek a permutation of the blocks, mapping the points increasingly
within each matched pair of blocks. All intervals in
\Cref{prop:overview-matching} can in fact be chosen as unions of whole blocks, including both the intervals in $\mathcal I$ \emph{and} the
parts $J_1,\ldots,J_\ell$. Consequently, the partial mapping $g$
already corresponds to a partial permutation of the blocks.
We call the remaining blocks \emph{free}.

Formally, the block mapping induced by a permutation
$\pi:[M]\to[M]$ sends
\[
 x_{(i-1)p+r}\longmapsto y_{(\pi(i)-1)p+r}
 \qquad(i\in[M],\ r\in[p]).
\]
As explained, we can ensure that the partial
mapping $g$ from \Cref{prop:overview-matching} is induced by a partial
permutation $\pi'$ on blocks, hence we seek an extension of $\pi'$ to
a permutation whose induced block mapping has small expansion in both directions.

For any block mapping extending $g$, consecutive points within the
same block already satisfy the desired expansion bound. Indeed, if
the block is not free, the block mapping agrees with $g$, so property (4)
of \Cref{prop:overview-matching} will bound the expansion by
$\widetilde O(n^{1/3})$. If the block is free, property (2) implies
that it contains no tiny gaps, so consecutive points are separated
by at least $\widetilde\Omega(n^{-4/3})$. Their images are consecutive
points in another free block, and hence are at distance at most
$\widetilde O(1/n)$. Therefore, such pairs also satisfy the desired expansion bound $\widetilde O(n^{1/3})$. Thus, in our choice of $\pi$ we only need to handle consecutive points lying in different blocks.

Unfortunately, oscillations cannot be applied directly in this setting. Ignoring polylogarithmic factors, a jump of approximately $n^{1/3}$ indices may separate the images of two consecutive points by approximately
$n^{-2/3}$. To allow such jumps at arbitrary boundaries between free blocks while keeping the expansion in each direction within our $n^{1/3}$ budget, we would need every such gap to be at least approximately $1/n$. However, gaps between adjacent free
blocks can be substantially smaller. This is how
difficulty~(\ref{item:overview-varying-gaps}) arises in the block formulation.

To resolve this, we apply a preparation step that enlarges the gaps between adjacent free blocks while preserving $g$. Formally, our rearrangement
lemma, \Cref{lem:rearrangement}, shows that, with high probability,
there exist configurations $X'=\{x'_1<\cdots<x'_n\}$ and
$Y'=\{y'_1<\cdots<y'_n\}$, and bijections $R_X:X\to X'$ and
$R_Y:Y\to Y'$, satisfying:
\begin{enumerate}
\item For every $i\in U$, $R_X(x_i)=x'_i=x_i$ and
$R_Y(y_i)=y'_i=y_i$. Thus the properties of $g$ guaranteed by
\Cref{prop:overview-matching} also hold for $X'$ and $Y'$.
\item Whenever blocks $i$ and $i+1$ are free, their boundary gap
satisfies $x'_{ip+1}-x'_{ip}=\widetilde\Omega(1/n)$, and similarly
in $Y'$.
\item The maps $R_X,R_Y$ and their inverses have Lipschitz constants
$\widetilde O(1)$.
\end{enumerate}
This preparation step resolves the difficulty~(\ref{item:overview-varying-gaps})
at only a polylogarithmic cost in distortion.

We now briefly explain how we address difficulty~(\ref{item:overview-multiple-displacements}). To overcome the difficulty, we use a more flexible oscillation construction,
\Cref{lem:oscillating-completion}, that extends
the prescribed partial permutation $\pi'$ to a permutation $\pi$ of
$[M]$ with expansion $\widetilde O(n^{1/3})$ in both directions.

Property (1) of \Cref{prop:overview-matching} bounds the size of each
part $J_j$ by $3H=\widetilde O(n^{2/3})$, limiting the distance the
oscillations should travel. It also bounds the length of each prescribed
interval by $\widetilde O(n^{1/3})$. This second bound
controls the size of a jump when an oscillation passes across a prescribed interval.

Property (5) ensures that every prescribed interval and its partner
lie in the same part $J_j$, allowing us to construct the extension
separately within each part. Finally, consecutive intervals of
$\mathcal I$ are separated by at least $D/2$ indices, and the same
margin separates them from the endpoints of their containing part.
This separation provides sufficiently long sequences of free blocks
in which to carry out the oscillations. The full lemma appears in
\Cref{sec:upper-bound} as \Cref{lem:oscillating-completion}.

\section{The Equally Spaced Model}
\label{sec:permutation-model}

Throughout this section, we equip $[n]$ with the metric $d(i,j)=|i-j|$.  We call this
the equally spaced model and refer to its elements as indices.  Under the identification
$i\leftrightarrow i/n$, a permutation $\pi$ corresponds to the bijection
$i/n\mapsto\pi(i)/n$ of the equally spaced points $\{i/n:i\in[n]\} \subseteq [0,1]$.  The distortion of the permutation (w.r.t. $d$) and the map on $\{i/n:i\in[n]\} \subseteq [0,1]$ associated with it, is the same by \Cref{obs:similarity-invariance}.  Thus, the model captures the geometry
at the typical spacing scale $1/n$ and above. There are three advantages to working in this model.
\begin{enumerate}
\item It is deterministic and isolates the permutation induced by a bijection $f:X\to Y$, where $X,Y\sim\mathcal O$, from irregularities in the consecutive gaps.

\item It makes the \emph{displacement} obstruction transparent: if $\pi(i)$ is far from
both $i$ and $n+1-i$ (in which case we say that $\pi$ displaces $i$), then $\pi$ must have large distortion.  This explains why matching a microscopic gap of $X$ to a distant microscopic gap of $Y$ typically introduces large distortion.

\item It isolates the role of oscillations, which allow large displacement while controlling distortion.  This mechanism is transparent when all adjacent gaps have the same length; its adaptation to random gaps requires a multiscale construction.
\end{enumerate}

The only isometries of this space are the identity $\operatorname{id}(i)=i$ and the
reversal $\rho(i)=n+1-i$.

\paragraph{Distortion, displacement.}
For a permutation $\pi\in S_n$ we define:
\[
 \Lip(\pi)=\max_{1\le i<n}|\pi(i+1)-\pi(i)|,
 \qquad
 \Lip(\pi^{-1})=\max_{1\le j<n}|\pi^{-1}(j+1)-\pi^{-1}(j)|.
\]
The \emph{distortion} of a permutation $\pi$ is defined by: \[\dist(\pi)=\Lip(\pi)\Lip(\pi^{-1}).\]
The \emph{displacement} of a permutation $\pi$ is defined by: \[
 \operatorname{disp}(\pi)
 :=
 \min_{\tau\in\{\operatorname{id},\rho\}}
 \max_{i\in[n]}|\pi(i)-\tau(i)|.
\]

The distortion measures how much $\pi$ can shrink or expand distances between indices. The displacement measures (in $\ell^{\infty}$) the distance of $\pi$ to some global isometry. The following theorem is key to the proof of the lower bound \emph{in the random setting}.

\begin{theorem}
\label{thm:shift-obstruction}
There is an absolute constant $c>0$ such that every $\pi\in S_n$ satisfies
\[
 \dist(\pi)\ge c\,\operatorname{disp}(\pi).
\]
\end{theorem}

We first begin with the following claim.

\begin{claim}
\label{cl:far_from_isometry}
Let $\pi\in S_n$ and $r>0$. Suppose that $\Lip(\pi)<r$ and
$\operatorname{disp}(\pi)>20r$. Then there exists $i\in[n]$ such that
\begin{equation}
|\pi(i)-i|\ge r
\qquad\text{and}\qquad
|\pi(i)-\rho(i)|\ge r.
\label{eq:is_far_from_isometry}    
\end{equation}
\end{claim}

\begin{proof}
Suppose towards a contradiction that no index $i\in [n]$ satisfies~\eqref{eq:is_far_from_isometry}. Label each $i\in[n]$ by $\operatorname{id}$ if
$|\pi(i)-i|<r$, and by $\rho$ otherwise. Since no index satisfies~\eqref{eq:is_far_from_isometry}, every index $i$ has label $\tau$ then $|\pi(i)-\tau(i)|<r$.

Let us show that the assumptions force $r>1$ and $n>20r$. Indeed, the images of two consecutive points are distinct, so $\Lip(\pi)\ge1$ and thus $r > \Lip(\pi)\ge1$. Also, every two points
of $[n]$ are at distance at most $n-1$, so
$\operatorname{disp}(\pi)\le n-1$ and hence $n > \operatorname{disp}(\pi) > 20r$.

Denote $t=2r$ and divide $[n]$ into
\[
L=\left[1,\frac n2-t\right],\qquad
M=\left(\frac n2-t,\frac n2+t\right),\qquad
R=\left[\frac n2+t,n\right],
\]
where the intervals are understood as subsets of $[n]$.

Let us show that the label is constant on each of $L$ and $R$. Otherwise, there exist consecutive indices
$i$ and $i+1$ in either $L$ or $R$ that have different labels. Thus, one of
$\pi(i)$ and $\pi(i+1)$ is within distance $r$ of an identity image and the other
is within distance $r$ of a reversal image. In either case, among the identity image and the reversal image one lies in $L$ and the other lies in $R$. As the distance between $L$ and $R$ is at least $2t$ we have
\[
|\pi(i+1)-\pi(i)|\geq2t-2r=2r,
\]
which contradicts $\Lip(\pi)<r$.

Now, let us show that all indices of both $L$ and $R$ have the same label. Wlog, assume towards contradiction that $L$ is labeled by
$\operatorname{id}$ and $R$ by $\rho$. Then
$\pi(L\cup R)\subseteq\{j\in[n]:j<n/2\}$. Indeed, for $i\in L$ we have
$\pi(i)<i+r\le n/2-r$ and for $i \in R$ we have
$\pi(i)<\rho(i)+r\le n/2+1-r<n/2$. 
Thus, the size of the image $\pi(L \cup R)$ is less than $n/2$. As $\pi$ is bijection, we conclude $|L \cup R| < n/2$. On the other hand, from the definitions of $L$ and $R$ we have
\[
|L\cup R|\ge n-4r-1>\left\lceil\frac n2\right\rceil,
\]
which is a contradiction.

Let $\tau\in\{\operatorname{id},\rho\}$ be the common \emph{global} label of all indexes from $L\cup R$.
For every $i\in M$,
\[
|i-\rho(i)|
  \le2\left|i-\frac n2\right|+1
  <4r+1.
\]
For $i \in M$, let $\tau_i\in\{\operatorname{id},\rho\}$ denote the map given in the label of $i$. Then, for all $i \in M$, we have 
\begin{align*}
    |\pi(i) - \tau(i)| &\leq |\pi(i) - \tau_i(i)| + |\tau_i(i) - \tau(i)|\\
    &\leq |\pi(i) - \tau_i(i)| + |\operatorname{id}(i) - \rho(i)|\\
    &< r + (4r+1) \leq 5r+1. \\
\end{align*}
Also, from the definition of $\tau$, we have $|\pi(i) - \tau(i)| \leq r$ for all $i \in L \cup R$. Thus, for all $i$ we have
$|\pi(i)-\tau(i)|<6r$ and, hence, $\operatorname{disp}(\pi)<6r$ --- a contradiction with
$\operatorname{disp}(\pi)>20r$.
\end{proof}

\begin{claim}
\label{cl:being_far_from_isometry_forces_distrotion}
Let $\pi\in S_n$, let $i,j\in[n]$, and let $d>0$. If $\pi(i)=j$ and
\[
|j-i|\ge d
\qquad\text{and}\qquad
|j-\rho(i)|\ge d,
\]
then $\dist(\pi)\ge d/200$.
\end{claim}

\begin{figure}[t]
\centering
\definecolor{BYRred}{RGB}{210,63,61}
\definecolor{BYRyellow}{RGB}{244,196,48}
\definecolor{BYRgold}{RGB}{145,101,0}
\definecolor{BYRblue}{RGB}{57,105,180}
\resizebox{\linewidth}{!}{%
\begin{tikzpicture}[
    x=0.48cm,
    y=1cm,
    cell/.style={
        draw=black!45,
        line width=0.35pt,
        rounded corners=0.7pt
    }
]
    \def\barheight{0.42}

    \foreach \x in {0,...,7} {
        \filldraw[cell,fill=BYRred]
            ({\x+0.09},2.90) rectangle ({\x+0.91},{2.90+\barheight});
    }
    \foreach \x in {8,...,12} {
        \filldraw[cell,fill=BYRyellow]
            ({\x+0.09},2.90) rectangle ({\x+0.91},{2.90+\barheight});
    }
    \foreach \x in {13,...,23} {
        \filldraw[cell,fill=BYRblue]
            ({\x+0.09},2.90) rectangle ({\x+0.91},{2.90+\barheight});
    }

    \node[anchor=east] at (-0.6,3.11) {$\mathcal C_k$};
    \node[text=BYRred,font=\small\bfseries] at (4,3.58) {red};
    \node[text=BYRgold,font=\small\bfseries] at (10.5,3.58) {yellow};
    \node[text=BYRblue,font=\small\bfseries] at (18.5,3.58) {blue};

    \node[below=3pt,font=\small] at (0.5,2.90) {$1$};
    \node[below=3pt,font=\small] at (8.5,2.90) {$i$};
    \node[below=3pt,font=\small] at (12.5,2.90) {$i+k$};
    \node[below=3pt,font=\small] at (23.5,2.90) {$n$};
    \draw[<->,line width=0.6pt] (8.5,3.98) -- (12.5,3.98)
        node[midway,above=1pt] {$k$};

    \draw[->,line width=0.8pt] (18.5,2.62) -- (18.5,1.70)
        node[midway,right] {$\pi$};

    \foreach \x in {0,...,4} {
        \filldraw[cell,fill=BYRred]
            ({\x+0.09},1.00) rectangle ({\x+0.91},{1.00+\barheight});
    }
    \foreach \x/\c in {5/BYRyellow,6/BYRblue,7/BYRblue,8/BYRyellow,
        9/BYRred,10/BYRyellow,11/BYRred,12/BYRred,13/BYRyellow,
        14/BYRblue,15/BYRblue,16/BYRyellow} {
        \filldraw[cell,fill=\c]
            ({\x+0.09},1.00) rectangle ({\x+0.91},{1.00+\barheight});
    }
    \foreach \x in {17,...,23} {
        \filldraw[cell,fill=BYRblue]
            ({\x+0.09},1.00) rectangle ({\x+0.91},{1.00+\barheight});
    }

    \node[anchor=east] at (-0.6,1.21) {$\mathcal C_k^\pi$};
    \node[text=BYRred,font=\small\bfseries] at (2.5,1.66) {red};
    \node[text=BYRblue,font=\small\bfseries] at (20.5,1.66) {blue};

    \draw[dashed,black!75,line width=0.75pt,rounded corners=1pt]
        (5,0.87) rectangle (17,1.55);
    \node[below=3pt,font=\small] at (0.5,1.00) {$1$};
    \node[below=3pt,font=\small] at (5.5,1.00) {$a$};
    \node[below=3pt,font=\small] at (13.5,1.00) {$j$};
    \node[below=3pt,font=\small] at (16.5,1.00) {$b$};
    \node[below=3pt,font=\small] at (23.5,1.00) {$n$};

    \draw[<->,line width=0.6pt] (5.5,0.28) -- (16.5,0.28)
        node[midway,below=2pt] {$w=b-a$};
\end{tikzpicture}
}
\caption{The colorings $\mathcal C_k$ and $\mathcal C_k^\pi$, where
$k$ is minimal such that $\mathcal C_k^\pi$ has no collision. Set
$j=\pi(i)$, and let $[a,b]$ be the smallest interval containing all
yellow points of $\mathcal C_k^\pi$, with width $w=b-a$. Since
$j\in[a,b]$, the absence of collisions, together with the sizes of the
red and blue classes, forces $w=\Omega(d)$; Moreover, the yellow interval in $\mathcal C_k$ has diameter $k$, so $\Lip(\pi)\ge w/k$, while the collision in
$\mathcal C_{k-1}^\pi$ gives $\Lip(\pi^{-1})\ge k+1$. Hence
$\dist(\pi)=\Omega(d)$.}
\label{fig:coloring-transport}
\end{figure}
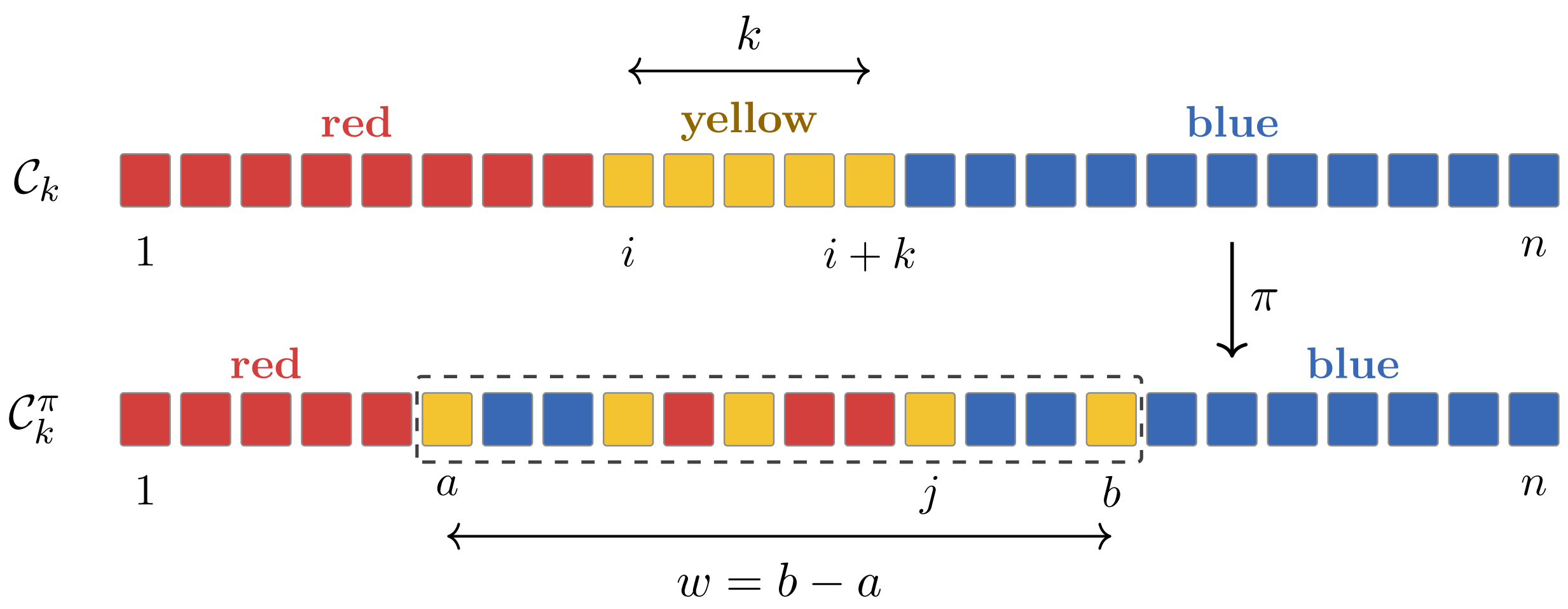

\begin{proof}
For $d\le200$, the claim follows from $\dist(\pi)\ge1$, so assume that
$d>200$. At least one of $i$ and $j$ lies at distance at least $d/4$ from
both endpoints of $[n]$: if $i$ is close to an endpoint, the two separation
assumptions force $j$ away from both endpoints. Replacing $\pi$ by
$\pi^{-1}$ and interchanging $i$ and $j$, if necessary, preserves the
assumptions and the distortion. We may therefore assume that
\[
i-1\ge\frac d4
\qquad\text{and}\qquad
n-i\ge\frac d4.
\]

Put $K=\lfloor d/100\rfloor$. For each integer $0\le k\le K$, let
$\mathcal C_k$ denote the coloring of the following three intervals
\[
[1,i-1],\qquad [i,i+k],\qquad [i+k+1,n]
\]
into red, yellow, and blue respectively. Transfer this coloring through $\pi$
by setting
\[
\mathcal C_k^\pi(x)=\mathcal C_k(\pi^{-1}(x)).
\]
We call two consecutive points colored red and blue, in either order, a
collision. A collision in $\mathcal C_k^\pi$ implies
$\Lip(\pi^{-1})\ge k+2$, since then there exist a red index $r$ and a blue index $b$ satisfying $|\pi(r) - \pi(b)| = 1$.

Let us show that the coloring $\mathcal C_0^\pi$ has a collision. Otherwise, $j = \pi(i)$ is the unique yellow point, and the intervals on its two sides are monochromatic
and have opposite colors. If the $[1, j-1]$ interval is red then counting the number of red points
gives $j-1=i-1$, and hence $j=i$. If it is blue then counting the number of blue points gives
$j-1=n-i$, and hence $j=\rho(i)$. Both conclusions contradict the
assumptions of the claim.

If $\mathcal C_K^\pi$ has a collision, then
$\Lip(\pi^{-1})\ge K+2>d/100$, and the claim follows. Otherwise, let
$1\le k\le K$ be the minimal index for which $\mathcal C_k^\pi$ has no
collision. Then $\mathcal C_{k-1}^\pi$ has a collision, and therefore
\[
\Lip(\pi^{-1})\ge k+1.
\]
We now show that $\Lip(\pi) \geq d/6k$. The reader is instructed to see \Cref{fig:coloring-transport} for illustration. Let $[a,b]$ be the smallest interval containing all yellow points of
$\mathcal C_k^\pi$, and put $w=b-a$. Since $j=\pi(i)$ is yellow, we have
$j\in[a,b]$. Let us show that $w>d/6$. Suppose otherwise. The red and blue
classes contain, respectively, $i-1$ and $n-i-k$ points, and both
quantities are larger than $d/5$. Since $[a,b]$ contains fewer than $d/5$
points, both red and blue colors of $\mathcal C_k^\pi$ occur outside of it. Since there are no collisions the two intervals outside of $[a,b]$ are monochromatic and have opposite
colors. If the $[1, a-1]$ interval is red, counting red points gives
$a-1\le i-1\le b$. Since $j\in[a,b]$, this yields
$|j-i|\le w+1<d$, a contradiction. If the $[1, a-1]$ interval is blue, then
$a-1\le n-i-k\le b$. Since
$\rho(i)=(n-i-k)+k+1$, we obtain
\[
|j-\rho(i)|\le w+k+1
   \le \frac d6+\frac d{100}+1<d,
\]
again a contradiction. Thus $w>d/6$.

The yellow interval in $\mathcal{C}_k$ has diameter $k$, so
$\Lip(\pi) \geq w/k$. Combining this with
$\Lip(\pi^{-1})\ge k+1$ gives
\[
\dist(\pi)
 \ge \frac{w}{k}(k+1)
 > w
 > \frac d6.
\]
In particular, $\dist(\pi)\ge d/200$.
\end{proof}

\begin{proof}[Proof of \Cref{thm:shift-obstruction}]
Put $D=\operatorname{disp}(\pi)$. If $D=0$, there is nothing to prove. If
$\Lip(\pi)\ge D/40$, then
\[
\dist(\pi)
 =\Lip(\pi)\Lip(\pi^{-1})
 \ge\frac D{40}.
\]

Suppose therefore that $\Lip(\pi)<D/40$, and put $r=D/40$. Since
$D>20r$, \Cref{cl:far_from_isometry} gives an index $i\in[n]$ satisfying
\[
|\pi(i)-i|\ge r
\qquad\text{and}\qquad
|\pi(i)-\rho(i)|\ge r.
\]
Applying \Cref{cl:being_far_from_isometry_forces_distrotion} with $j=\pi(i)$ and $d=r$ gives
\[
\dist(\pi)\ge\frac r{200}=\frac D{8000}.
\]
Thus the theorem holds with $c=1/8000$.
\end{proof}

\section{Proof of the Lower Bound}
\label{sec:lower-bound}

Set $T=\lfloor n^{2/3}\rfloor$. Let $X, Y \subseteq \mathbb{R}$ be sets of $n$ points. For each $Z\in\{X,Y\}$, let $\lambda_Z$
be an index satisfying $\Delta_{\lambda_Z}^Z=\Delta_{\min}^Z$, and
define
\[
I_Z^{\operatorname{id}}
=[\lambda_Z-T,\lambda_Z+T+1]\cap[n-1],
\qquad
I_Z^\rho
=[n-\lambda_Z-T,n-\lambda_Z+T+1]\cap[n-1].
\]

The following claim summarizes several useful properties that hold whp for both $X$ and $Y$.

\begin{claim}
\label{cl:lower_bound_properties_summarised}
    Let $X, Y \sim \mathcal{O}$ be iid. There exist universal constants $c, C>0$ such that, whp:
\begin{enumerate}
\item For each $Z=\{z_1<\cdots<z_n\}\in\{X,Y\}$,
\[
z_j-z_i\ge\frac{j-i}{Cn}
\quad\text{whenever }j-i\ge C\log n,
\qquad
\max_{i\in[n-1]}\Delta_i^Z\le\frac{C\log n}{n}.
\]

\item For each $\tau\in\{\operatorname{id},\rho\}$,
\begin{align*}
\Delta_{\min}^X,\Delta_{\min}^Y
    &\le\frac{C\log n}{n^2},\\
\min_{i\in I_X^\tau}\Delta_i^Y
    &\ge\frac{c}{nT\log n},
\qquad
\min_{i\in I_Y^\tau}\Delta_i^X
    \ge\frac{c}{nT\log n}.
\end{align*}
\end{enumerate}
\end{claim}

\begin{proof}
    The first property follows from \Cref{cl:global-rank-spacing,cr:maximum-gap}. For the second property we apply \Cref{cl:window-minimum-gap,cr:minimum-gap}, which says that in a window of length $T$, the minimum gap in $Y$ is $\geq \Omega(1/nT\log n)$ whp. The crucial point is that even though the interval of indices $I_X^{\tau}$ is random -- it is \emph{independent} of the randomness of $Y$, hence by first fixing $X$ we may apply \Cref{cl:window-minimum-gap,cr:minimum-gap} for $Y$, and vice versa. 
\end{proof}

Let $\mathcal{E}$ denote the event that both properties from \Cref{cl:lower_bound_properties_summarised} hold. We will show that whenever $\mathcal{E}$ holds, any mapping between $X$ and $Y$ has sufficiently large distortion. We begin with the following claim, which allows us to lower bound the Lipschitz constant of an arbitrary bijection $f:X\to Y$ by the Lipschitz constant of its associated permutation.

\begin{claim}
\label{cl:lipschitz-bound}
Let $X=\{x_1<\cdots<x_n\leq 1\}$ and $Y=\{0 \leq y_1<\cdots<y_n\leq 1\}$ be sets of points satisfying $\mathcal{E}$. Let $f:X\to Y$ be a
bijection and let $\pi\in S_n$ be the permutation satisfying $f(x_i)=y_{\pi(i)}$ for all $i \in [n]$. Then, uniformly over $f$,
\[
\Lip(f)\gtrsim\frac{\Lip(\pi)}{\log n}.
\]
\end{claim}

\begin{proof}
Put $r=\Lip(\pi)$, and choose $i\in[n-1]$ such that
$|\pi(i+1)-\pi(i)|=r$. If $r\ge C\log n$, then
\Cref{cr:adjacent-gap} and the first property from $\mathcal{E}$ give
\[
\Lip(f)
 \ge
\frac{|f(x_{i+1})-f(x_{i})|}
      {x_{i+1} - x_{i}} =\frac{|y_{\pi(i+1)}-y_{\pi(i)}|}{x_{i+1}-x_i}
 \ge\frac{r/(Cn)}{C\log n/n}
 =\frac{r}{C^2\log n}.
\]
If $r<C\log n$, then a constant lower bound on $\Lip(f)$ suffices. For all
sufficiently large $n$, the first property of $\mathcal{E}$ applied to $i=1$ and $j=n$
gives
\[
\operatorname{diam}(Y)=y_n-y_1\ge\frac{n-1}{Cn}\ge\frac{1}{2C}.
\]
Since $f$ is onto and $\operatorname{diam}(X)\le1$, it follows that
\[
\Lip(f)\ge\frac{\operatorname{diam}(Y)}{\operatorname{diam}(X)}
\ge\frac{1}{2C}
\ge\frac{r}{2C^2\log n}.
\]
\end{proof}

\begin{proof}[Proof of the lower bound in \Cref{thm:main}]
Due to \Cref{cl:lower_bound_properties_summarised}, we can assume that $X$ and $Y$ satisfy $\mathcal{E}$. Consider any bijection
$f:X\to Y$ and let $\pi\in S_n$ be the permutation satisfying $f(x_i)=y_{\pi(i)}$ for all $i \in [n]$. Denote
$D=\operatorname{disp}(\pi)$. Recall that $T=\lfloor n^{2/3}\rfloor$.

Suppose first that $D\le T$. Choose
$\tau\in\{\operatorname{id},\rho\}$ such that
$|\pi(i)-\tau(i)|\le T$ for every $i\in[n]$. We have 
\begin{align}
    \Lip(f) \ge
\frac{|f(x_{\lambda_X+1})-f(x_{\lambda_X})|}
      {x_{\lambda_X+1} - x_{\lambda_X}} = \frac{|y_{\pi(\lambda_X+1)}-y_{\pi(\lambda_X)}|}
      {x_{\lambda_X+1} - x_{\lambda_X}}=
 \frac{|y_{\pi(\lambda_X+1)}-y_{\pi(\lambda_X)}|}
      {\Delta_{\min}^X} \geq \frac{\min_{i\in I_X^\tau}\Delta_i^Y}{\Delta_{\min}^X}
      \label{eq:lower_bound_first_case}
\end{align}
where the last inequality holds, as both
$\pi(\lambda_X)$ and $\pi(\lambda_X+1)$ belong to $I_X^\tau$. Due to the second property of $\mathcal{E}$ we further have
\begin{equation*}
\eqref{eq:lower_bound_first_case} \gtrsim
 \frac{n}{T(\log n)^2}.
\end{equation*}

Observe that the permutation $\pi^{-1}$ is within distance $T$ of the same $\tau$,  hence $\operatorname{disp}(\pi^{-1})=\operatorname{disp}(\pi)$.
Therefore, we may apply the same argument to $f^{-1}:Y\to X$. This gives
$\Lip(f^{-1})\gtrsim n/(T(\log n)^2)$, and hence
\[
\dist(f)\gtrsim \frac{n^2}{T^2(\log n)^4}.
\]

Suppose now that $D>T$. Applying \Cref{cl:lipschitz-bound} to $f$ and
$f^{-1}$, and then applying \Cref{thm:shift-obstruction}, gives
\[
\dist(f)
 \gtrsim\frac{\dist(\pi)}{(\log n)^2}
 \gtrsim\frac{D}{(\log n)^2}
 >\frac{T}{(\log n)^2}.
\]

\noindent Combining the two cases and using $T=\lfloor n^{2/3}\rfloor$, we obtain
\[
\dist(f)
 \gtrsim
 \min\left\{
   \frac{n^2}{T^2(\log n)^4},
   \frac{T}{(\log n)^2}
 \right\}
 \gtrsim\frac{n^{2/3}}{(\log n)^4}.
\]
Since the two properties hold simultaneously with high probability and
$f$ was arbitrary, taking the minimum over all bijections proves the
theorem.
\end{proof}
\section{The upper bound}
\label{sec:upper-bound}

The main goal of this section is to prove the following theorem.

\begin{theorem}[Upper bound]\label{thm:upper-bound}
Let $X$ and $Y$ be independent samples of $n$ uniform points in $[0,1]$.
With high probability,
\[
 c(X,Y)\le n^{2/3}\log^{O(1)}n.
\]
\end{theorem}

\subsection{Set-up}

Assume throughout \Cref{sec:upper-bound} that $n$ is sufficiently large. Let
$X=\{x_1<\cdots<x_n\}$ and $Y=\{y_1<\cdots<y_n\}$ be independent
samples of $n$ uniform points in $[0,1]$.
For $Z=\{z_1<\cdots<z_n\}\in\{X,Y\}$ and $I\subseteq[n]$, write
$Z_I=\{z_i:i\in I\}$.
By Proposition~2.1, in order to show Theorem~6.1 we may assume that $X$ and $Y$ are defined as follows. Let $(E_i^X)_{i=0}^n$ and $(E_i^Y)_{i=0}^n$ be independent $\operatorname{Exp}(1)$ variables. Define the gaps of consecutive points of $X$ and $Y$ by: for $Z\in\{X,Y\}$, set $z_1=\frac{E_0^Z}{\Sigma_Z}$
\[
 \Delta_i^Z:=z_{i+1}-z_i=\frac{E_i^Z}{\Sigma_Z}
 \quad(i\in[n-1]),
 \qquad \Sigma_Z=\sum_{j=0}^n E_j^Z.
\]
By \Cref{cl:exp_sum_tail_bounds}, $n/2\le\Sigma_X,\Sigma_Y\le2n$ with high probability.

We now introduce the main parameters for the proof, and explain their role. Put
\[
 L=\lceil\log n\rceil,\qquad
 H=\lceil n^{2/3}L^4\rceil,\qquad p=C_0L,
\]
where $C_0$ is a sufficiently large even positive integer.

Between $X$ and $Y$, we will construct a specific kind of mapping,
called a \emph{block mapping}. The parameter $p=\Theta(\log n)$ is
the block size: we partition each sample into blocks of $p$ consecutive
points, with a possibly shorter last block. A block mapping sends each
full block of $X$ increasingly onto a full block of $Y$, keeping its
points together.

To construct the block mapping, we will also find a partition
$[n]=J_1\mathbin{\dot\cup}\cdots\mathbin{\dot\cup}J_m$ into intervals
of between $H$ and $3H$ consecutive indices. Thus each part has
$\Theta(n^{2/3}\log^4 n)$ indices. For each $J_i$, we will construct
a suitable matching between the corresponding blocks of $X_{J_i}$ and $Y_{J_i}$ separately.

We now introduce and explain a second set of parameters, used to design
the block mapping \emph{within each} $J_i$. Set
\[
 s=L^3,\qquad D=Ap\lceil\sqrt H\rceil,\qquad
 t=1+\lceil\log_sD\rceil
\]
where $A$ is a sufficiently large absolute integer.

Since we study distortion up to polylogarithmic factors, we record
exponential gaps $E_i^Z$ only up to a factor of $s$.
Concretely, for $x>0$, define
$e_s(x)=\max\{0,\lfloor\log_s(1/x)\rfloor\}$.
The integer $e_s(x)$ records the scale of $x$, with larger exponents
corresponding to smaller gaps. For the gaps of our samples, write
$e_i^Z=e_s(E_i^Z)$ for $Z\in\{X,Y\}$ and $i\in[n-1]$. Then
\[
 e_i^Z=a\quad\Longleftrightarrow\quad
 s^{-(a+1)}<E_i^Z\le s^{-a}\qquad(a\ge1).
\]
Thus each positive exponent $a$ represents a class of gaps lying
between two consecutive powers of $1/s$. Within this class, replacing
an exponential gap by the upper endpoint $s^{-a}$ changes its value
by at most a factor of $s$. In particular, two exponential gaps with
the same positive exponent differ by at most this factor.
All larger gaps, those with $E_i^Z>s^{-1}$, are grouped into the
remaining class $e_i^Z=0$.

\begin{samepage}
As explained in \Cref{sec:proof-overview}, consecutive points of $X$ or $Y$
whose gap is smaller than $n^{-4/3}$ require extra care.
We use $t$ to define the precise cutoff in terms of the gap exponents.

\begin{definition}[Tiny indices]\label{def:tiny-indices}
An index $i\in[n-1]$ is \emph{tiny} if
$\max\{e_i^X,e_i^Y\}\ge t$, or equivalently, if
$\min\{E_i^X,E_i^Y\}\le s^{-t}$.
\end{definition}

For $Z\in\{X,Y\}$, the condition $E_i^Z\le s^{-t}$ is equivalent to
$\Delta_i^Z\le1/(\Sigma_Zs^t)$.
On the event $n/2\le\Sigma_X,\Sigma_Y\le2n$, this distance threshold
is within a factor of $2$ of $1/(ns^t)$.
The bounds $1/(ns^2D)\le1/(ns^t)\le1/(nsD)$, together with
$s=L^3$ and $D=\Theta(n^{1/3}L^3)$, show that this threshold is
smaller than $n^{-4/3}$ by a polylogarithmic factor.
\par
\end{samepage}

We partition tiny indices into \emph{components}. Namely, we consider the partition of tiny indices into maximal subsets (which we call components) such that every pair of tiny indices $i< j$ with $j - i \leq D$ lies in the same component. The low-distortion map $f:X\to Y$ we aim to construct will satisfy the following property: For every component $C$ there exists an integer $k \in [-n, n]$ such that, for all tiny indices $i \in C$ the map $f:X\to Y$ satisfies $f(x_i) = y_{i+k}$ and $f(x_{i+1}) = y_{i+1+k}$. In other words, $f$ applies the same translation to all indices from $C$ when transferring them from $X$ to $Y$. Without this constraint, two nearby tiny gaps could be matched to distant gap in $J_i$, stretching the points between them too much. 
We aim to first build a partial mapping between $X$ and $Y$ that maps some intervals $\mathcal{I} = \{I_1, \ldots, I_{2\ell}\}$ to each other in a way that all components of tiny indices are contained in these intervals. During this step, we will also require the intervals from $\mathcal{I}$ to be separated by at least $D/2$ points. Then, the choice $D=Ap\lceil\sqrt H\rceil=\Theta(n^{1/3}L^3)$ will give us enough
space between the intervals to extend $f$ to the whole set $X$ while keeping its distortion small.

\subsection{Matching the tiny gaps}

\begin{definition}[Blocks]\label{def:blocks}
Set $M=\lceil n/p\rceil$. For $j\in[M]$, define
\[
 B_j=\{(j-1)p+1,\ldots,\min\{jp,n\}\},
 \qquad B_j^Z=Z_{B_j}\quad(Z\in\{X,Y\}).
\]
The blocks $B_1,\ldots,B_M$ partition $[n]$ into consecutive intervals,
each of size $p$ except possibly the last.
We call $B_M$ and $B_M^Z$ the \emph{terminal blocks}, even when they
have $p$ elements.
\end{definition}

\begin{definition}[Partial block mappings]\label{def:block-mappings}
Let $S,T\subseteq[M-1]$, and let $\sigma:S\to T$ be a bijection.
The \emph{partial block mapping induced by $\sigma$} is the map
$f:\bigcup_{j\in S}B_j^X\to\bigcup_{j\in T}B_j^Y$
that sends each $B_j^X$, $j\in S$, increasingly onto $B_{\sigma(j)}^Y$.
Explicitly,
\[
 f(x_{(j-1)p+a})=y_{(\sigma(j)-1)p+a}
 \qquad(j\in S,\ a\in[p]).
\]
Given a partial block mapping $f$ induced by $\sigma:S\to T$, we say
that $B_j^X$ is \emph{free} if $j\notin S$, and that $B_j^Y$ is
\emph{free} if $j\notin T$.
We call $\sigma$ a \emph{partial involution} if $S=T$ and
$\sigma^2=\mathrm{id}_S$.
\end{definition}

Given a permutation $\pi$ of $[M]$ extending $\sigma$ and satisfying
$\pi(M)=M$, define $F:X\to Y$ by mapping each $B_j^X$ increasingly
onto $B_{\pi(j)}^Y$. We call $F$ the \emph{block extension of $f$
induced by $\pi$}.

\begin{proposition}[Matching the tiny gaps]\label{prop:matching-patterns}
With high probability, there exist $\ell\in\mathbb N_0$, a family
$\mathcal I=\{I_1,\bar I_1,\ldots,I_\ell,\bar I_\ell\}$ of
$2\ell$ nonempty intervals in $[n]$, where each interval is a union of nonterminal blocks,
and a partial block mapping
$f:X_U\to Y_U$, where $U=\bigcup_{I \in \mathcal I}I$, with the following properties.
\begin{enumerate}
\item\label{it:ub-matching-involution} The mapping $f$ is induced by a partial involution and maps
$X_{I_i}$ in increasing order onto $Y_{\bar I_i}$ for every $i\in[\ell]$.
\item\label{it:ub-matching-tiny} For every tiny index $j$, the pair $\{j,j+1\}$ lies in some
$I\in\mathcal I$.
\item\label{it:ub-matching-lipschitz} For every $I\in\mathcal I$,
\[
 \max\bigl\{\Lip(f|_{X_I}),\Lip(f^{-1}|_{Y_I})\bigr\}
 =O(Ls^{t+1}).
\]
\item\label{it:ub-matching-separation} The intervals in $\mathcal I$ are pairwise separated by at least
$D/2$, and each has size $O(LD)$.
\item\label{it:ub-matching-margins} There is a partition $[n]=J_1\mathbin{\dot\cup}\cdots
\mathbin{\dot\cup}J_m$ into consecutive intervals with
$H\le|J_j|\le3H$ for every $j\in[m]$. For every $i\in[\ell]$, there is $j\in[m]$ such that both $I_i$ and $\bar I_i$ lie in $J_j$ and are at distance at least
$D/2$ from each endpoint of $J_j$.
\end{enumerate}
\end{proposition}

\subsection{Patterns and covers}

Recall that $e_i^Z$ denotes the discretized exponent of the $i$-th
exponential gap in $Z$. We introduce the definition of a rooted pattern. The rooted pattern is a formal tuple that we use to record some of the gap exponents $e_i^Z$ and some distances between them. Namely, the rooted pattern consists of the root $r \in [n]$, and the pattern, which is a tuple of distances $d_1, \ldots d_{k-1} \in [n]$ and of integers $\alpha_1^X,\alpha_1^Y, \ldots \alpha_k^Y$ for some $k > 0$. Roughly, the rooted pattern occurs if for all $i \in \{0, \ldots, k-1\}$ the discretized exponential gaps $e_{r+d_1+\ldots+d_i}^X$ and $e_{r+d_1+\ldots+d_i}^Y$ coincide with $\alpha_{i+1}^X$ and $\alpha_{i+1}^Y$. In that sense, a rooted pattern that occurs records several exponential gaps and the number of indices between them. Notice that the root of the pattern represents the index of the first gap of the pattern. Let us first define patterns without roots, then rooted patterns.

\begin{definition}[Patterns]\label{def:patterns}
A pattern of size $k\ge1$ is a sequence
\[
 P=(\alpha_1^X,\alpha_1^Y,d_1,\alpha_2^X,\alpha_2^Y,
 \ldots,d_{k-1},\alpha_k^X,\alpha_k^Y),
\]
where $\alpha_j^X,\alpha_j^Y\in\mathbb N_0$ are the prescribed gap
exponents at the $j$-th recorded position, and each $d_j$ is a positive
integer giving the difference between the $j$-th and $(j+1)$-st recorded indices.
We write $|P|=k$.
\end{definition}

\noindent
We now introduce some terminology for patterns that will be used later
in the proof. Let $P$ be a pattern of size $k$, and set $a_1=0$ and
$a_j=\sum_{h<j}d_h$ for $2\le j\le k$.

\begin{itemize}
\item The \emph{shape} of $P$ is $(d_1,\ldots,d_{k-1})$.
For $k=1$, the shape is empty.

\item The \emph{span} of $P$ is
$\operatorname{span}(P)=a_k=\sum_{j=1}^{k-1}d_j$.
For $k=1$, the span is zero.

\item A \emph{rooted pattern} is a pair
$\mathbf P=(r,P)$, where $r\in[n-1]$ and
$r+\operatorname{span}(P)\le n-1$. We call $r$ its root.

\item The \emph{support} of $\mathbf P=(r,P)$ is
$\operatorname{supp}(\mathbf P)=\{r+a_1, r+a_2,\ldots,r+a_k\}$,
the set of recorded gap indices.

\item The \emph{interval} of $\mathbf P=(r,P)$ is
$\operatorname{int}(r,P)=[r,r+\operatorname{span}(P)]\cap\mathbb Z$.
It contains all gap indices between the first and last recorded positions.

\item We say that $(r,P)$ \emph{occurs} if
$e_{r+a_j}^Z=\alpha_j^Z$ for every $j\in[k]$ and $Z\in\{X,Y\}$.
We say that $P$ occurs in an interval $I\subseteq[n]$ at root $r$ if
$(r,P)$ occurs and $[r,r+\operatorname{span}(P)+1]\cap\mathbb Z\subseteq I$.

\end{itemize}

\begin{definition}[Covers]\label{def:covers}
Let $P$ and $\bar P$ be patterns of size $k$ with the same shape,
and denote their exponent pairs by $(\alpha_j^X,\alpha_j^Y)$ and
$(\bar\alpha_j^X,\bar\alpha_j^Y)$, respectively.
We say that $\bar P$ \emph{covers} $P$ if
$\bar\alpha_j^X\le\alpha_j^Y$ and $\bar\alpha_j^Y\le\alpha_j^X$
for every $j\in[k]$.
The cover has \emph{low expansion} if
\[
 (\alpha_j^Y-t)_+\le\bar\alpha_j^X\le\alpha_j^Y,
 \qquad (\alpha_j^X-t)_+\le\bar\alpha_j^Y\le\alpha_j^X
 \qquad(j\in[k]),
\]
where $(a)_+=\max\{a,0\}$.
\end{definition}

Before we justify the name \emph{low-expansion} for covers, we record the following regularity bounds on the gap sizes $\Delta_i^Z=E_i^Z/\Sigma_Z$.

\begin{observation}\label{obs:gap-regularity}
With high probability, simultaneously for $Z\in\{X,Y\}$, we have
$n/2\le\Sigma_Z\le2n$ and $E_i^Z\le4L$ for every $0\le i\le n$.
In particular, every gap is at most $8L/n$.
\end{observation}

\begin{proof}
By \Cref{cl:exp_sum_tail_bounds}, $n/2\le\Sigma_Z\le2n$ for both samples $X$ and $Y$ with high
probability. Since $\mathbb P(E_i^Z>4L)=e^{-4L}$, a union bound over
the $2(n+1)$ exponentials gives failure probability
$2(n+1)e^{-4L}=o(1)$ for the remaining bound.
\end{proof}

Consider a rooted pattern $\mathbf P = (r, P)$ of size $k$. Recall, if $\mathbf P$ occurs, it records the gap exponents in $X$ and $Y$ at $k$ indices $i_1< \ldots <i_k$. In our strategy, we aim to match the indices $i_1, i_1+1, i_2, i_2+1, \ldots i_k+1$ in increasing order to the corresponding indices of some other rooted pattern $\mathbf{\bar P} = (\bar r, \bar P)$ with the same size and shape as $\mathbf P$. We aim to choose $\mathbf{\bar P}$ to be a low-expansion cover of $\mathbf P$ located near $\mathbf P$,  i.e. so that $|\bar r - r|$ is small. As a consequence, the natural mapping from $\mathbf P$ to $\mathbf{\bar P}$ will have a small Lipschitz constant due to \cref{obs:gap-regularity}. In addition, we will need to extend this natural mapping to some intervals $I \supseteq \operatorname{int}(\mathbf P) $ and $\bar I \supseteq \operatorname{int}(\mathbf{\bar P})$ in a way that keeps the Lipschitz constants small. We aim to show that the increasing extension also has small Lipschitz constant as long as all tiny gaps of $I$ and $\bar I$ are recorded by the rooted patterns $\mathbf{P}$ and $\mathbf{\bar P}$ and some other natural conditions hold.

\begin{claim}\label{obs:cover-lipschitz}
With high probability, the following holds.
Let $I,\bar I\subseteq[n]$ be intervals of consecutive point indices
with $|I|=|\bar I|$.
Suppose that $P$ occurs in $I$ at root $r$, and a low-expansion cover
$\bar P$ of $P$ occurs in $\bar I$ at root $\bar r$, where
\[
 r-\min I=\bar r-\min\bar I.
\]
Suppose that all tiny gap indices $i$ with $i,i+1\in I$ lie in
$\operatorname{supp}(r,P)$, and similarly for $\bar I$ and
$(\bar r,\bar P)$.
Then the increasing\footnotemark{} bijections $X_I\to Y_{\bar I}$ and
$X_{\bar I}\to Y_I$, together with their inverses, have Lipschitz
constant $O(Ls^{t+1})=\widetilde O(n^{1/3})$.
\end{claim}

\footnotetext{The bijections map $x_{\min I+j}$ to $y_{\min\bar I+j}$,
and $x_{\min\bar I+j}$ to $y_{\min I+j}$, respectively,
for $0\le j<|I|$.}

\begin{proof}
Assume that the event described in \Cref{obs:gap-regularity} holds.
The definition of $e_i^Z$ and the bound $E_i^Z\le4L$ give
\begin{equation}
    s^{-(e_i^Z+1)}<E_i^Z\le4Ls^{-e_i^Z}
 \qquad(i\in[n-1],\ Z\in\{X,Y\}),
 \label{eq:E_i^Z_bound}
\end{equation}
including when $e_i^Z=0$.

Fix intervals and patterns satisfying the assumptions.
For $i,i+1\in I$, put $\bar{\imath}=i+\bar r-r$.
The condition on the roots implies that the increasing bijection
$X_I\to Y_{\bar I}$ sends $x_i,x_{i+1}$ to
$y_{\bar{\imath}},y_{\bar{\imath}+1}$.

First, assume that $i \not \in \operatorname{supp}(P)$. Then $\bar{\imath}\not \in \operatorname{supp}(\bar P)$ as well, since the patterns $P$ and $\bar P$ have the same shape. Thus, neither $i$ nor $\bar{\imath}$ are tiny, since all tiny indices inside the two intervals
are recorded in $\operatorname{supp}(P)$ and $\operatorname{supp}(\bar P)$, by assumption. Consequently, due to the event from \Cref{obs:gap-regularity}, both gaps
$\Delta_i^X,\Delta_{\bar{\imath}}^Y$ lie between $1/(2ns^t)$ and $8L/n$,
so their ratio in either direction is at most $16Ls^t$.

Now suppose that $i \in \operatorname{supp} P$. The low-expansion condition gives
$|e_i^X-e_{\bar{\imath}}^Y|\le t$.
We have
\[
 \frac{\Delta_{\bar{\imath}}^Y}{\Delta_i^X}
 =\frac{E_{\bar{\imath}}^Y/\Sigma_Y}{E_i^X/\Sigma_X}
 =\frac{E_{\bar{\imath}}^Y}{E_i^X}\,
   \frac{\Sigma_X}{\Sigma_Y}
 \le16Ls^{t+1}.
\]
The inequality follows from~\eqref{eq:E_i^Z_bound} and the bounds on $\Sigma_X,\Sigma_Y$ in
\Cref{obs:gap-regularity}. The same bound holds for the ratio in the other direction.

Since the bijection is increasing, it preserves the order of the points
and sends consecutive points to consecutive points. By \Cref{cl:adjacent-gap}, the Lipschitz constants of the bijection and its inverse are each attained on a pair of consecutive points, so both are at most $16Ls^{t+1}$.

The same argument applies to both increasing bijections and
their inverses, giving the claimed bounds, since
$Ls^{t+1}=\widetilde O(n^{1/3})$. 
\end{proof}

We also record a uniform bound on the number of tiny indices in short
intervals.

\begin{claim}\label{cl:tiny-indices-in-search-intervals}
With high probability, every interval $J\subseteq[n-1]$ of length at
most $3H$ contains at most $12H/s^t$ tiny indices.
\end{claim}

\begin{proof}
Each index is tiny with probability $1-e^{-2s^{-t}}\le2/s^t$,
independently of the other indices. Thus, for each fixed interval $J$
of length at most $3H$, the number $N_J$ of tiny indices is binomial
with mean at most $6H/s^t$. A binomial tail bound and a union bound
over the at most $n^2$ intervals give, for an absolute constant $c>0$,
\[
 \mathbb P\bigl(\exists J\subseteq[n-1],\ |J|\le3H:
                   N_J>12H/s^t\bigr)
 \le n^2e^{-cH/s^t}=o(1),
\]
where $J$ ranges over intervals and we used $H/s^t\gg\log n$.
\end{proof}

\subsection{Properties of typical patterns}\label{subsec:ub-typical-patterns}

To prove \Cref{prop:matching-patterns}, we will choose pairs of intervals
to which \Cref{obs:cover-lipschitz} applies. In particular, to apply the
claim, the corresponding patterns must record all tiny gaps inside
these intervals. We also require the intervals to be separated by at
least $D/2$ gap indices for later stages of the proof. Consequently,
we have the following definition.

\begin{samepage}
\begin{definition}[Components]\label{def:components}
A \emph{component} (of $X$ \emph{and} $Y$) is a nonempty subset
$C=\{c_1<\cdots<c_k\}$ of $[n-1]$ such that:
\begin{enumerate}
\item\label{it:ub-component-tiny} The index $c_i$ is tiny for every $i\in[k]$.
\item\label{it:ub-component-distances} $c_{i+1}-c_i\le D$ for every $i\in[k-1]$.
\item\label{it:ub-component-isolation} Every tiny index in $[\min C-D,\max C+D]$ belongs to $C$.
\end{enumerate}
For a component $C=\{c_1<\cdots<c_k\}$, let
$P_C$ be the pattern with $\alpha_j^Z=e_{c_j}^Z$ for
$j\in[k]$, $Z\in\{X,Y\}$, and $d_j=c_{j+1}-c_j$ for $j\in[k-1]$.
\end{definition}
\end{samepage}

We call an occurrence $(r,P)$ \emph{isolated} if every tiny index within
distance $D$ of its gap-index interval $\operatorname{int}(r,P)$ belongs
to its support.

For each component $C$, the pattern $P_C$ occurs at root $\min C$, with
support $C$ and gap-index interval $[\min C,\max C]$.
By \Cref{it:ub-component-isolation}, this occurrence is isolated.

We first show that occurring patterns whose recorded indices are tiny
and sufficiently close have at most $L$ recorded positions.

\begin{claim}\label{cl:pattern-size}
Let $P=(\alpha_1^X,\alpha_1^Y,d_1,\alpha_2^X,\alpha_2^Y,
\ldots,d_{k-1},\alpha_k^X,\alpha_k^Y)$ be a pattern, and assume that
\begin{enumerate}
\item $d_j\in[12D]$ for every $j\in[k-1]$;
\item $\max\{\alpha_j^X,\alpha_j^Y\}\ge t$ for every $j\in[k]$.
\end{enumerate}
Then, with high probability, \emph{no} such pattern $P$ with $k>L$ occurs at an admissible root.
\end{claim}

\begin{proof}
For each $i\in[n-1]$, the independence of $E_i^X$ and $E_i^Y$ gives
$\mathbb P(i\text{ is tiny})=1-e^{-2s^{-t}}\le2s^{-t}$.
It suffices to consider $k=L+1$; otherwise, take the subpattern
consisting of the first $L+1$ positions of $P$.

We now take a union bound over all such patterns.
There are $12D$ choices for each $d_j$, hence at most $(12D)^L$ shapes.
There are at most $n$ choices for the root $r$.
For a fixed shape and an admissible root, the probability that all
$L+1$ recorded gap indices are tiny is at most $(2s^{-t})^{L+1}$,
by independence across these indices.
Therefore, by a union bound, the probability that one of these
patterns occurs is at most
\[
 n(12D)^L(2s^{-t})^{L+1}
 \le 2ns^{-t}(24/s)^L=o(1),
\]
where we used $D/s^t\le1/s$, $s=L^3$, and $L=\lceil\log n\rceil$.
Hence, with high probability, none of the aforementioned patterns
occurs, completing the proof.
\end{proof}

Let $\mathcal P$ be the family of patterns $P$ of size $1\le k\le L$
with $d_j\le D$ for every $j\in[k-1]$, and with
$0\le\alpha_j^X,\alpha_j^Y\le L$ and
$\max\{\alpha_j^X,\alpha_j^Y\}\ge t$ for every $j\in[k]$.

\begin{corollary}\label{cor:component-size}
With high probability, $P_C\in\mathcal P$ for every component $C$.
\end{corollary}

\begin{proof}
By \Cref{cl:pattern-size}, with high probability every component $C$ has
size at most $L$: its pattern $P_C$ occurs and satisfies the hypotheses of
that claim. Moreover, for every $i\in[n-1]$ and $Z\in\{X,Y\}$,
\[
 \mathbb P(e_i^Z>L)
 =\mathbb P\bigl(E_i^Z\le s^{-(L+1)}\bigr)
 \le s^{-(L+1)}\le n^{-3},
\]
for sufficiently large $n$. A union bound over the gap indices in both
samples shows that all exponents are at most $L$ with probability
at least $1-2/n^2$. The remaining conditions defining $\mathcal P$
hold by the definition of a component, proving the corollary.
\end{proof}

We aim to bound the size of $\mathcal{P}$. Notice that every pattern $P \in \mathcal{P}$ can be constructed by choosing its size, shape, and exponents. For each size $k \in [L]$, there are at most $D^{k-1}$ choices of the distances
and at most $(L+1)^{2k}$ choices of the exponents. Hence
\begin{equation}\label{eq:pattern-family-cardinality}
 |\mathcal P|\le\sum_{k=1}^{L}D^{k-1}(L+1)^{2k}
 =\exp\bigl((1/3+o(1))L^2\bigr),
\end{equation}
where we used $\log D=(1/3+o(1))L$ and $\log(L+1)=o(L)$.

For each component pattern, we will seek a low-expansion cover with the
largest probability of occurring at a fixed root.
For $a\in\mathbb N_0$, let $\gamma_a=\mathbb P(e_s(E)=a)$,
where $E\sim\operatorname{Exp}(1)$.
Thus $\gamma_0=e^{-1/s}$ and
$\gamma_a=e^{-s^{-(a+1)}}-e^{-s^{-a}}$ for $a\ge1$.
For a fixed pattern $P$ and a fixed root $r$ for which $(r,P)$ is defined,
the recorded exponential variables are independent, so let us denote
\[
 \operatorname{prob}(P):= \mathbb P\bigl((r,P)\text{ occurs}\bigr)
=\prod_{j=1}^{|P|}\gamma_{\alpha_j^X}\gamma_{\alpha_j^Y}.
\]

For sufficiently large $n$, $\gamma_a$ is strictly decreasing in $a$.
Thus, among the low-expansion covers of $P$, the occurrence probability
is maximized by choosing each exponent as small as allowed.
This motivates the following definition.

\begin{definition}[Most probable cover]\label{def:most-probable-cover}
For a pattern $P$, we denote its most probable low-expansion cover
by $\operatorname{MPC}(P)$;
it has the same size and shape as $P$ and exponent pairs
$((\alpha_j^Y-t)_+,(\alpha_j^X-t)_+)$ for $j\in[|P|]$.
Put $q_P=\operatorname{prob}(\operatorname{MPC}(P))$.
\end{definition}

We next compare the occurrence probability of each pattern in
$\mathcal P$ with that of its most probable cover.

\begin{claim}\label{cl:pattern-probability-ratio}
There is an absolute constant $C>0$ such that every pattern
$P\in\mathcal P$ of size $k$ satisfies
\[
 \frac{\operatorname{prob}(P)}{q_P}\le\left(\frac{C}{s^t}\right)^k.
\]
\end{claim}
\begin{proof}
Recall that $\gamma_0=\Theta(1)$ and
$\gamma_a=\Theta(s^{-a})$ uniformly for $a\ge1$.
Since $\gamma_a$ is decreasing, the ratio
$\gamma_a/\gamma_{(a-t)_+}$ is at most $1$ for every $a\ge0$.
For $a\ge t$, the same ratio is at most $C/s^t$, for an absolute
constant $C$.

By the definition of the most probable cover,
\[
 \frac{\operatorname{prob}(P)}{q_P}
 =\prod_{j=1}^{k}
 \frac{\gamma_{\alpha_j^X}}{\gamma_{(\alpha_j^X-t)_+}}
 \frac{\gamma_{\alpha_j^Y}}{\gamma_{(\alpha_j^Y-t)_+}}.
\]
At each recorded position, at least one of $\alpha_j^X,\alpha_j^Y$
is at least $t$, since $P\in\mathcal P$.
Thus each pair of factors is at most $C/s^t$, proving the claim.
\end{proof}

For a pattern $P$ of size $k$ with $d_j\le D$ for every $j\in[k-1]$, put
\[
 w_P=C_1L^2\left(\frac{pk^2}{q_P}+kD\right),
\]
where $C_1$ is a fixed, sufficiently large absolute constant.
In \Cref{subsec:ub-matching-proof}, $w_P$ will be the length of an interval in which
we search for an isolated occurrence of $\operatorname{MPC}(P)$.
The first term accounts for the occurrence probability of the cover,
and the second provides separation for the matching intervals.
The following lemma bounds the first contribution globally.

\begin{lemma}\label{lem:global-component-search-cost}
With high probability,
\[
 C_1L^2p\sum_C\frac{|C|^2}{q_{P_C}}\le H/L,
\]
where the sum is over all components.
\end{lemma}

\begin{proof}
Before bounding the sum in the lemma, which is \emph{over components},
we bound a related sum over all patterns in $\mathcal P$. This bound
expresses that patterns with large $1/q_P$ are sufficiently unlikely
to occur. Specifically, we show that
\[
 \sum_{P\in\mathcal P}|P|^2
 \frac{\operatorname{prob}(P)}{q_P}=O(L^2/s^t).
\]
There are at most $D^{k-1}(L+1)^{2k}$ patterns of size $k$ in $\mathcal{P}$:
each of the $k-1$ successive distances has at most $D$ choices,
and each of the $2k$ exponents has at most $L+1$ choices.
By \Cref{cl:pattern-probability-ratio},
\begin{align*}
 \sum_{P\in\mathcal P}|P|^2\frac{\operatorname{prob}(P)}{q_P}
 &\le\sum_{k=1}^{L}k^2D^{k-1}(L+1)^{2k}
       \left(\frac{C}{s^t}\right)^k\\
 &=\frac{C(L+1)^2}{s^t}
   \sum_{k=1}^{L}k^2\left(\frac{CD(L+1)^2}{s^t}\right)^{k-1}
 =O(L^2/s^t).
\end{align*}
We only need to explain the last equality. Since $s^t\ge sD$ and $s=L^3$, we have
$CD(L+1)^2/s^t=O(1/L)$. For sufficiently large $n$, the last
sum is at most $\sum_{k\ge1}k^2 2^{-(k-1)}=O(1)$.

To bound the sum over components, define the random variable $S$ by:
\[
 S=C_1L^2p\sum_{P\in\mathcal P}\frac{|P|^2}{q_P}
 \sum_{r=1}^{n-1-\operatorname{span}(P)}
 \mathbf 1_{\{(r,P)\text{ occurs}\}},
\]
where an inner sum is empty if there is no valid root.
For each valid root $r$, the probability of $(r,P)$ occurring is
$\operatorname{prob}(P)$. By linearity of expectation and
the preceding estimate,
\[
 \mathbb E S\le C_1L^2pn\sum_{P\in\mathcal P}|P|^2
 \frac{\operatorname{prob}(P)}{q_P}
 =O(nL^4p/s^t)=O(nL^5/s^t).
\]
Applying Markov's inequality, we get
\[
 \mathbb P(S>H/L)
 \le\frac{L\mathbb E S}{H}
 =O\left(\frac{nL^6}{Hs^t}\right)
 =O\left(\frac{nL^3}{HD}\right)=O(L^{-4})=o(1),
\]
using $s^t\ge sD$, $s=L^3$, and $HD=\Theta(nL^7)$.

By \Cref{cor:component-size}, with high probability every
component pattern belongs to $\mathcal P$. On this event, each component
$C$ gives a distinct rooted occurrence $(\min C,P_C)$ counted by $S$.
Hence $C_1L^2p\sum_C|C|^2/q_{P_C}\le S$.
Intersecting this event with $S\le H/L$ proves the lemma.
\end{proof}

\subsection{Proof of \texorpdfstring{\Cref{prop:matching-patterns}}{the matching proposition}}\label{subsec:ub-matching-proof}

We use the properties established in \Cref{subsec:ub-typical-patterns} to construct
the matching. We first partition the indices $[n]=J_1\mathbin{\dot\cup}\cdots\mathbin{\dot\cup}J_m$ while keeping components
away from the boundaries of the partition $\{J_i\}_{i=1}^m$. We then find isolated covers, allocate disjoint
search intervals, and select and align pairs of patterns.

\subsubsection{Partition into intervals}

\begin{claim}\label{cl:component-partition}
With high probability, there is a partition
$[n]=J_1\mathbin{\dot\cup}\cdots\mathbin{\dot\cup}J_m$ into consecutive
intervals of sizes between $H$ and $3H$ such that, for every component
$C$, some $J_j$ satisfies
\[
 [\min C,\max C+1]\subseteq[\min J_j+4D,\max J_j-4D].
\]
\end{claim}

\begin{proof}
We start with a deterministic partition
$[n]=\widehat J_1\mathbin{\dot\cup}\cdots\mathbin{\dot\cup}\widehat J_m$
into $m=\lfloor n/(2H)\rfloor$ consecutive intervals whose sizes differ
by at most one, each of size $2H+o(H)$.
It is easy to see that such a partition always exists.
Write $\widehat J_j=[b_{j-1}+1,b_j]$, where
$0=b_0<b_1<\cdots<b_m=n$.
We now slightly perturb the boundaries $b_j$ to new boundaries $b'_j$
and define $J_j=[b'_{j-1}+1,b'_j]$.

Condition on the event in \Cref{cl:pattern-size}.
Keep the two outer boundaries fixed: $b'_0=0$ and $b'_m=n$.
For each $j\in[m-1]$, if $[b_j-5D,b_j+5D]$ contains no tiny index,
set $b'_j=b_j$.
Otherwise, all tiny indices in this interval lie within distance $10D$
of one another, so they belong to one maximal sequence of tiny indices
with successive distances at most $12D$.
By \Cref{cl:pattern-size}, this sequence has at most $L$ indices.
Let $u_j$ be the last index of this sequence, and set $b'_j=u_j+5D+1$.
Then $0\le b'_j-b_j\le12LD=o(H)$.
Moreover, every tiny index is either at most $u_j$ or greater than
$u_j+12D$, so $b'_j$ is more than $5D$ from every tiny index.

Since the original intervals have size $2H+o(H)$ and every boundary
moves by $o(H)$, we have $0=b'_0<b'_1<\cdots<b'_m=n$, and the intervals
$J_j=[b'_{j-1}+1,b'_j]$ form a partition with sizes between $H$ and $3H$.
By construction, there are no tiny indices in the $5D$-neighborhood
of $b'_j$ for any $j\in[m-1]$.

It remains to verify that components in $J_1$ and $J_m$ lie at the
required distance from the outer endpoints of these intervals,
corresponding to the boundaries $b'_0=0$ and $b'_m=n$.
To handle these cases, observe that a simple union bound over the first
and last $8D$ gap indices shows that none of them is tiny with probability
$1-O(D/s^t)=1-o(1)$.
This gives the remaining margins for $J_1$ and $J_m$, completing the proof.
\end{proof}

\subsubsection{Isolated covers}

Recall the definition of the most probable cover and the search length $w_P$ defined in
\Cref{subsec:ub-typical-patterns}. We show that intervals of length at least
$w_P$ contain suitable isolated occurrences.

The following simple claim shows that from an interval $V$ of gap indices
we can extract many pairwise disjoint sets with the same shape as a given
pattern $P$. Their roots lie in a prescribed residue class modulo $p$, and
the $D$-neighborhoods of their corresponding intervals lie in $V$.
Concretely:

\begin{claim}[Disjoint supports]\label{cl:disjoint-pattern-supports}
Let $P$ be a pattern of size $k$ with $d_j\le D$ for every $j\in[k-1]$,
let $K\ge1$ be an integer, and let
$V\subseteq[n-1]$ be an interval of gap indices with
\[
 |V|\ge pk^2K+(k+1)D.
\]
For every $r_0\in\{0,\ldots,p-1\}$, there are $K$ rooted patterns
$(r_1,P),\ldots,(r_K,P)$ such that:
\begin{enumerate}
\item\label{it:ub-supports-disjoint} Their supports are pairwise disjoint.
\item\label{it:ub-supports-residue} $r_i\equiv r_0\pmod p$ for every $i\in[K]$.
\item\label{it:ub-supports-interior} The $D$-neighborhood of $\operatorname{int}(r_i,P)$ is contained
in $V$ for every $i\in[K]$.
\end{enumerate}
\end{claim}

\begin{proof}
Fix $r_0\in\{0,\ldots,p-1\}$.
Let $V'$ denote the set of roots $r$ such that the $D$-neighborhood of
$\operatorname{int}(r,P)$ lies in $V$. Notice that $V'$ is an interval of consecutive
indices. We have
\[
 |V'|=|V|-\operatorname{span}(P)-2D
 \ge |V|-(k+1)D\ge pk^2K,
\]
where we used $\operatorname{span}(P)\le(k-1)D$, since each $d_j\le D$.
Let $V''=\{r\in V':r\equiv r_0\pmod p\}$.
Since $V'$ is an interval, $|V''|\ge k^2K$.

We claim that there are $K$ roots in $V''$ whose corresponding supports
are pairwise disjoint.
Indeed, define a graph $G$ on $V''$, joining distinct vertices $u,v\in V''$
if the supports of $(u,P)$ and $(v,P)$ intersect.
Two such supports intersect only if $u-v=a_i-a_j$ for some $i,j\in[k]$
(recall that $a_i=\sum_{h<i}d_h$, as defined after \Cref{def:patterns}).
There are at most $k^2$ such differences, including $0$, so the maximum
degree of $G$ satisfies $\Delta(G)<k^2$.

Greedily choosing a vertex of $G$ and deleting it and its neighbors gives an
independent set of size at least
$|V''|/(\Delta(G)+1)\ge |V''|/k^2\ge K$.
Choose any $K$ vertices $r_1,\ldots,r_K$ in this independent set.
It is easy to verify by construction that \Cref{it:ub-supports-disjoint,it:ub-supports-residue,it:ub-supports-interior} hold for these roots, completing the proof.
\end{proof}

We now apply \Cref{cl:disjoint-pattern-supports} to search for isolated
occurrences of most probable covers. We first fix a pattern in the family
$\mathcal P$; a union bound will
then give the conclusion simultaneously for all component patterns.

\begin{claim}\label{cl:isolated-cover-probability}
Fix a pattern $P\in\mathcal P$, a residue class modulo $p$, and an interval
$W\subseteq[n-1]$ with $|W|\ge w_P$.
Then, with probability at least $1-2^{-10L^2}$, there is an isolated
occurrence of $\operatorname{MPC}(P)$ whose root lies in the prescribed
residue class and whose interval's $D$-neighborhood is contained in $W$.
\end{claim}

\begin{proof}
Put $k=|P|$, $q=q_P$, and $K=\lceil4/q\rceil$, so $4\le Kq\le5$.
Split $W$ into intervals of length
\[
 m=pk^2K+(k+1)D\le5(pk^2/q+kD),
\]
discarding the remainder. By the definition of $w_P$, there are at least
$10L^2$ such intervals for sufficiently large $C_1$.
It suffices to show that each fixed subinterval contains a suitable
isolated occurrence, with its entire $D$-neighborhood inside that
subinterval, with probability at least $1/2$.
Indeed, these events depend only on the exponential variables inside
their respective subintervals and are therefore independent. The
probability that all fail is then at most $2^{-10L^2}$.

Fix one such subinterval $U$.
Apply \Cref{cl:disjoint-pattern-supports} to $\operatorname{MPC}(P)$
to choose deterministic roots $r_1,\ldots,r_K$ in the prescribed residue class.
Their supports are pairwise disjoint, and the $D$-neighborhoods of their
intervals lie in $U$.
Let $A_i$ be the event that $(r_i,\operatorname{MPC}(P))$ occurs.
The events $A_i$ are independent, each with probability $q$, so the
probability that none occurs is at most $e^{-Kq}\le e^{-4}$.

Let $T_i$ be the set of gap indices in the $D$-neighborhood of
$\operatorname{int}(r_i,\operatorname{MPC}(P))$ that do not belong to
its support. This neighborhood contains $\operatorname{span}(P)+1+2D$
indices, of which $k$ belong to the support. Thus
\[
 |T_i|=\operatorname{span}(P)+1+2D-k\le(k+1)D,
\]
since $\operatorname{span}(P)\le(k-1)D$.
Let $B_i$ be the event that some index in $T_i$ is tiny.
For each $j\in T_i$, the variables $E_j^X,E_j^Y$ are independent of
$A_i$. The index $j$ is tiny if at least one of these variables is at
most $s^{-t}$, so $\mathbb P(j\text{ is tiny}\mid A_i)\le2s^{-t}$.
Recall that $\mathbb P(A_i)=q$ by definition.
Applying a union bound over $T_i$ gives
\begin{align*}
 \mathbb P(A_i\cap B_i)
 =\mathbb P(A_i)\,\mathbb P(B_i\mid A_i)
 =q\,\mathbb P(B_i\mid A_i)
 &\le q\sum_{j\in T_i}\mathbb P(j\text{ is tiny}\mid A_i)\\
 &\le\frac{2q|T_i|}{s^t}
 \le\frac{2q(k+1)D}{s^t}.
\end{align*}
By a union bound, the probability that no candidate gives an isolated
occurrence is at most
\[
 e^{-4}+\sum_{i=1}^{K}\mathbb P(A_i\cap B_i)
 \le e^{-4}+\frac{10(k+1)D}{s^t}
 =e^{-4}+O(L^{-2})<\frac12,
\]
where we used $k\le L$ and $D/s^t\le1/s=L^{-3}$.
This proves the claim.
\end{proof}

\begin{lemma}[Isolated covers]\label{lem:isolated-covers}
With high probability, for every component $C$, every residue class modulo
$p$, and every interval $V\subseteq[n-1]$ with $|V|\ge w_{P_C}$,
there is an isolated occurrence of $\operatorname{MPC}(P_C)$ with root
in that residue class, whose interval and its $D$-neighborhood are contained
in $V$.
\end{lemma}

\begin{proof}
By \Cref{cl:isolated-cover-probability}, a union bound over
$P\in\mathcal P$, the $p$ residue classes, and the at most $n^2$ intervals $V$
gives total failure probability at most
\[
 |\mathcal P|pn^2\,2^{-10L^2}
 \le\exp\bigl((1/3-10\log 2+o(1))L^2\bigr)=o(1),
\]
by \Cref{eq:pattern-family-cardinality}. By \Cref{cor:component-size},
with high probability every component pattern belongs to $\mathcal P$,
so the lemma follows.
\end{proof}

\subsubsection{Selecting and aligning pairs}

We now choose component--cover pairs and use them to define the partial
block mapping. We first assign each component an isolated cover so that
the cover intervals are mutually separated. A cover may still lie close
to another component; we resolve these conflicts by selecting pairs in
decreasing order of the components' largest exponents.

\begin{lemma}\label{lem:selecting-pattern-pairs}
With high probability, the partition $J_1,\ldots,J_m$ from
\Cref{cl:component-partition} admits pairs of occurring rooted patterns
$\mathbf P_i,\bar{\mathbf P}_i$, for $i\in[\ell]$ and some
$\ell\in\mathbb N_0$, such that:
\begin{enumerate}
\item\label{it:ub-selection-cover} Each $\mathbf P_i$ records a component, and
$\bar{\mathbf P}_i$ is an isolated occurrence of its most probable
cover, with root congruent to that of $\mathbf P_i$ modulo $p$.
\item\label{it:ub-selection-separation} The $2\ell$ gap-index intervals are pairwise separated by more
than $D$, and every tiny index belongs to one of their supports.
\item\label{it:ub-selection-margins} For each $i\in[\ell]$, both corresponding gap-index intervals
lie in $[\min J_j+4D,\max J_j-4D]$ for some $j\in[m]$.
\end{enumerate}
\end{lemma}

\begin{proof}
Assume that the events described in
\Cref{cl:component-partition,cl:tiny-indices-in-search-intervals,lem:global-component-search-cost,lem:isolated-covers} hold. Fix the
partition $J_1,\ldots,J_m$ from \Cref{cl:component-partition}. We claim the following auxiliary property: for every $j \in [m]$ it holds that
\begin{equation}
    \sum_{C\subseteq J_j}|C| \leq 12H/s^t.
    \label{eq:sum_of_components_sizes}
\end{equation}
That is because the components $C\subseteq J_j$ in the above sum are disjoint and consist of tiny indices. Therefore the sum on the left-hand side counts the number of tiny indices in $J_j$, which is at most $12H/s^t$ due to the event from \Cref{cl:tiny-indices-in-search-intervals}.

We allocate to every component $C$ an interval of length
$\lceil w_{P_C}\rceil$ in which we want to allocate $C$'s points to an isolated
occurrence of $\operatorname{MPC}(P_C)$ in the interval.
For a fixed part $J_j$, the components are disjoint sets of tiny indices.
Since $w_{P_C}\ge D\ge1$ we have
\begin{align*}
 \sum_{C\subseteq J_j}\bigl(\lceil w_{P_C}\rceil+D\bigr)
 &\le3\sum_{C\subseteq J_j}w_{P_C}\\
 &=3C_1L^2p\sum_{C\subseteq J_j}\frac{|C|^2}{q_{P_C}}
   +3C_1L^2D\sum_{C\subseteq J_j}|C|\\
 &\stackrel{\cref{lem:global-component-search-cost}}\le3H/L+3C_1L^2D\sum_{C\subseteq J_j}|C|\\
 &\stackrel{\eqref{eq:sum_of_components_sizes}}{\le}3H/L+36C_1L^2DH/s^t=O(H/L),
\end{align*}
where the last equality uses $D/s^t\le L^{-3}$.
Thus intervals $V_C$ of lengths $\lceil w_{P_C}\rceil$, with $D$ unused
gap indices between consecutive intervals, fit inside
$[\min J_j+4D,\max J_j-4D-1]$: this interval has at least
$H-8D-1$ indices, and $D=o(H)$.
The intervals $V_C$ are separated by more than $D$, also across
distinct parts because of the margins.

By \Cref{lem:isolated-covers}, in each $V_C$ choose an isolated
occurrence $\mathbf Q_C$ of $\operatorname{MPC}(P_C)$ whose root is
congruent to $\min C$ modulo $p$ and whose gap-index interval and its
$D$-neighborhood lie in $V_C$.
Observe that even when the intervals $\{V_C\}_C$ are mutually separated,
they need not be separated from the components. For example, if
$E_i^X\approx 1/n$, a low-expansion cover must match this gap to a gap
of size at most $\widetilde O(n^{-2/3})$. The latter gap is still tiny,
so its index belongs to a component. The following observation will
allow us to resolve such conflicts by selecting pairs in decreasing
order of the components' largest exponents.

For every component $C$, put $h(C)=\max_{i\in C}\max\{e_i^X,e_i^Y\}$.
Fix a component $C$, and let $C'$ be any component whose gap-index
interval is within distance $D$ of the gap-index interval of
$\mathbf Q_C$. We claim that
\[
 C'\subseteq\operatorname{supp}(\mathbf Q_C),
 \qquad h(C')<h(C).
\]
Indeed, $\mathbf Q_C$ is isolated by \Cref{lem:isolated-covers}, so
every tiny index within distance $D$ of
$\operatorname{supp}(\mathbf Q_C)$ belongs to this support.
Successive indices in both $C'$ and $\operatorname{supp}(\mathbf Q_C)$
are at distance at most $D$. Thus the assumed proximity of their
gap-index intervals implies that some index of $C'$ is within
distance $D$ of $\operatorname{supp}(\mathbf Q_C)$.
Since all indices of $C'$ are tiny, it follows that
$C'\cap\operatorname{supp}(\mathbf Q_C)\ne\varnothing$.
Starting from an index in this intersection, each successive index
of $C'$ in either direction is tiny and within distance $D$ of the
previous one. Isolation therefore puts each of these indices in
$\operatorname{supp}(\mathbf Q_C)$, proving
$C'\subseteq\operatorname{supp}(\mathbf Q_C)$.

By the definition of $\operatorname{MPC}(P_C)$, each exponent recorded
by $\mathbf Q_C$ is of the form $(a-t)_+$, where $a$ is an exponent
recorded by $P_C$. Since $a\le h(C)$ and $h(C)\ge t$, every such
exponent is at most $h(C)-t$. We have shown that every index of $C'$
belongs to the support of $\mathbf Q_C$, so both exponents at every
index of $C'$ satisfy this bound. Thus
\[
 h(C')\le h(C)-t<h(C).
\]

We now describe how to match components to low-expansion covers so
that all components are covered.
Process the components in decreasing order of $h(C)$, breaking ties
arbitrarily. Let $(\mathbf P_{C_1},\mathbf Q_{C_1}),\ldots,
(\mathbf P_{C_j},\mathbf Q_{C_j})$ be the pairs selected so far,
where initially $j=0$, and let $C$ be the next component in this order.
If the gap-index interval of $C$ is within distance $D$ of the
gap-index interval of $\mathbf Q_{C_i}$ for some $i\in[j]$, skip $C$.
Otherwise, set $C_{j+1}=C$ and select the pair
$(\mathbf P_{C_{j+1}},\mathbf Q_{C_{j+1}})$, where
$\mathbf P_C=(\min C,P_C)$.

We now verify the three properties in the lemma. Let $\ell$ be the
number of selected pairs, and relabel them as
$(\mathbf P_i,\bar{\mathbf P}_i)$ for $i\in[\ell]$.

\Cref{it:ub-selection-cover} holds by construction: each selected $\mathbf P_C$
records a component, and $\mathbf Q_C$ is an isolated occurrence
of $\operatorname{MPC}(P_C)$ with root congruent to that of
$\mathbf P_C$ modulo $p$.

\Cref{it:ub-selection-separation} follows from the selection rule and the containment and
strict decrease proved above: the selected intervals are pairwise
separated by more than $D$, and every skipped component is contained in an
already selected cover. Thus every tiny index belongs to a selected support.

For \Cref{it:ub-selection-margins}, consider a selected pair generated by
$C\subseteq J_j$. The component's gap-index
interval has the required margins by \Cref{cl:component-partition}.
The cover's gap-index interval has the same margins as $V_C\subseteq
[\min J_j+4D,\max J_j-4D]$.
\end{proof}

\begin{proof}[Proof of \Cref{prop:matching-patterns}]
Assume that the high-probability events of
\Cref{lem:selecting-pattern-pairs,cor:component-size,obs:cover-lipschitz} hold.
By \cref{lem:selecting-pattern-pairs}, we select pairs of rooted patterns $\mathbf P_i,\bar{\mathbf P}_i$ such that $\mathbf P_i$ records a component of tiny indices and $\bar{\mathbf P}_i$ is a most probable cover of $\mathbf P_i$. For each selected pair of rooted patterns $\mathbf P_i,\bar{\mathbf P}_i$, let $I_i$ be the smallest union of consecutive blocks containing $j, j+1$ for all indices $j$ recorded by $\mathbf P_i$, and let $\bar I_i$ be defined the same way with respect to $Y$ and $\bar{\mathbf P}_i$.

By \Cref{cor:component-size}, each component has at most $L$ indices
and span at most $(L-1)D$. Hence $|I_i|,|\bar I_i|\le LD+2p+1=O(LD)$.
The original gap-index intervals of $\mathbf P_i$ and $\bar{\mathbf P}_i$  are separated by at least $D$ and have margins at least $4D$ in their respective parts. Thus the enlarged intervals are separated by at least $D-2p-1\ge D/2$ and retain margins
at least $4D-p-1\ge D/2$. They contain no terminal block, since that
block lies within $p$ indices of $n$, whereas these intervals remain
more than $p$ indices from $n$. This proves the interval requirements
in \Cref{it:ub-matching-separation,it:ub-matching-margins}.

Define $f$ by the increasing bijections $X_{I_i}\to Y_{\bar I_i}$
and $X_{\bar I_i}\to Y_{I_i}$ for each $i$.
The intervals $I_i$ and $\bar I_i$ are disjoint unions of nonterminal blocks,
and corresponding blocks of $I_i$ and $\bar I_i$ are exchanged in pairs. Thus $f$ is a partial
block mapping induced by a partial involution, proving \Cref{it:ub-matching-involution}.
Every tiny index $j$ lies in a selected support, so $j,j+1$ lie together
in its enlarged interval, proving \Cref{it:ub-matching-tiny}.

Finally, both occurrences in each pair are isolated: one records a
component, and the other was chosen to be isolated. Since the extensions
are shorter than $p<D$, every tiny gap index inside either enlarged
interval lies in the support of its pattern.
The paired patterns are low-expansion covers, and their roots have
the same relative positions in the enlarged intervals.
Therefore \Cref{obs:cover-lipschitz} applies and gives
\Cref{it:ub-matching-lipschitz}, completing the proof.
\end{proof}

\subsection{Rearrangement lemma}

We now rearrange the free points (those whose indices lie in free blocks)
so that gaps between consecutive free blocks have length
$\widetilde{\Theta}(1/n)$, while fixing every prescribed point.

The main idea is that we need not construct the matching directly between
$X$ and $Y$. We may first pass to perturbed copies $X'$ and $Y'$ that
have essentially the same metric: each copy is related to the original
sample by a bijection under which every distance is preserved up to a
polylogarithmic factor. A matching between these copies can then be
transferred back to $X$ and $Y$ with only a polylogarithmic increase in
distortion. Concretely, we seek copies in the following \emph{canonical form}.

\begin{definition}[Canonical form]\label{def:canonical-form}
Let $\mathcal I$ be a family of disjoint intervals in $[n]$, each a union
of nonterminal blocks. Put $U=\bigcup_{I\in\mathcal I}I$, and suppose
that there is at least one block outside $U$ between consecutive intervals.
Let $f:X_U\to Y_U$ be a partial block mapping.
We say that $(X,Y,f;\mathcal I)$ is in
\emph{canonical form} if, for $Z\in\{X,Y\}$ and $i\in[n-1]$,
\begin{enumerate}
\item\label{it:ub-canonical-upper} $\Delta_i^Z\le s^2/n$;
\item\label{it:ub-canonical-nontiny} if $\{i,i+1\}\nsubseteq U$, then $\Delta_i^Z\ge1/(2ns^t)$;
\item\label{it:ub-canonical-boundary} if $i,i+1\notin U$ and $p\mid i$, then $\Delta_i^Z\ge1/(sn)$.
\end{enumerate}
\end{definition}

Notice that \Cref{it:ub-canonical-upper} requires every gap to have length
$\widetilde O(1/n)$. This property holds whp for the original samples by \Cref{cr:maximum-gap} and will remain true after rearrangement.
\Cref{it:ub-canonical-nontiny} reflects that, in our construction, we pass to
canonical form \emph{only after} obtaining the partial mapping $f$ from
\Cref{prop:matching-patterns}. Every tiny gap then lies inside a
prescribed interval. Thus every gap with at least one endpoint indexed
outside $U$ is nontiny, giving the lower bound in \Cref{it:ub-canonical-nontiny} on the
event $\Sigma_X,\Sigma_Y\le2n$.
\Cref{it:ub-canonical-boundary} is the property we seek to enforce: every gap
between consecutive free blocks has length $\widetilde\Omega(1/n)$.
In the original random samples, gaps at block boundaries can be much
smaller.

\noindent\begin{minipage}{\linewidth}
To achieve canonical form, we will use the following claim.

\begin{claim}[Deterministic rearrangement]\label{cl:deterministic-rearrangement}
Let $A=\{a_1<\cdots<a_\ell\}\subseteq\mathbb R$, where $\ell\ge3$,
and let $T>0$. Suppose at least a $0.9$-fraction of the gaps
$a_{i+1}-a_i$, $i\in[\ell-1]$, are at least $T$, including the first
and last gaps.

Then for every $1\le k<\ell$, there exist an ordered set
$A'=\{a'_1<\cdots<a'_\ell\}$ and a bijection $R:A\to A'$ such that:
\begin{enumerate}
\item\label{it:ub-rearrange-endpoints} $R$ fixes the endpoints, which remain the endpoints of $A'$, i.e.,
$R(a_1)=a'_1=a_1$ and \mbox{$R(a_\ell)=a'_\ell=a_\ell$}.
\item\label{it:ub-rearrange-target} $a'_{k+1}-a'_k\ge T$.
\item\label{it:ub-rearrange-multiset} The multiset of gap lengths is preserved.
\item\label{it:ub-rearrange-lipschitz} We have $\max\{\Lip(R),\Lip(R^{-1})\}\le(a_\ell-a_1)/T$.
\end{enumerate}
\end{claim}
\end{minipage}

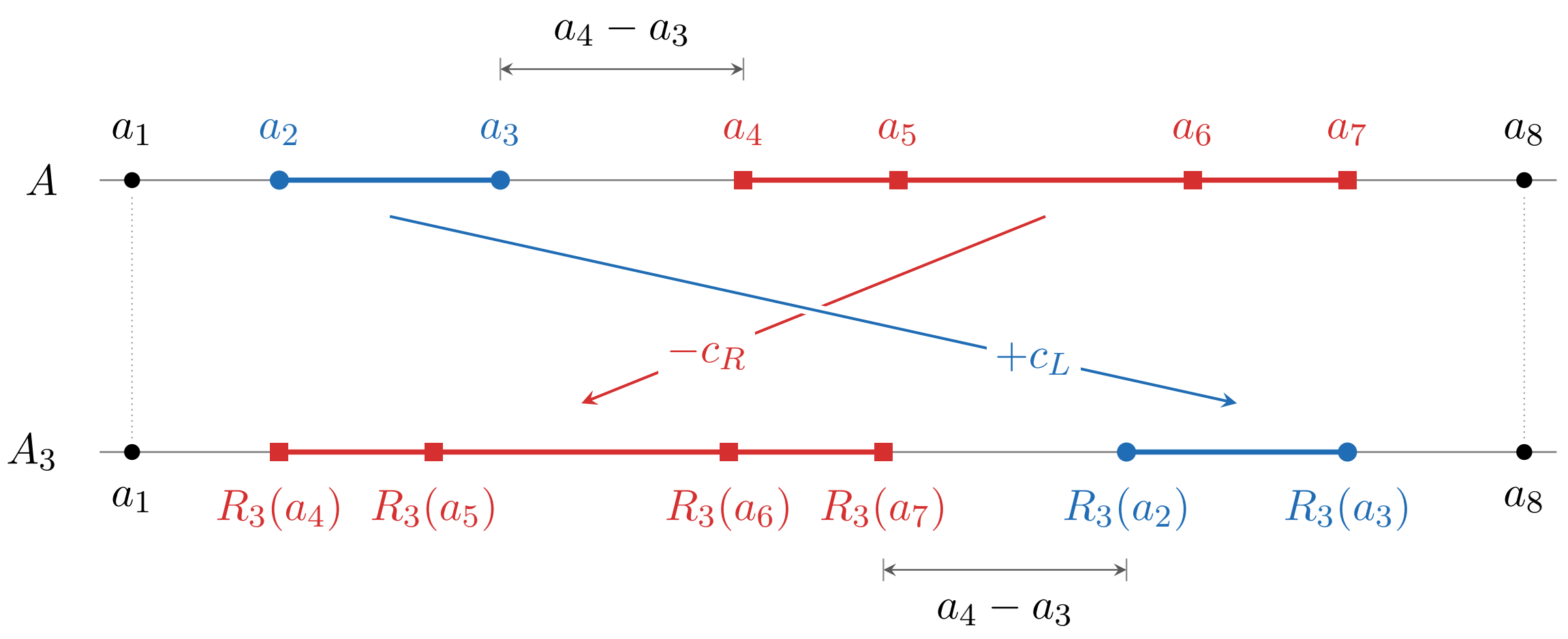
\begin{figure}[H]
\centering
\scalebox{1.1}{
\begin{tikzpicture}[x=1.3cm,y=1cm,>=stealth,font=\normalsize]
\definecolor{leftgroup}{HTML}{206DB5}
\definecolor{rightgroup}{HTML}{D62F2F}

\foreach \x in {0,9.45}
  \draw[black!35,densely dotted] (\x,0.12)--(\x,2.28);
\foreach \y in {0,2.4}
  \draw[black!45,line width=0.5pt] (-0.22,\y)--(9.67,\y);
\node[anchor=east] at (-0.4,2.4) {$A$};
\node[anchor=east] at (-0.4,0) {$A_3$};

\draw[leftgroup,line width=1.3pt] (1,2.4)--(2.5,2.4);
\draw[rightgroup,line width=1.3pt] (4.15,2.4)--(8.25,2.4);
\foreach \x/\i in {0/1,9.45/8} {
  \fill (\x,2.4) circle[radius=2pt];
  \node[above=5pt] at (\x,2.4) {$a_{\i}$};
}
\foreach \x/\i in {1/2,2.5/3} {
  \fill[leftgroup] (\x,2.4) circle[radius=2.4pt];
  \node[above=5pt,text=leftgroup] at (\x,2.4) {$a_{\i}$};
}
\foreach \x/\i in {4.15/4,5.2/5,7.2/6,8.25/7} {
  \node[fill=rightgroup,rectangle,inner sep=0pt,minimum size=4.5pt] at (\x,2.4) {};
  \node[above=5pt,text=rightgroup] at (\x,2.4) {$a_{\i}$};
}

\draw[->,rightgroup,line width=0.7pt] (6.2,2.08)--(3.05,0.43)
  node[pos=0.73,fill=white,inner sep=2pt] {$-c_R$};
\draw[->,leftgroup,line width=0.7pt,
      preaction={draw=white,line width=3pt}]
  (1.75,2.08)--(7.5,0.43)
  node[pos=0.76,fill=white,inner sep=2pt] {$+c_L$};

\draw[rightgroup,line width=1.3pt] (1,0)--(5.1,0);
\draw[leftgroup,line width=1.3pt] (6.75,0)--(8.25,0);
\foreach \x/\i in {0/1,9.45/8} {
  \fill (\x,0) circle[radius=2pt];
  \node[below=5pt] at (\x,0) {$a_{\i}$};
}
\foreach \x/\i in {6.75/2,8.25/3} {
  \fill[leftgroup] (\x,0) circle[radius=2.4pt];
  \node[below=5pt,text=leftgroup] at (\x,0) {$R_3(a_{\i})$};
}
\foreach \x/\i in {1/4,2.05/5,4.05/6,5.1/7} {
  \node[fill=rightgroup,rectangle,inner sep=0pt,minimum size=4.5pt] at (\x,0) {};
  \node[below=5pt,text=rightgroup] at (\x,0) {$R_3(a_{\i})$};
}

\draw[<->,black!65] (2.5,3.38)--(4.15,3.38)
  node[midway,above=2pt,text=black] {$a_4-a_3$};
\foreach \x in {2.5,4.15}
  \draw[black!45] (\x,3.28)--(\x,3.48);
\draw[<->,black!65] (5.1,-1.04)--(6.75,-1.04)
  node[midway,below=2pt,text=black] {$a_4-a_3$};
\foreach \x in {5.1,6.75}
  \draw[black!45] (\x,-0.94)--(\x,-1.14);
\end{tikzpicture}
}
\caption{Eight points before (top) and after (bottom) applying $R_3$.
In the application, $A$ is a set of consecutive points from one sample,
$X$ or $Y$; the rows show its original and rearranged configurations.
The endpoints $a_1,a_8$ are fixed. The blue interval translates right by
$c_L=a_7-a_3$, and the red interval translates left by $c_R=a_4-a_2$,
preserving all distances within each interval. The marked gap has length
$a_4-a_3$ in both rows.}
\label{fig:rearrangement-r3}
\end{figure}

\begin{proof}
For $1\le j\le\ell-2$, put
\[
 c_L=a_{\ell-1}-a_j,
 \qquad
 c_R=a_{j+1}-a_2.
\]
Define the $j$th rearrangement map $R_j$ by $R_j(a_1)=a_1$,
$R_j(a_\ell)=a_\ell$, and
\[
 R_j(a_i)=
 \begin{cases}
 a_i+c_L, & 2\le i\le j,\\
 a_i-c_R, & j+1\le i\le\ell-1.
 \end{cases}
\]
Set $A_j=R_j(A)=\{b_1^{(j)}<\cdots<b_\ell^{(j)}\}$.
We claim that, for some $j$, $A_j$ is the desired set $A'$ and
$R:=R_j$ is the desired map.

It is easy to check that \Cref{it:ub-rearrange-endpoints,it:ub-rearrange-multiset}
hold for every $j$ by construction; see \Cref{fig:rearrangement-r3}.

To verify \Cref{it:ub-rearrange-target}, we choose $j$ avoiding both
of the following conditions:
\[
 \begin{aligned}
 a_{j+1}-a_j&<T &&\text{(short cut gap)},\\
 b_{k+1}^{(j)}-b_k^{(j)}&<T &&\text{(short target gap)}.
 \end{aligned}
\]
Let $r$ be the number of gaps in $A$ smaller than $T$; then $r\le(\ell-1)/10$.
At most $r$ values of $j$ satisfy the first condition.
Among the remaining values of $j$, at most $r$ satisfy the second:
each original short gap occupies the target gap position for at most
one value of $j$. Since $2r<\ell-2$, we may choose $j$ satisfying
\Cref{it:ub-rearrange-target} and $a_{j+1}-a_j\ge T$.
Fix such a $j$ and set $A'=A_j$ and $R=R_j$.

It remains to verify \Cref{it:ub-rearrange-lipschitz}. Distances within
$\{a_2,\ldots,a_j\}$ and within $\{a_{j+1},\ldots,a_{\ell-1}\}$
are preserved by $R$. For every other pair of distinct points,
both its distance in $A$ and the distance between its images in
$A'$ lie between $T$ and $a_\ell-a_1$. Thus both Lipschitz
constants are at most $(a_\ell-a_1)/T$.
\end{proof}

We will apply \Cref{cl:deterministic-rearrangement} locally, and
\emph{separately} to $X$ and $Y$. The following claim gives the required
gap bounds.

\begin{claim}\label{cl:rearrangement-gap-bounds}
Put $T=1/(sn)$. With high probability, the following hold simultaneously
for $Z\in\{X,Y\}$:
\begin{enumerate}
\item\label{it:ub-gap-bounds-density} Among the $p-1$ gaps between any $p$ consecutive points,
at least a $0.9$-fraction have length at least $T$.
\item\label{it:ub-gap-bounds-terminal} The last $p$ gaps, $\Delta_{n-p}^Z,\ldots,\Delta_{n-1}^Z$,
are at least $T$.
\end{enumerate}
\end{claim}

\begin{proof}
By \Cref{obs:gap-regularity}, $\Sigma_X,\Sigma_Y\le2n$ with high
probability. When $\Sigma_Z\le2n$, a gap smaller than $T$ must have
$E_i^Z<2/s$. The events $E_i^Z<2/s$ are independent as $i$ and $Z$
vary, each occurring with probability at most $2/s$.

To prove \Cref{it:ub-gap-bounds-density}, it suffices to take a union bound over the at most
$n$ windows of $p-1$ consecutive gaps in each sample. The probability
that some window contains at least $(p-1)/10$ gaps smaller than $T$,
while $\Sigma_X,\Sigma_Y\le2n$, is at most
\[
 2n\,2^{p-1}\left(\frac{2}{s}\right)^{\lceil(p-1)/10\rceil}
 =n e^{-\Omega(p\log s)}=o(1).
\]
For \Cref{it:ub-gap-bounds-terminal}, a union bound over the last $p$ gaps in both samples
bounds the probability of failure while $\Sigma_X,\Sigma_Y\le2n$
by $2p\cdot(2/s)=4p/s=o(1)$. This proves the claim.
\end{proof}

\noindent\begin{minipage}{\linewidth}
\begin{claim}\label{cl:rearrangement-windows}
Put $T=1/(sn)$ and fix $Z\in\{X,Y\}$.
On the event of \Cref{cl:rearrangement-gap-bounds}, for each block
boundary $b=jp$ with $\Delta_b^Z<T$, there is an interval $A_b$ of
at most $p$ consecutive points containing $z_b,z_{b+1}$.
The intervals $A_b$ are pairwise disjoint, and at least a $0.9$-fraction
of the gaps in each are at least $T$, including the first and last.
\end{claim}
\end{minipage}

\begin{proof}
Fix $Z\in\{X,Y\}$. The gap preceding the terminal block, at
$b=(M-1)p$, is among the last $p$ gaps and hence is at least $T$ by
\Cref{it:ub-gap-bounds-terminal}. Consequently, every
boundary with gap smaller than $T$ lies between two full blocks.
For each such boundary $b=jp$, let $\widehat A_b$ consist of the last
$p/2$ points of $B_j^Z$ and the first $p/2$ points of $B_{j+1}^Z$.

By \Cref{it:ub-gap-bounds-density}, $\widehat A_b$ has
at most $(p-1)/10$ short gaps. Each half has $p/2-1>(p-1)/10$ gaps,
so it contains a gap of length at least $T$.
Trim $\widehat A_b$ at both ends until its first and last gaps are at
least $T$, obtaining $A_b$.
Since both the left and right parts of $\widehat A_b$ contain a gap of
length at least $T$, $A_b$ contains $z_b,z_{b+1}$.
Since we removed only short gaps, at least a $0.9$-fraction of the
gaps of $A_b$ are at least $T$.
Finally, the $\widehat A_b$ are pairwise disjoint, and hence so are
the $A_b$.
\end{proof}

\noindent\begin{minipage}{\linewidth}
\begin{lemma}[Rearrangement]\label{lem:rearrangement}
With high probability, for every $f,\mathcal I$ satisfying the conclusions
of \Cref{prop:matching-patterns}, with $U=\bigcup_{I\in\mathcal I}I$,
there exist ordered sets
$X'=\{x'_1<\cdots<x'_n\}$, $Y'=\{y'_1<\cdots<y'_n\}$ and
bijections $R_Z:Z\to Z'$, for $Z\in\{X,Y\}$, such that
\[
 R_Z(z_i)=z'_i=z_i\qquad(i\in U),
\]
and $\max\{\Lip(R_Z),\Lip(R_Z^{-1})\}=O(L^5)$.
Moreover, $f'=R_Y\circ f\circ R_X^{-1}$ is induced by the same partial
involution as $f$, and $(X',Y',f';\mathcal I)$ is in canonical form.
\end{lemma}
\end{minipage}

\begin{proof}
Put $T=1/(sn)$. Work on the events of
\Cref{obs:gap-regularity,cl:rearrangement-gap-bounds},
and fix $f,\mathcal I$ satisfying the conclusions of
\Cref{prop:matching-patterns}. We construct $R_X$ and $R_Y$ separately;
fix $Z\in\{X,Y\}$.

For each boundary $b$ between free blocks with $\Delta_b^Z<T$,
take $A_b$ from \Cref{cl:rearrangement-windows}. These intervals are
pairwise disjoint and contain only free points, i.e., points in free blocks (recall \Cref{def:block-mappings}).
Applying \Cref{cl:deterministic-rearrangement} to each $A_b$, we obtain
a map $R_b$ such that the gap at index $b$ after rearrangement is at
least $T$, and
\[
 \max\{\Lip(R_b),\Lip(R_b^{-1})\}
 \le\frac{\operatorname{diam}(A_b)}{T}
 \le8psL=O(L^5),
\]
where we used $|A_b|\le p$ and \Cref{obs:gap-regularity}.

Apply all these maps simultaneously, fixing every point outside the
intervals $A_b$, to obtain $R_Z$, and put $Z'=R_Z(Z)$.
Since the $A_b$ are pairwise disjoint and each $R_b$ fixes its endpoints, \Cref{cl:adjacent-gap} bounds $\Lip(R_Z)$ by the maximum of $1$ and the $\Lip(R_b)$;
the same holds for the inverses. Consequently,
\[
 \max\{\Lip(R_Z),\Lip(R_Z^{-1})\}\le8psL=O(L^5).
\]

Each $R_b$ rearranges only free points within $A_b$, so the points with indices in $U$ retain their positions and indices. Hence
$R_Z(z_i)=z'_i=z_i$ for $i\in U$, and $f'$ is induced by the same
partial block permutation as $f$.

We verify the three conditions of canonical form.
\begin{enumerate}
\item By \Cref{it:ub-rearrange-multiset}, and \Cref{cl:deterministic-rearrangement} the
multiset of gap lengths is preserved, so every gap remains at most
$8L/n\le s^2/n$.

\item By \Cref{it:ub-matching-tiny} of \Cref{prop:matching-patterns}, $E_i^Z>s^{-t}$
whenever $\{i,i+1\}\nsubseteq U$. Since $\Sigma_Z\le2n$, each such
gap initially has length at least $1/(2ns^t)$. These gaps are either
unchanged or permuted within intervals of free points, so the lower
bound persists.

\item Each boundary gap shorter than $T$ between free blocks becomes
at least $T$. Every other block-boundary gap is unchanged, since each
$R_b$ fixes its endpoints and acts within the two blocks adjacent
to $b$. Thus every boundary between free blocks has gap at least $T$.
\end{enumerate}
Hence $(X',Y',f';\mathcal I)$ is in canonical form.
\end{proof}

\subsection{Reduction to permutations}

We now use \Cref{lem:rearrangement} to show that, after the pattern
matching step of \Cref{prop:matching-patterns}, it remains to find a
suitable block extension of the partial mapping $f$ constructed there.
We recall the relevant setup and conclusions in the following proposition.

For $I\subseteq[n]$, recall that $B(I)=\{j\in[M]:B_j\subseteq I\}$.
For a family $\mathcal I$ of prescribed intervals with union $U$,
we have $B(U)=\bigcup_{I\in\mathcal I}B(I)$.

\begin{proposition}\label{prop:reduction-to-permutations}
Let $X,Y\subseteq[0,1]$ have size $n$, and let the family
$\mathcal I$ and the partial block mapping $f:X_U\to Y_U$ satisfy
the conclusions of \Cref{prop:matching-patterns}, where
$U=\bigcup_{I\in\mathcal I}I$.
Assume that the conclusion of \Cref{lem:rearrangement} holds for
$(X,Y,f;\mathcal I)$.
Let
$\sigma:B(U)\to B(U)$ be the partial involution inducing $f$
as in \Cref{def:block-mappings}.
Suppose that there exists a permutation $\pi$ of $[M]$ such that
\begin{enumerate}
\item $\pi$ extends $\sigma$ and fixes $M$.
\item $\max\{\Lip(\pi),\Lip(\pi^{-1})\}=O(D)$.
\item\label{it:neighbour_identitiy} For every $I\in\mathcal I$, putting $a=\min B(I)$ and
$b=\max B(I)$, we have
\[
 \pi(a-1)=\sigma(a)-1,\qquad
 \pi(b+1)=\sigma(b)+1.
\]
\end{enumerate}
Then there is a bijection $F:X\to Y$ extending $f$ with
$\dist(F)=\widetilde O(n^{2/3})$.
\end{proposition}

\begin{proof}
Our plan is to use the block mapping $F'$ induced by $\pi$ between
$X'=R_X(X)$ and $Y'=R_Y(Y)$, where $R_X,R_Y$ are provided by
\Cref{lem:rearrangement}. After bounding the Lipschitz constants
of $F'$ and its inverse, we transfer this map back by setting
$F=R_Y^{-1}\circ F'\circ R_X$. The polylogarithmic Lipschitz bounds
for the rearrangement maps and their inverses then give the desired
distortion bound for $F$.

Apply \Cref{lem:rearrangement} to obtain $X',Y'$, bijections
$R_X,R_Y$, and $f'=R_Y\circ f\circ R_X^{-1}$ such that
$(X',Y',f';\mathcal I)$ is in canonical form.
The map $f'$ is induced by the same $\sigma$.
Since the rearrangement fixes every point whose index lies in $U$, \Cref{it:ub-matching-lipschitz} gives
\[
 \max\{\Lip(f'|_{X'_I}),\Lip((f')^{-1}|_{Y'_I})\}
 =O(Ls^{t+1})\qquad(I\in\mathcal I).
\]
Let $F':X'\to Y'$ be the block mapping induced by $\pi$.

We bound the expansion of consecutive points $x'_i,x'_{i+1}$ under $F'$.
We distinguish four cases according to their locations.
\begin{enumerate}
\item \emph{Both indices lie in $U$.}
Since the prescribed intervals $\mathcal{I}$ are separated, both indices belong to
the same $I\in\mathcal I$. The map $F'$ agrees with $f'$ there,
so their expansion is $O(Ls^{t+1})$ by \Cref{it:ub-matching-lipschitz} of \Cref{prop:matching-patterns}.

\item \emph{Both points lie in the same free block.}
Their images lie in the same block of $Y'$ and are consecutive,
since $F'$ maps each block increasingly. \Cref{it:ub-canonical-upper,it:ub-canonical-nontiny} of \Cref{def:canonical-form} (canonical form) give
\[
 x'_{i+1}-x'_i\ge\frac{1}{2ns^t},\qquad
 |F'(x'_{i+1})-F'(x'_i)|\le\frac{s^2}{n}.
\]
Thus their expansion is at most $(s^2/n)/(1/(2ns^t))=2s^{t+2}$.

\item \emph{Exactly one index lies in $U$.}
The points lie across a boundary of some $I\in\mathcal I$.
Put $a=\min B(I)$ and $b=\max B(I)$.
Here we use the crucial neighbor identities $\pi(a-1)=\sigma(a)-1$
and $\pi(b+1)=\sigma(b)+1$. They ensure that $F'$ sends the blocks
immediately before and after $X'_I$ to the blocks immediately before
and after $f'(X'_I)$, respectively. Since each block is mapped increasingly,
the two points across either boundary have consecutive images.
\Cref{it:ub-canonical-upper,it:ub-canonical-nontiny} of \Cref{def:canonical-form} gives
$x'_{i+1}-x'_i\ge1/(2ns^t)$ and
$|F'(x'_{i+1})-F'(x'_i)|\le s^2/n$, so their expansion is at most
$2s^{t+2}$.

\item \emph{The points lie in different free blocks.}
Here $i=jp$ for some $j\in[M-1]$.
Crucially, because we first applied the rearrangement, \Cref{it:ub-canonical-boundary} of \Cref{def:canonical-form} gives
\[
 x'_{jp+1}-x'_{jp}\ge\frac{1}{sn}.
\]
The image indices are $\pi(j)p$ and $(\pi(j+1)-1)p+1$.
This formula also holds when $j+1=M$, since $\pi(M)=M$.
The upper gap bound in \Cref{it:ub-canonical-upper} therefore gives
\[
 \frac{|F'(x'_{jp+1})-F'(x'_{jp})|}{x'_{jp+1}-x'_{jp}}
 \le s^3p\bigl(|\pi(j+1)-\pi(j)|+1\bigr)
 \le2s^3p\Lip(\pi).
\]
\end{enumerate}
Combining the cases and applying \Cref{cl:adjacent-gap} yields
\[
 \Lip(F')\le\max\{O(Ls^{t+1}),2s^{t+2},2s^3p\Lip(\pi)\}.
\]

\noindent By symmetry, the same bound holds for $\Lip((F')^{-1})$ with
$\pi^{-1}$ in place of $\pi$.

Substituting $\max\{\Lip(\pi),\Lip(\pi^{-1})\}=O(D)$,
$s^t\le s^2D$, $s=L^3$, and $p=C_0L$ into these bounds gives
\[
 \max\{\Lip(F'),\Lip((F')^{-1})\}=\widetilde O(D).
\]
Finally, transfer the mapping back to the original samples by setting
$F=R_Y^{-1}\circ F'\circ R_X$.
Since $F'$ extends $f'$, the map $F$ extends $f$.
The $O(L^5)$ Lipschitz bounds for the
rearrangement maps and their inverses imply
\[
 \dist(F)=\Lip(F)\Lip(F^{-1})
 =\widetilde O(D^2)=\widetilde O(n^{2/3}).\qedhere
\]
\end{proof}

\subsection{Oscillations}

\begin{lemma}[Oscillations on intervals]\label{lem:interval-oscillations}
Let $m\ge1$ and let $I_1,\ldots,I_m$ be nonempty disjoint integer intervals listed in
increasing order, with $|I_j|>10m$. Write $\ell_j=\min I_j$,
$u_j=\max I_j$, and let $\tau$ be an involution of $[m-1]$.
Then there exists an involution $\pi$ of $I_1\cup\cdots\cup I_m$ satisfying
\[
 \begin{aligned}
 \pi(\ell_1)&=\ell_1,&
 \pi(\ell_j)&=\ell_{\tau(j-1)+1} &&(2\le j\le m),\\
 \pi(u_m)&=u_m,&
 \pi(u_j)&=u_{\tau(j)} &&(1\le j<m).
 \end{aligned}
\]
If $i,i+1\in I_j$ and $\pi(i)\in I_k$, $\pi(i+1)\in I_{k'}$,
then $|k-k'|\le1$.
\end{lemma}

\begin{figure}[H]
\centering
\begingroup
\definecolor{paironetwo}{HTML}{206DB5}
\definecolor{paironethree}{HTML}{D62F2F}
\definecolor{pairtwothree}{HTML}{D5A000}
\definecolor{ink}{HTML}{293440}
\definecolor{muted}{HTML}{89929C}
\resizebox{.99\linewidth}{!}{%
\begin{tikzpicture}[x=1cm,y=.70cm,font=\normalsize,text=ink,
  imagearrow/.style={->,>=stealth,line width=.5pt},
  point/.style={circle,minimum size=13pt,inner sep=0pt,fill=white,
    draw=muted!85,line width=.5pt,font=\normalsize},
  endpoint/.style={point,draw=ink,line width=.9pt}]
\useasboundingbox (-.10,1.55) rectangle (16.45,15.55);

\foreach \i in {1,2,3}{
  \pgfmathsetmacro{\a}{.7+5.45*(\i-1)}
  \pgfmathsetmacro{\b}{\a+4.24}
  \coordinate (left\i) at (\a,14.05);
  \coordinate (right\i) at (\b,14.05);
  \coordinate (leftimage\i) at (\a,11.50);
  \coordinate (rightimage\i) at (\b,11.50);
  \foreach \y in {14.05,11.50}{
    \draw[muted!85,line width=.5pt] (\a,\y)--(\b,\y);
    \node[fill=white,inner sep=1pt] at ({(\a+\b)/2},\y) {$I_{\i}$};
    \fill[ink] (\a,\y) circle[radius=1.8pt];
    \fill[ink] (\b,\y) circle[radius=1.8pt];
  }
  \node[anchor=south] at (\a,14.27) {$\ell_{\i}$};
  \node[anchor=south] at (\b,14.27) {$u_{\i}$};
  \node[anchor=north] at (\a,11.25) {$\ell_{\i}$};
  \node[anchor=north] at (\b,11.25) {$u_{\i}$};
}
\foreach \y in {14.05,11.50}{
  \draw[muted!60,densely dotted] (4.99,\y)--(6.10,\y);
  \draw[muted!60,densely dotted] (10.44,\y)--(11.55,\y);
}
\node[font=\normalsize,text=ink] at (5.545,14.43) {$1$};
\node[font=\normalsize,text=ink] at (10.995,14.43) {$2$};
\draw[<->,>=stealth,draw=ink,line width=.5pt]
  (5.545,14.95)--node[fill=white,inner sep=2pt] {$\tau=(1,2)$} (10.995,14.95);
\foreach \from/\to in {left1/leftimage1,right1/rightimage2,right2/rightimage1,
  left2/leftimage3,left3/leftimage2,right3/rightimage3}{
  \draw[imagearrow,ink,shorten <=1.8pt,shorten >=1.8pt]
    (\from)--(\to);
}

\foreach \i in {1,2,3}{
  \pgfmathsetmacro{\a}{.7+5.45*(\i-1)}
  \pgfmathsetmacro{\b}{\a+4.24}
  \draw[muted!70,line width=.45pt] (\a,9.12)--(\a,9.24)--(\b,9.24)--(\b,9.12);
  \node at ({(\a+\b)/2},9.51) {$I_{\i}$};
  \draw[muted!70,line width=.45pt] (\a,3.20)--(\a,3.08)--(\b,3.08)--(\b,3.20);
  \node at ({(\a+\b)/2},2.78) {$I_{\i}$};
  \foreach \k in {1,...,9}{
    \coordinate (source-\i-\k) at ({\a+.53*(\k-1)},8.65);
    \coordinate (target-\i-\k) at ({\a+.53*(\k-1)},3.70);
  }
}
\foreach \i/\k/\lab in {
  1/1/1,1/2/1,1/3/1,1/4/2,1/5/2,1/6/3,1/7/3,1/8/3,1/9/2,
  2/1/3,2/2/2,2/3/1,2/4/1,2/5/2,2/6/3,2/7/3,2/8/2,2/9/1,
  3/1/2,3/2/1,3/3/1,3/4/1,3/5/2,3/6/2,3/7/3,3/8/3,3/9/3}{
  \ifnum\k=1\def\pointstyle{endpoint}\else
    \ifnum\k=9\def\pointstyle{endpoint}\else\def\pointstyle{point}\fi\fi
  \node[\pointstyle] (sourcepoint-\i-\k) at (source-\i-\k) {$\lab$};
  \node[\pointstyle] (targetpoint-\i-\k) at (target-\i-\k) {$\lab$};
}
\foreach \i/\k in {1/2,1/3,2/2,2/5,2/8,3/7,3/8}{
  \draw[imagearrow,draw=muted!45,line width=.4pt]
    (sourcepoint-\i-\k.south)--(targetpoint-\i-\k.north);
}
\foreach \i/\a/\j/\b/\col in {
  1/4/2/3/paironetwo,1/5/2/4/paironetwo,
  1/6/3/2/paironethree,1/7/3/3/paironethree,1/8/3/4/paironethree,
  2/6/3/5/pairtwothree,2/7/3/6/pairtwothree}{
  \draw[imagearrow,draw=\col]
    (sourcepoint-\i-\a.south)--(targetpoint-\j-\b.north);
  \draw[imagearrow,draw=\col]
    (sourcepoint-\j-\b.south)--(targetpoint-\i-\a.north);
}
\foreach \i/\a/\j/\b in {1/1/1/1,3/9/3/9,1/9/2/9,2/9/1/9,2/1/3/1,3/1/2/1}{
  \draw[imagearrow,draw=ink,line width=.6pt]
    (sourcepoint-\i-\a.south)--(targetpoint-\j-\b.north);
}
\node[anchor=east,text=muted,font=\normalsize] at (.34,6.18) {$\pi$};
\foreach \x/\i/\j/\col in {.80/1/2/paironetwo,6.25/1/3/paironethree,11.40/2/3/pairtwothree}{
  \draw[draw=\col,line width=.65pt] (\x,1.98)--({\x+.52},1.98);
  \node[anchor=west,font=\small] at ({\x+.65},1.98)
    {$C'_{\i\j}\leftrightarrow C'_{\j\i}$};
}
\end{tikzpicture}
}
\endgroup
\caption{An illustration of the construction in
\Cref{lem:interval-oscillations} for $\tau=(1,2)$, with nine points per
interval and $r_{ij}=3$ for $i,j=1,2,3$. The endpoint assignments (black) extend to
an involution by pairing the point of $C'_{ij}$ with those of $C'_{ji}$ (colored arrows) and fixing
$C'_{ii}$ (gray arrows).}
\label{fig:interval-oscillations}
\end{figure}
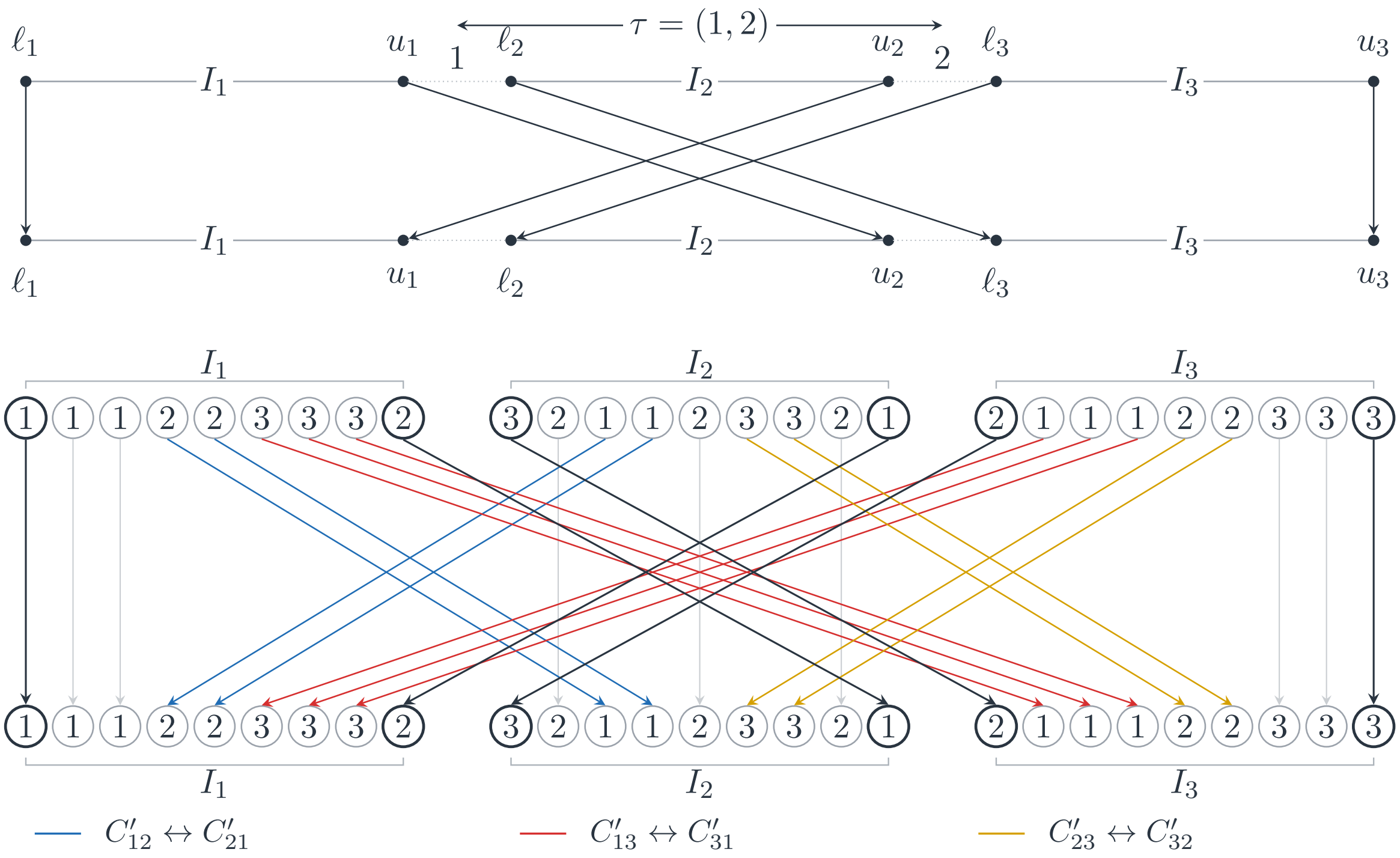

\begin{proof}
The reader may consult \Cref{fig:interval-oscillations} throughout the proof.
Let $s_i,t_i$ be the indices of the intervals to which the left and
right endpoints of $I_i$ must be mapped, respectively. Thus
\[
 \begin{aligned}
 s_1&=1,& s_i&=\tau(i-1)+1 &&(2\le i\le m),\\
 t_m&=m,& t_i&=\tau(i) &&(1\le i<m).
 \end{aligned}
\]
\noindent We construct $\pi$ in two steps: first a labeling of $I_1\cup\cdots\cup I_m$ by $[m]$, then an involution respecting that labeling.

\emph{Constructing the labeling.}
Given a labeling of $I_1\cup\cdots\cup I_m$ by $[m]$, write
\[
 C_{ij}=\{x\in I_i:x\text{ is labeled }j\}.
\]
Set
\[
 r_{ij}=3\quad(i\ne j),\qquad r_{ii}=|I_i|-3(m-1).
\]
Then $r_{ij}=r_{ji}\ge3$ and $\sum_j r_{ij}=|I_i|$.
We will arrange that $|C_{ij}|=r_{ij}$ for every $i,j\in[m]$.
For each $i$, concatenate the three walks
\[
 (s_i,s_i-1,\ldots,1),\qquad
 (1,2,\ldots,m),\qquad
 (m,m-1,\ldots,t_i).
\]
Each label occurs between one and three times. Since $r_{ij}\ge3$,
insert consecutive repetitions so that label $j$ occurs exactly $r_{ij}$ times.
The resulting sequence has length $|I_i|$; use it to label the points of
$I_i$ in increasing order.
\begin{samepage}
The resulting labeling has the following properties:
\begin{enumerate}
\item\label{it:ub-labeling-endpoints} The points $\ell_i$ and $u_i$ have labels $s_i$ and $t_i$, respectively.
\item\label{it:ub-labeling-neighbors} Consecutive points of $I_i$ have labels differing by at most one.
\item\label{it:ub-labeling-counts} $|C_{ij}|=r_{ij}=r_{ji}=|C_{ji}|$ for every $i,j\in[m]$.
\end{enumerate}
\end{samepage}

\emph{Defining the involution.}
We now use this labeling to construct $\pi$.
First prescribe $\pi(\ell_i)=\ell_{s_i}$ and $\pi(u_i)=u_{t_i}$,
and put $C'_{ij}=C_{ij}\setminus\{\ell_i,u_i\}$.
Since $\tau$ is an involution, $s_i=j$ if and only if $s_j=i$, and
$t_i=j$ if and only if $t_j=i$. Thus removing the endpoints preserves
the equality in \Cref{it:ub-labeling-counts}, so $|C'_{ij}|=|C'_{ji}|$.

For each $i<j$, match $C'_{ij}$ with $C'_{ji}$ arbitrarily and let
$\pi$ exchange each matched pair. Let $\pi$ fix every point of $C'_{ii}$.
Together with the endpoint assignments, this defines an involution
with the required endpoint values. Every point labeled $j$ is mapped
into $I_j$, so \Cref{it:ub-labeling-neighbors} gives the last assertion of this lemma.
\end{proof}

\begin{lemma}[Oscillating completion]\label{lem:oscillating-completion}
Let $J$ be an interval of integers, and let $A_1,\ldots,A_h$ be
nonempty pairwise disjoint subintervals of $J$, where $h\ge0$.
Let $\sigma$ be an involution on $A_1\cup\cdots\cup A_h$ such that
$\sigma$ maps each $A_i$ increasingly and bijectively onto some $A_j$.

Suppose that, for positive integers $b,q$ with $q^2\ge16|J|$,
\begin{itemize}
\item For every $i$, $|A_i|\le b$.
\item $J\setminus(A_1\cup\cdots\cup A_h)=B_1\cup\cdots\cup B_{h+1}$,
where the $B_j$ are intervals interlacing the $A_i$ in the order
$B_1,A_1,B_2,\ldots,A_h,B_{h+1}$, and each $B_j$ has size at least $q$.
\end{itemize}
Then $\sigma$ extends to an involution $\pi$ of $J$ fixing both
endpoints and satisfying
\[
 \Lip(\pi)\le4q+b.
\]
Moreover, whenever $\sigma(A_i)=A_j$, the map $\pi$ sends the point
immediately before $A_i$ to the point immediately before $A_j$,
and the point immediately after $A_i$ to the point immediately
after $A_j$.
\end{lemma}

\begin{proof}
Partition each $B_j$ into pieces of sizes between $q$ and $2q-1$:
take pieces of size $q$ and add the
remainder to the last piece. List all pieces in increasing order as
$I_1,\ldots,I_m$, and write $\ell_j=\min I_j$, $u_j=\max I_j$.
Since $mq\le|J|$ and $q^2\ge16|J|$, we have
$|I_j|\ge q\ge16m>10m$.

We now define $\pi$. We first define an involution $\tau$ on $[m-1]$,
which we will use with \Cref{lem:interval-oscillations} to construct $\pi$.
We set $\tau(j)=k$ if $A_i$ lies between $I_j$ and
$I_{j+1}$ and $\sigma(A_i)$ lies between $I_k$ and $I_{k+1}$.
For every remaining gap between $I_j$ and $I_{j+1}$,
both pieces lie in the same $B_r$; set $\tau(j)=j$.
Since $\sigma$ is an involution, so is $\tau$.
Apply \Cref{lem:interval-oscillations} to obtain an involution $\pi$
on $I_1\cup\cdots\cup I_m$, and extend $\pi$ by $\sigma$ on
$A_1\cup\cdots\cup A_h$. This defines an involution of $J$ extending $\sigma$.

For the last assertion of the lemma, consider the images of the points
immediately before and after each $A_i$.
By the endpoint conditions in \Cref{lem:interval-oscillations},
$\pi$ fixes both endpoints of $J$. By our choice of $\tau$, whenever
$\sigma(A_i)=A_j$, these conditions also give
\begin{equation}\label{eq:ub-completion-neighbors}
 \begin{aligned}
 \pi(\min A_i-1)&=\min A_j-1,\\
 \pi(\max A_i+1)&=\max A_j+1.
 \end{aligned}
\end{equation}

Finally, we bound the Lipschitz constant.
Fix consecutive points $x,x+1\in J$. We distinguish three cases according
to the partition of $J$ into the intervals $A_i$ and $I_j$.
\begin{enumerate}
\item \emph{Both points lie in one $A_i$.}
Since $\pi$ agrees with the increasing bijection $\sigma$ on $A_i$,
their images are consecutive, so $|\pi(x+1)-\pi(x)|=1$.

\item \emph{The points lie in different intervals of the partition.}
If one belongs to an $A_i$, \Cref{eq:ub-completion-neighbors} and the fact that
$\sigma$ maps $A_i$ increasingly onto its image show that their images
are consecutive. Otherwise, the points lie across a subdivision cut
within one $B_r$. At such a cut, $\tau(j)=j$, so
the involution supplied by \Cref{lem:interval-oscillations} fixes both $u_j=x$ and $\ell_{j+1}=x+1$.
Thus the image distance is again $1$.

\item \emph{Both points lie in one $I_j$.}
By \Cref{lem:interval-oscillations}, their images lie in the same interval or in two
successive intervals among $I_1,\ldots,I_m$. Each has size less than $2q$,
and successive intervals are separated by at most one $A_i$, of size
at most $b$. Hence $|\pi(x+1)-\pi(x)|\le4q+b$.
\end{enumerate}
In all cases the image distance is at most $4q+b$, so \Cref{cl:adjacent-gap} gives $\Lip(\pi)\le4q+b$.
\end{proof}

\subsection{Proof of the upper-bound theorem}

\begin{proof}[Proof of \Cref{thm:upper-bound}]
Work on the event where \Cref{prop:matching-patterns,lem:rearrangement}
both hold, which has probability $1-o(1)$. Choose $f,\mathcal I$
from \Cref{prop:matching-patterns}. For $I\subseteq[n]$, recall that $B(I)=\{j\in[M]:B_j\subseteq I\}$. Denote $U=\bigcup_{I\in\mathcal I}I$,
and let $\sigma:B(U)\to B(U)$ be the partial involution inducing $f$.
By \Cref{prop:reduction-to-permutations}, it suffices to extend
$\sigma$ to an involution $\pi$ fixing $M$, satisfying the neighbor
identities from \Cref{prop:reduction-to-permutations}, and having $\Lip(\pi)=O(D)$. Since $\pi$ is an involution, the inverse bound $\Lip(\pi^{-1})=O(D)$ holds as well.

Let $J_1,\ldots,J_m$ be the partition from \Cref{it:ub-matching-margins} of \Cref{prop:matching-patterns}. Each prescribed pair $I,\bar I$
is contained in a common part $J_j$.
We will apply \Cref{lem:oscillating-completion} to intervals of block
indices. Assign each block $B_k$ to the part containing its largest
index, and let
\[
 K_j=\{k\in[M]:\max B_k\in J_j\}.
\]
The sets $K_1,\ldots,K_m$ partition $[M]$ into intervals of block
indices. Since $|J_j|\le3H$, each has size at most $3H/p+2$.

If a prescribed pair $I,\bar I$ lies in $J_j$, every block forming
either interval lies entirely in $J_j$. Thus $B(I),B(\bar I)\subseteq K_j$,
and $\sigma$ maps $B(I)$ increasingly onto $B(\bar I)$.
Consequently, $\sigma$ restricts to an involution on $B(U)\cap K_j$ and, below, we aim to show that $\sigma|_{B(U)\cap K_j}$ can be extended to an involution $\pi_j$ fixing the end-points of $K_j$ and satisfying $\Lip(\pi_j) = O(D)$.

Fix a part $K_j$. We verify the remaining
hypotheses of \Cref{lem:oscillating-completion}.
By \Cref{it:ub-matching-separation}, each prescribed
interval $I$ has size $O(LD)$. Since $I$ is a union of blocks of size $p$,
we have $|B(I)|=|I|/p=O(LD/p)$.
Choose a positive integer $b=O(LD/p)=O(D)$ bounding the size $B(I)$ for all $I\in \mathcal \mathcal I$ from above.

Next, we bound the sizes of the free intervals in $K_j$. Successive prescribed intervals of  $\mathcal{I}$ contained in $J_j$ are separated by at least $D/2$ points by \Cref{it:ub-matching-separation} of \Cref{prop:matching-patterns}. Then, \Cref{it:ub-matching-margins} of that proposition gives the same margin between the first
prescribed interval and the left endpoint of $J_j$, and between the last
prescribed interval and the right endpoint of $J_j$.
In passing to block indices, we allow for a loss of at most two blocks
from rounding at the boundaries.
Thus each maximal interval of free indices in $K_j$ has size at least
$D/(2p)-2$ .

We choose $q=\lfloor D/(4p)\rfloor$, leaving room for this rounding loss.
Since $D/p=A\lceil\sqrt H\rceil$ tends to infinity, for sufficiently
large $n$ we have $D/(2p)-2\ge D/(4p)\ge q$.
This verifies the lower bound on the sizes of the free intervals required
by \Cref{lem:oscillating-completion}.

Finally, we verify that $q^2\ge16|K_j|$.
Recall that $|K_j|\le3H/p+2\le4H$ for sufficiently large $n$, so it
suffices to show that $q\ge8\sqrt H$.
Substituting $D=Ap\lceil\sqrt H\rceil$ into the definition of $q$ gives
\[
 q=\left\lfloor\frac{A\lceil\sqrt H\rceil}{4}\right\rfloor
 \ge\frac{A\sqrt H}{4}-1.
\]
Taking $A\ge64$, the right-hand side is at least $8\sqrt H$, as required.

Apply \Cref{lem:oscillating-completion} to each part $K_j$, taking all prescribed
intervals $B(I), B(\bar I) \subseteq K_j$ as the sets $A_1, \ldots, A_h$, to obtain $\pi_j:K_j \to K_j$ with
$\Lip(\pi_j)\le4q+b=O(D)$.

Define $\pi$ by $\pi|_{K_j}=\pi_j$ for every $j$ and let us bound $\Lip(\pi)$.
As the local maps $\pi|_{K_j}$ are involutions, $\pi$ is an involution of $[M]$ extending $\sigma$. Also, the map $\pi$ fixes all endpoints of all parts $K_j$, including $M$. Since each prescribed block interval $B(I)$ and its neighbors lie in one part $K_j$, $\pi$ also satisfies the neighbor identities from \Cref{it:neighbour_identitiy} of \Cref{prop:reduction-to-permutations}, as $\pi_j$ satisfies them as well. Due to \Cref{cl:adjacent-gap}, in order to bound the Lipschitz constants of $\pi$ we need to bound the quantity $\frac{|\pi(x) - \pi(y)|}{|x - y|}$ for consecutive indices $x, y \in [M]$. If consecutive indices $x,y \in [M]$ lie within one part $K_j$ then the quantity is bounded from above by $\Lip(\pi_j)=O(D)$. Otherwise, $x$ and $y$ lie in different parts then they both lie at the boundaries of their parts and, as a consequence $\pi(x) = x$ and $\pi(y)=y$, making the quantity equal $1 = O(D)$.
Thus \Cref{cl:adjacent-gap} gives $\Lip(\pi)= O(D)$, and, as we noted at the beginning of the proof, the same bound holds for $\pi^{-1}=\pi$.

Given $\pi$, \Cref{prop:reduction-to-permutations} gives a bijection $F:X\to Y$ with
\[
 c(X,Y)\le\dist(F)\le n^{2/3}\log^{O(1)}n.\qedhere
\]
\end{proof}

\bibliographystyle{plainnat}
\bibliography{references}
\appendix
\section{Woven permutation minors}\label{sec:woven}

The goal of this appendix is to show that an oscillating structure is necessary
for low-distortion extensions of a given partial mapping with large displacement.
For simplicity, we focus on permutations $f\in S_n$ with a common Lipschitz
bound $k$ in both directions and displacement $\Omega(k^2)$. In this regime,
the distortion $\Lip(f)\Lip(f^{-1})$ and the displacement
$\operatorname{disp}(f)$ are within a constant factor of each other.

\begin{definition}[Minors]\label{def:minor}
Let $f\in S_n$ and let $S=\{x_1<\cdots<x_N\}\subseteq[n]$ be nonempty.
Write $f(S)=\{y_1<\cdots<y_N\}$.
The \emph{minor of $f$ induced by $S$} is the permutation
$\theta\in S_N$ determined by
\[
 f(x_i)=y_{\theta(i)}\qquad(i\in[N]).
\]
Thus $\theta(i)=j$ if $f$ maps the $i$th smallest element of $S$
to the $j$th smallest element of $f(S)$.

A permutation $\theta\in S_N$ is a \emph{box minor of $f$} if it is induced
by a nonempty set $S=\{x\in U:f(x)\in V\}$ for integer intervals
$U,V\subseteq[n]$. Equivalently, we retain exactly the points
$(x,f(x))$ in the rectangle $U\times V$ and replace each coordinate
by its rank among the retained coordinates.
\end{definition}

We now define a family of oscillating permutations, called
\emph{woven permutations}. We show that every permutation $f\in S_n$
with $\Lip(f),\Lip(f^{-1})\le k$ and displacement $\Omega(k^2)$
contains a specific kind of box minor, of order $\Theta(k^2)$, which we call a woven permutation.

\begin{definition}[Coordinate orders]\label{def:coordinate-orders}
Let $P\subseteq\mathbb R^2$ be a finite nonempty set with distinct
$x$-coordinates and distinct $y$-coordinates. Its \emph{$x$-order} is the order on $P$ induced by increasing $x$ values. Similarly, define the
\emph{$y$-order} by increasing $y$ values. The \emph{permutation induced by $P$}
sends the position of each point in $x$-order to its position in $y$-order. 
\end{definition}

For example, if $P=\{(1,30),(4,10),(7,20)\}$, then the points
in $x$-order have $y$-ranks $3,1,2$, respectively. Thus $P$
induces the permutation $\theta\in S_3$ given by
$\theta(1)=3$, $\theta(2)=1$, and $\theta(3)=2$.

\begin{definition}[Woven permutations]\label{def:grid}
Let $P\subseteq\mathbb R^2$ be a finite nonempty set with distinct
$x$-coordinates and distinct $y$-coordinates, and let $m\ge1$ be an integer.
Partition a rectangle containing $P$ into $m$ vertical \emph{$x$-strips}
and $m$ horizontal \emph{$y$-strips}, with no point on a boundary of a strip.
The intersection of two strips is a \emph{cell} of the resulting $m\times m$ grid.
We call $P$ \emph{woven} with respect to this grid if it satisfies the following two properties:
\begin{enumerate}
 \item Every cell contains at least one point of $P$;
 \item Consecutive points in either the $x$-order or the $y$-order occupy either: (1)  the same cell, (2) cells sharing a side, or (3) cells sharing a corner.
\end{enumerate}
We call the permutation induced by $P$ a \emph{woven permutation}
on this grid.
\end{definition}

We say that the woven representation of $P$ has \emph{strip bound} $b\ge1$ if every $x$-strip and every $y$-strip contains at most $b$ points of $P$.

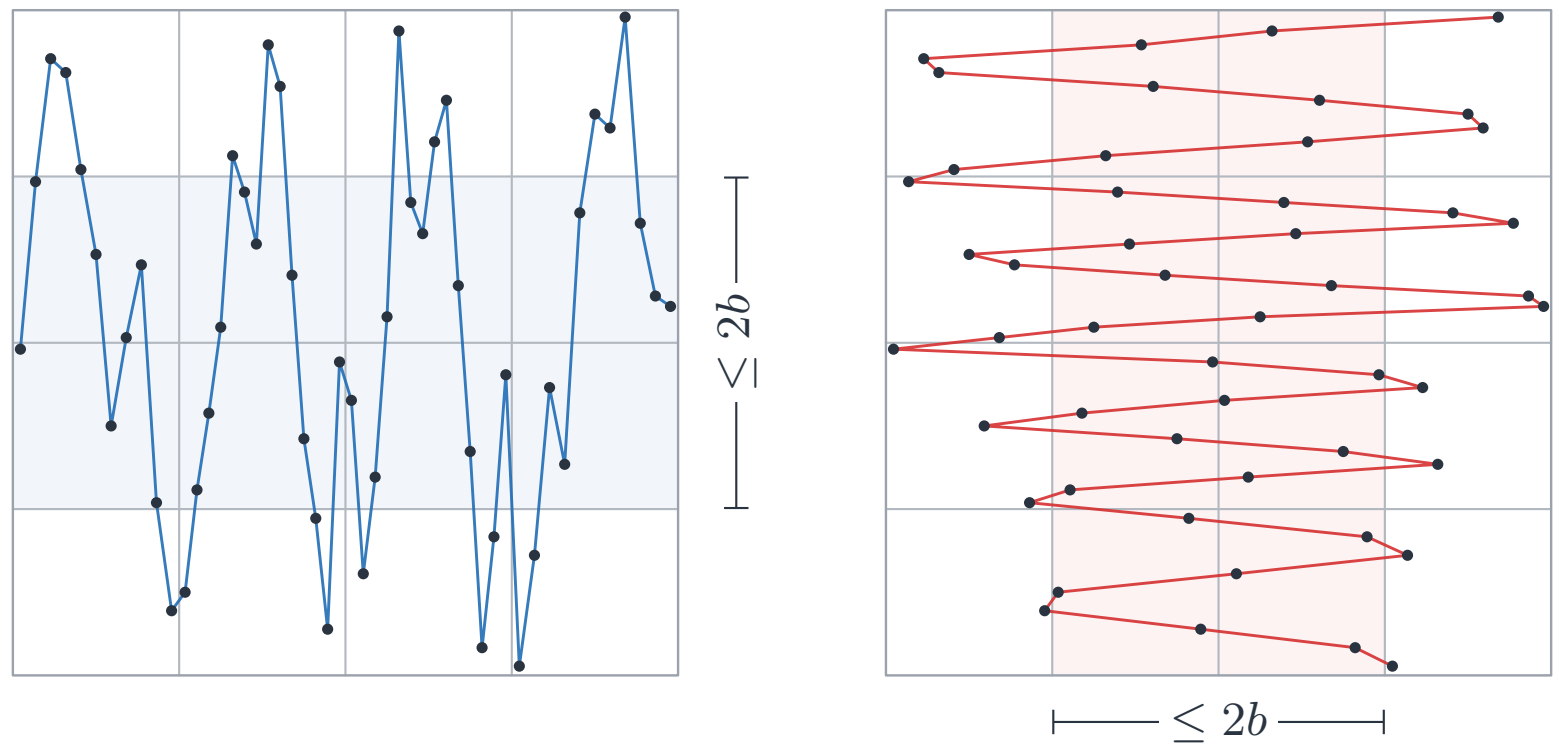
\begin{figure}[H]
\centering
\begin{tikzpicture}[x=1.1cm,y=1.1cm,line cap=round,line join=round]
\definecolor{pointInk}{HTML}{293440}
\definecolor{gridInk}{HTML}{89929C}
\definecolor{xInk}{HTML}{206DB5}
\definecolor{yInk}{HTML}{D62F2F}
\begin{scope}[xshift=0.00000cm]
\fill[xInk!6] (0,1) rectangle (4,3);
\draw[gridInk!65,line width=.4pt] (1,0) -- (1,4);
\draw[gridInk!65,line width=.4pt] (0,1) -- (4,1);
\draw[gridInk!65,line width=.4pt] (2,0) -- (2,4);
\draw[gridInk!65,line width=.4pt] (0,2) -- (4,2);
\draw[gridInk!65,line width=.4pt] (3,0) -- (3,4);
\draw[gridInk!65,line width=.4pt] (0,3) -- (4,3);
\draw[gridInk!85,line width=.45pt] (0,0) rectangle (4,4);
\draw[xInk,opacity=.9,line width=.55pt] (0.045455,1.961538) -- (0.136364,2.968750) -- (0.227273,3.708333) -- (0.318182,3.625000) -- (0.409091,3.041667) -- (0.500000,2.531250) -- (0.590909,1.500000) -- (0.681818,2.031250) -- (0.772727,2.468750) -- (0.863636,1.038462) -- (0.954545,0.388889) -- (1.035714,0.500000) -- (1.107143,1.115385) -- (1.178571,1.576923) -- (1.250000,2.093750) -- (1.321429,3.125000) -- (1.392857,2.906250) -- (1.464286,2.593750) -- (1.535714,3.791667) -- (1.607143,3.541667) -- (1.678571,2.406250) -- (1.750000,1.423077) -- (1.821429,0.944444) -- (1.892857,0.277778) -- (1.964286,1.884615) -- (2.035714,1.653846) -- (2.107143,0.611111) -- (2.178571,1.192308) -- (2.250000,2.156250) -- (2.321429,3.875000) -- (2.392857,2.843750) -- (2.464286,2.656250) -- (2.535714,3.208333) -- (2.607143,3.458333) -- (2.678571,2.343750) -- (2.750000,1.346154) -- (2.821429,0.166667) -- (2.892857,0.833333) -- (2.964286,1.807692) -- (3.045455,0.055556) -- (3.136364,0.722222) -- (3.227273,1.730769) -- (3.318182,1.269231) -- (3.409091,2.781250) -- (3.500000,3.375000) -- (3.590909,3.291667) -- (3.681818,3.958333) -- (3.772727,2.718750) -- (3.863636,2.281250) -- (3.954545,2.218750);
\fill[pointInk] (0.045455,1.961538) circle[radius=1.05pt];
\fill[pointInk] (0.136364,2.968750) circle[radius=1.05pt];
\fill[pointInk] (0.227273,3.708333) circle[radius=1.05pt];
\fill[pointInk] (0.318182,3.625000) circle[radius=1.05pt];
\fill[pointInk] (0.409091,3.041667) circle[radius=1.05pt];
\fill[pointInk] (0.500000,2.531250) circle[radius=1.05pt];
\fill[pointInk] (0.590909,1.500000) circle[radius=1.05pt];
\fill[pointInk] (0.681818,2.031250) circle[radius=1.05pt];
\fill[pointInk] (0.772727,2.468750) circle[radius=1.05pt];
\fill[pointInk] (0.863636,1.038462) circle[radius=1.05pt];
\fill[pointInk] (0.954545,0.388889) circle[radius=1.05pt];
\fill[pointInk] (1.035714,0.500000) circle[radius=1.05pt];
\fill[pointInk] (1.107143,1.115385) circle[radius=1.05pt];
\fill[pointInk] (1.178571,1.576923) circle[radius=1.05pt];
\fill[pointInk] (1.250000,2.093750) circle[radius=1.05pt];
\fill[pointInk] (1.321429,3.125000) circle[radius=1.05pt];
\fill[pointInk] (1.392857,2.906250) circle[radius=1.05pt];
\fill[pointInk] (1.464286,2.593750) circle[radius=1.05pt];
\fill[pointInk] (1.535714,3.791667) circle[radius=1.05pt];
\fill[pointInk] (1.607143,3.541667) circle[radius=1.05pt];
\fill[pointInk] (1.678571,2.406250) circle[radius=1.05pt];
\fill[pointInk] (1.750000,1.423077) circle[radius=1.05pt];
\fill[pointInk] (1.821429,0.944444) circle[radius=1.05pt];
\fill[pointInk] (1.892857,0.277778) circle[radius=1.05pt];
\fill[pointInk] (1.964286,1.884615) circle[radius=1.05pt];
\fill[pointInk] (2.035714,1.653846) circle[radius=1.05pt];
\fill[pointInk] (2.107143,0.611111) circle[radius=1.05pt];
\fill[pointInk] (2.178571,1.192308) circle[radius=1.05pt];
\fill[pointInk] (2.250000,2.156250) circle[radius=1.05pt];
\fill[pointInk] (2.321429,3.875000) circle[radius=1.05pt];
\fill[pointInk] (2.392857,2.843750) circle[radius=1.05pt];
\fill[pointInk] (2.464286,2.656250) circle[radius=1.05pt];
\fill[pointInk] (2.535714,3.208333) circle[radius=1.05pt];
\fill[pointInk] (2.607143,3.458333) circle[radius=1.05pt];
\fill[pointInk] (2.678571,2.343750) circle[radius=1.05pt];
\fill[pointInk] (2.750000,1.346154) circle[radius=1.05pt];
\fill[pointInk] (2.821429,0.166667) circle[radius=1.05pt];
\fill[pointInk] (2.892857,0.833333) circle[radius=1.05pt];
\fill[pointInk] (2.964286,1.807692) circle[radius=1.05pt];
\fill[pointInk] (3.045455,0.055556) circle[radius=1.05pt];
\fill[pointInk] (3.136364,0.722222) circle[radius=1.05pt];
\fill[pointInk] (3.227273,1.730769) circle[radius=1.05pt];
\fill[pointInk] (3.318182,1.269231) circle[radius=1.05pt];
\fill[pointInk] (3.409091,2.781250) circle[radius=1.05pt];
\fill[pointInk] (3.500000,3.375000) circle[radius=1.05pt];
\fill[pointInk] (3.590909,3.291667) circle[radius=1.05pt];
\fill[pointInk] (3.681818,3.958333) circle[radius=1.05pt];
\fill[pointInk] (3.772727,2.718750) circle[radius=1.05pt];
\fill[pointInk] (3.863636,2.281250) circle[radius=1.05pt];
\fill[pointInk] (3.954545,2.218750) circle[radius=1.05pt];
\draw[pointInk,line width=.4pt,|-|] (4.35,1) -- (4.35,3);
\node[font=\footnotesize,text=pointInk,rotate=90,fill=white,inner sep=2pt] at (4.35,2) {$\le 2b$};
\end{scope}
\begin{scope}[xshift=5.77500cm]
\fill[yInk!6] (1,0) rectangle (3,4);
\draw[gridInk!65,line width=.4pt] (1,0) -- (1,4);
\draw[gridInk!65,line width=.4pt] (0,1) -- (4,1);
\draw[gridInk!65,line width=.4pt] (2,0) -- (2,4);
\draw[gridInk!65,line width=.4pt] (0,2) -- (4,2);
\draw[gridInk!65,line width=.4pt] (3,0) -- (3,4);
\draw[gridInk!65,line width=.4pt] (0,3) -- (4,3);
\draw[gridInk!85,line width=.45pt] (0,0) rectangle (4,4);
\draw[yInk,opacity=.9,line width=.55pt] (3.045455,0.055556) -- (2.821429,0.166667) -- (1.892857,0.277778) -- (0.954545,0.388889) -- (1.035714,0.500000) -- (2.107143,0.611111) -- (3.136364,0.722222) -- (2.892857,0.833333) -- (1.821429,0.944444) -- (0.863636,1.038462) -- (1.107143,1.115385) -- (2.178571,1.192308) -- (3.318182,1.269231) -- (2.750000,1.346154) -- (1.750000,1.423077) -- (0.590909,1.500000) -- (1.178571,1.576923) -- (2.035714,1.653846) -- (3.227273,1.730769) -- (2.964286,1.807692) -- (1.964286,1.884615) -- (0.045455,1.961538) -- (0.681818,2.031250) -- (1.250000,2.093750) -- (2.250000,2.156250) -- (3.954545,2.218750) -- (3.863636,2.281250) -- (2.678571,2.343750) -- (1.678571,2.406250) -- (0.772727,2.468750) -- (0.500000,2.531250) -- (1.464286,2.593750) -- (2.464286,2.656250) -- (3.772727,2.718750) -- (3.409091,2.781250) -- (2.392857,2.843750) -- (1.392857,2.906250) -- (0.136364,2.968750) -- (0.409091,3.041667) -- (1.321429,3.125000) -- (2.535714,3.208333) -- (3.590909,3.291667) -- (3.500000,3.375000) -- (2.607143,3.458333) -- (1.607143,3.541667) -- (0.318182,3.625000) -- (0.227273,3.708333) -- (1.535714,3.791667) -- (2.321429,3.875000) -- (3.681818,3.958333);
\fill[pointInk] (0.045455,1.961538) circle[radius=1.05pt];
\fill[pointInk] (0.136364,2.968750) circle[radius=1.05pt];
\fill[pointInk] (0.227273,3.708333) circle[radius=1.05pt];
\fill[pointInk] (0.318182,3.625000) circle[radius=1.05pt];
\fill[pointInk] (0.409091,3.041667) circle[radius=1.05pt];
\fill[pointInk] (0.500000,2.531250) circle[radius=1.05pt];
\fill[pointInk] (0.590909,1.500000) circle[radius=1.05pt];
\fill[pointInk] (0.681818,2.031250) circle[radius=1.05pt];
\fill[pointInk] (0.772727,2.468750) circle[radius=1.05pt];
\fill[pointInk] (0.863636,1.038462) circle[radius=1.05pt];
\fill[pointInk] (0.954545,0.388889) circle[radius=1.05pt];
\fill[pointInk] (1.035714,0.500000) circle[radius=1.05pt];
\fill[pointInk] (1.107143,1.115385) circle[radius=1.05pt];
\fill[pointInk] (1.178571,1.576923) circle[radius=1.05pt];
\fill[pointInk] (1.250000,2.093750) circle[radius=1.05pt];
\fill[pointInk] (1.321429,3.125000) circle[radius=1.05pt];
\fill[pointInk] (1.392857,2.906250) circle[radius=1.05pt];
\fill[pointInk] (1.464286,2.593750) circle[radius=1.05pt];
\fill[pointInk] (1.535714,3.791667) circle[radius=1.05pt];
\fill[pointInk] (1.607143,3.541667) circle[radius=1.05pt];
\fill[pointInk] (1.678571,2.406250) circle[radius=1.05pt];
\fill[pointInk] (1.750000,1.423077) circle[radius=1.05pt];
\fill[pointInk] (1.821429,0.944444) circle[radius=1.05pt];
\fill[pointInk] (1.892857,0.277778) circle[radius=1.05pt];
\fill[pointInk] (1.964286,1.884615) circle[radius=1.05pt];
\fill[pointInk] (2.035714,1.653846) circle[radius=1.05pt];
\fill[pointInk] (2.107143,0.611111) circle[radius=1.05pt];
\fill[pointInk] (2.178571,1.192308) circle[radius=1.05pt];
\fill[pointInk] (2.250000,2.156250) circle[radius=1.05pt];
\fill[pointInk] (2.321429,3.875000) circle[radius=1.05pt];
\fill[pointInk] (2.392857,2.843750) circle[radius=1.05pt];
\fill[pointInk] (2.464286,2.656250) circle[radius=1.05pt];
\fill[pointInk] (2.535714,3.208333) circle[radius=1.05pt];
\fill[pointInk] (2.607143,3.458333) circle[radius=1.05pt];
\fill[pointInk] (2.678571,2.343750) circle[radius=1.05pt];
\fill[pointInk] (2.750000,1.346154) circle[radius=1.05pt];
\fill[pointInk] (2.821429,0.166667) circle[radius=1.05pt];
\fill[pointInk] (2.892857,0.833333) circle[radius=1.05pt];
\fill[pointInk] (2.964286,1.807692) circle[radius=1.05pt];
\fill[pointInk] (3.045455,0.055556) circle[radius=1.05pt];
\fill[pointInk] (3.136364,0.722222) circle[radius=1.05pt];
\fill[pointInk] (3.227273,1.730769) circle[radius=1.05pt];
\fill[pointInk] (3.318182,1.269231) circle[radius=1.05pt];
\fill[pointInk] (3.409091,2.781250) circle[radius=1.05pt];
\fill[pointInk] (3.500000,3.375000) circle[radius=1.05pt];
\fill[pointInk] (3.590909,3.291667) circle[radius=1.05pt];
\fill[pointInk] (3.681818,3.958333) circle[radius=1.05pt];
\fill[pointInk] (3.772727,2.718750) circle[radius=1.05pt];
\fill[pointInk] (3.863636,2.281250) circle[radius=1.05pt];
\fill[pointInk] (3.954545,2.218750) circle[radius=1.05pt];
\draw[pointInk,line width=.4pt,|-|] (1,-.28) -- (3,-.28);
\node[font=\footnotesize,text=pointInk,fill=white,inner sep=2pt] at (2,-.28) {$\le 2b$};
\end{scope}
\end{tikzpicture}
\caption{An illustration of a set $P$ which induces a woven permutation. The $x$-order is depicted on the left, and the $y$-order appears on the right.}
\label{fig:filled-grid}
\end{figure}

\begin{observation}\label{obs:strip-bound}
Every woven permutation $\theta\in S_N$ on an $m\times m$ grid
with strip bound $b$ satisfies
\[
 m^2\le N\le mb,\qquad
 \Lip(\theta),\Lip(\theta^{-1})\le2b-1.
\]
\end{observation}
\begin{proof}
Since every cell is occupied, $N\ge m^2$. Each of the
$m$ $x$-strips contains at most $b$ points, so $N\le mb$. 

Let
$p_1,\ldots,p_N$ be the points in $x$-order. For $1\le i<N$, the points
$p_i$ and $p_{i+1}$ lie either in the same or in neighboring $y$-strips and these strips contain all the points that lie between $p_i$ and $p_{i+1}$ in $y$-order. As any set of at most two strips contains at most $2b$ points, we conclude $|\theta(i+1)-\theta(i)|\le2b-1$ and, thus, $\Lip(\theta)\le2b-1$. By exchanging coordinates we get
$\Lip(\theta^{-1})\le2b-1$.
\end{proof}

The following theorem shows that large displacement, together with
bounds on both directional Lipschitz constants, forces a large woven box minor.

\begin{theorem}\label{prop:main}
For every $\varepsilon>0$, there exists $c=c(\varepsilon)>0$ such that,
for all sufficiently large integers $k$, every $f\in S_n$ satisfying
\[
 \max\{\Lip(f),\Lip(f^{-1})\}\le k,
 \qquad \operatorname{disp}(f)\ge\varepsilon k^2
\]
has a box minor $\theta\in S_N$ that is a woven permutation
on an $m\times m$ grid
with at most $k$ points per strip, where
\[
 ck\le m\le k,\qquad m^2\le N\le mk.
\]
\end{theorem}

In particular, the representing woven point set has $\Theta_\varepsilon(k^2)$ points and $O_\varepsilon(m)$ points per strip.

\subsection{Proof of \texorpdfstring{\cref{prop:main}}{Theorem \getrefnumber{prop:main}}}

In this subsection, for integers $a\le b$ let us denote $[a,b]=\{a,a+1,\ldots,b\}$.
We begin by proving the following lemma.

\begin{lemma}\label{lem:windows}
Let $f\in S_n$, let $1\le k<n$ be an integer, and suppose
$\Lip(f^{-1})\le k$. Then
\[
 \operatorname{disp}(f)\le2\max_{1\le p\le n-k}\operatorname{diam} f([p,p+k]).
\]
\end{lemma}

\begin{proof}
Denote $W=\max_{1\le p\le n-k}\operatorname{diam} f([p,p+k])$.
We claim that for all $a \leq b \leq c$ from $[1, n]$ it holds that
\begin{equation}
    \min\{f(a),f(c)\}-W\le f(b)\le\max\{f(a),f(c)\}+W.
    \label{eq:f(b)_close_to_f(a)_and_f(c)}
\end{equation}
Wlog assume that $f(a) \leq f(c)$ and denote $A = (f(a), f(a)+1, \ldots, f(c))$. As $\Lip(f^{-1}) \leq k$, the nearest coordinates in the tuple $f^{-1}(A)$ differ by at most $k$. Since the first and the last elements of $f^{-1}(A)$ are $a$ and $c$ respectively and $a \leq b \leq c$, there exists $d \in f^{-1}(A)$ satisfying $|d - b| \leq k$. By the definition of $W$ we have $|f(d) - f(b)| \leq W$. Also, by the definition of the set $A$ we have $f(a)\leq f(d) \leq f(c)$. Combining the inequalities we conclude~\eqref{eq:f(b)_close_to_f(a)_and_f(c)}.

Wlog, we assume $f(1)<f(n)$, as we can substitute $f(i)$ by $f(n+1 -i)$ for all $i \in [n]$. Applying~\eqref{eq:f(b)_close_to_f(a)_and_f(c)} for $a=1$ and $c=n$ we conclude that $f(i) \geq f(1)-W$ for all $i \in [n]$. As $f(i) = 1$ for some $i \in[n]$, we have
$f(1)\le W+1$. Now, applying~\eqref{eq:f(b)_close_to_f(a)_and_f(c)} for $a = 1$ and for all $1 \leq b \leq c \leq n$ we get
\begin{equation}
\label{eq:f(b)_better_bound}
 f(b)\le\max\{f(1),f(c)\}+W\le f(c)+2W.
\end{equation}
Consider arbitrary $i \in [n]$. Notice that $f([1, i])$ contains an element that is at least $i$. Due to~\eqref{eq:f(b)_better_bound} applied to $c = i$, all elements in $f([1, i])$ are at most $f(i)+2W$. Hence, $f(i)\ge i-2W$. Similarly, since $f([i,n])$ contains $n-i+1$ distinct integers,
there exists $c\in[i,n]$ with $f(c)\le i$.
Applying~\eqref{eq:f(b)_better_bound} with $b=i$ gives
$f(i)\le f(c)+2W\le i+2W$.
Together with the preceding lower bound, this yields
$|f(i)-i|\le2W$ for every $i\in[n]$, and hence
$\operatorname{disp}(f)\le2W$, as required.
\end{proof}

\begin{lemma}\label{lem:filled-rectangle}
Let $f\in S_n$, and let $1\le k<n$ and $m\ge2$ be integers.
Suppose that $\Lip(f) \le k$ and $\Lip(f^{-1}) \le k$ and
$\operatorname{disp}(f)\ge128mk$. Then there exist intervals $U,V\subseteq[n]$, each of size $mk$,
such that, writing $U_1,\ldots,U_m$ and $V_1,\ldots,V_m$
for their respective partitions into consecutive blocks of
$k$ integers, we have
\[
 f(U_r)\cap V_s\ne\varnothing\qquad(r,s\in[m]).
\]
\end{lemma}

\begin{proof}
Set $R=mk$. By \cref{lem:windows}, choose $I=[a,a+k]$ with
$\operatorname{diam}f(I)\ge\operatorname{disp}(f)/2\ge64R$,
and set $J=[\min f(I),\max f(I)]$.
We first construct the sets $U$ and $V$ and then verify that every intersection
$f(U_r)\cap V_s$ is non-empty.


\paragraph{Step 1: finding $U$ and $V$.}
We claim that, after possibly substituting $f(i)$ by $f(n+1 -i)$ for all $i \in [n]$, the set
$J_-:=\{y\in J:f^{-1}(y)<a-R\}$ satisfies $|J_-|\ge2R+1$.
Indeed, the number of elements in $f^{-1}(J) \subseteq [n]$ that lie within distance
$R$ of $I$ is at most $|I| + 2R \leq  k+1+2R\le3R$. Of the remaining values $f^{-1}(J)$, at least
$(|J|-3R)/2\ge2R+1$ lie on the same side of $I$. So, wlog, we can assume that $|J_-|\ge2R+1$. 

Let $y_0:=\min J_-$ and $y_1:=\max J_-$, and set
\[
 U=[a-R,a-1],\qquad V=[y_0+1,y_0+R].
\]
Since $J_-\ne\varnothing$, we have $U\subseteq[n]$. Since
$y_1-y_0\ge|J_-|-1\ge2R$, we have $V\subseteq J\subseteq[n]$.
Both $U$ and $V$ have order $R$, as desired. Denote their consecutive disjoint $k$-point
blocks by $U_1,\ldots,U_m$ and $V_1,\ldots,V_m$, respectively. Wlog, we can assume that for all $1 \leq i < j \leq m$ all the elements of $U_i$ are smaller than all the elements of $U_j$ and all the elements of $V_i$ are smaller than all the elements of $V_j$.

\paragraph{Step 2: verifying that every cell is occupied.}
We claim that every $U_r$ contains points $x_r^-,x_r^+$ satisfying
\[
 f(x_r^-)\le y_0,\qquad f(x_r^+)\ge y_1-k.
\]
Fix $i\in\{0,1\}$. Let us first find
$z_i\in f(I)\cap[y_i-k,y_i]$. Let $r_1, r_2 \in I$ be the elements satisfying $f(r_1) = \min f(I)$ and $f(r_2) = \max f(I)$. Wlog, assume that $r_1 < r_2$, in the case $r_1 > r_2$ the proof is analogous. Consider the sequence $f(r_1)$, $f(r_1+1), \ldots,$ $f(r_2)$. As $\Lip (f) \leq k$, the consecutive elements in this sequence differ by at most $k$.
So, as $y_i\in[f(r_1),f(r_2)]$, the sequence contains an element $z_i \in [y_i - k, y_i]$ and, therefore, $f^{-1}(z_i) \in [r_1, r_2] \subseteq I$, as needed.

Now, consider another sequence
$f^{-1}(z_i)$, $f^{-1}(z_i+1)$, $\ldots,$ $f^{-1}(y_i)$.
Its first term belongs to $I=[a,a+k]$ and therefore is larger than all elements of $U$.
Its last term is smaller than $a-R=\min U$, since $y_i\in J_-$.
Consecutive terms in the sequence differ by at most $k$, by $\Lip(f^{-1})\le k$.
Hence the sequence has an element in every $k$-point block $U_r$ for all $r \in [m]$ as $f^{-1}(z_i) \geq \max U_r$ and $f^{-1}(y_i) \leq \min U_r$.
Each such an element $x$ satisfies $z_i\le f(x)\le y_i$.
For $i=0$, choose such a point as $x_r^-$; then $f(x_r^-)\le y_0$.
For $i=1$, choose such a point as $x_r^+$; then
$f(x_r^+)\ge z_1\ge y_1-k$.

Fix $r\in[m]$. Since $y_1-y_0\ge2R$ and $R\ge k$, we have
\[
 f(x_r^-)\le y_0<\min V,
 \qquad
 f(x_r^+)\ge y_1-k\ge y_0+R=\max V.
\]
Consider the sequence $f(x_r^-), f(x_r^-+1), \ldots, f(x_r^+) \subseteq f(U_r)$. It is sufficient for us to show that this sequence contains at least one element in every $V_s$ for all $s \in [m]$.
Every two consecutive elements in this sequence differ by at most $k$ as $\Lip(f)\leq k$. The first element of the sequence is smaller than $\min V$ and the last element is larger than $\max V$. Therefore, as $|V_s| = k$ for all $s \in [m]$, the sequence contains at least one element in every set $V_s$ for all $s \in [m]$, as needed.
\end{proof}

\begin{lemma}\label{lem:rectangle}
Let $f\in S_n$ and let $\theta\in S_N$ be its box minor induced by
intervals $U,V\subseteq[n]$. Then
\[
 \Lip(\theta)\le\Lip(f),\qquad
 \Lip(\theta^{-1})\le\Lip(f^{-1}).
\]
Also, for integers $k,m\ge1$, suppose that $\max\{\Lip(f),\Lip(f^{-1})\}\le k$, suppose that 
$U$ and $V$ are partitioned into consecutive disjoint blocks $U_1,\ldots,U_m$ and $V_1,\ldots,V_m$ of size $k$ respectively, and suppose that
$f(U_r)\cap V_s\ne\varnothing$ for all $r,s\in[m]$.
Then $\theta$ is woven on an $m\times m$ grid with strip bound $k$.
\end{lemma}

\begin{proof}
We first prove that $\Lip(\theta)\le\Lip(f)$ and
$\Lip(\theta^{-1})\le\Lip(f^{-1})$.
Let $a<b$ be consecutive elements of $U\cap f^{-1}(V)$. As $U$ is an interval, we have $[a,b]\subseteq U$. Since $V$ is an interval containing $f(a)$ and $f(b)$, it contains
$[\min\{f(a),f(b)\},\max\{f(a),f(b)\}]$. Since $a$ and $b$ are consecutive in $U\cap f^{-1}(V)$, we also have $f([a, b]) \cap V = \{f(a), f(b)\}$. In other words, no $c$ with $a<c<b$ has its image $f(c)$ in $V$.

Assume $f(a)<f(b)$. Consider the smallest $t \in [a, b-1]$ such that $f(t+1) \geq f(b)$. Then, $f(t)\le f(a)<f(b)\le f(t+1)$. Thus,
\begin{equation}
    |f(b)-f(a)|\le\Lip(f).
    \label{eq:lip_woven}
\end{equation}
This bound is proved the same way in the case $f(b)<f(a)$.

The points $a$ and $b$ correspond to consecutive indices
$i$ and $i+1$ in the domain of $\theta$.
Replacing every element of $U \cap f^{-1}(V)$ by the corresponding element in the domain of $\theta$ does not
increase the distance between elements. Hence,
\[
 |\theta(i+1)-\theta(i)|
 \le |f(b)-f(a)|
 \le \Lip(f).
\]
Since $a<b$ were arbitrary consecutive elements of $U \cap f^{-1}(V)$, we conclude
$\Lip(\theta)\le\Lip(f)$.
Applying the same argument to $f^{-1}$ and the rectangle
$V\times U$ gives $\Lip(\theta^{-1})\le\Lip(f^{-1})$.

Now, assume that $\max\{\Lip(f),\Lip(f^{-1})\}\le k$ and the consecutive
$k$-point blocks of $U$ and $V$ satisfy $f(U_i)\cap V_j\ne\varnothing$
for all $i,j\in[m]$. Wlog, we can assume that for all $1 \leq i < j \leq m$ all the elements of $U_i$ are smaller than all the elements of $U_j$ and all the elements of $V_i$ are smaller than all the elements of $V_j$. We claim that, for all $r,s ,r', s' \in [m]$, if there exist points $(a, f(a)) \in U_r\times V_s$ and $(b, f(b)) \in U_{r'}\times V_{s'}$ that are consecutive in either $x$-order or $y$-order then
\[
    |r-r'|\le1,\qquad |s-s'|\le1.
\]
Indeed, consider such consecutive points $(a,f(a))$ and $(b,f(b))$
in $x$-order, with $a<b$. If $r'\ge r+2$, choose
$c\in U_{r+1}\cap f^{-1}(V)$; then $a<c<b$, contradicting consecutiveness.
Hence $0\le r'-r\le1$. If $|s-s'|\ge2$, the $k$-point blocks
$V_s$ and $V_{s'}$ are separated by at least $k$ integers, giving
\[
 k+1\le |f(b)-f(a)|\stackrel{\eqref{eq:lip_woven}}{\le}\Lip(f)\le k,
\]
a contradiction. Thus $|s-s'|\le1$ as well, so the cells coincide or
share a side or corner. For the $y$-order the claim is proved the same way.

Finally, we verify that every cell is occupied and every strip contains
at most $k$ points. For all $r, s \in [m]$ the set $f(U_r) \cap V_s$ is non-empty and thus $\theta$ is a woven permutation. Each $U_r$ and each $V_s$ have size $k$.
Hence every strip contains at most $k$ points, so the woven representation of $\theta$ has strip bound $k$.
\end{proof}

\begin{proof}[Proof of \cref{prop:main}]
Set $m=\min\{k,\lfloor\varepsilon k/128\rfloor\}$ and
$c=\min\{1,\varepsilon/256\}$.
For sufficiently large $k$, we have $k<n$, $2\le m\le k$, and $m\ge ck$.
Since $\operatorname{disp}(f)\ge\varepsilon k^2\ge128mk$,
\cref{lem:filled-rectangle,lem:rectangle} give the required woven box minor.
Finally, \cref{obs:strip-bound} gives $m^2\le N\le mk$.
\end{proof}

\section{Infinite expectation}\label{sec:inf-exp}

\begin{claim}[Infinite expectation]
\label{clm:infinite-expectation}
For every $n\ge 3$ and iid $X,Y\sim\mathcal{O}$, we have
\[\mathbb{E}[c(X,Y)]=\infty.\]
\end{claim}

 \begin{proof}
By \Cref{prop:exponential-spacings}, we may write $\Delta_i^Z=E_i^Z/\Sigma^Z$ for $Z\in\{X,Y\}$ and every $i\in[n-1]$, where $E_0^X,\ldots,E_n^X,E_0^Y,\ldots,E_n^Y$ are iid exponential with mean $1$ and $\Sigma^Z=\sum_{j=0}^n E_j^Z$. In particular, $z_n-z_1=\sum_{i=1}^{n-1}E_i^Z/\Sigma^Z$.

Let $f:X\to Y$ be a bijection. The images of a closest pair in $X$ are distinct points of $Y$, so $\Lip(f)\ge\Delta_{\min}^Y/\Delta_{\min}^X$. Since $|f(x_n)-f(x_1)|\le y_n-y_1$, we also have $\Lip(f^{-1})\ge(x_n-x_1)/(y_n-y_1)$. Since $f$ is arbitrary, and the normalization factors $\Sigma^X,\Sigma^Y$ cancel in the expression below, we get,
\[
 c(X,Y)
 \ge\frac{\Delta_{\min}^Y}{y_n-y_1}\cdot\frac{x_n-x_1}{\Delta_{\min}^X}
 =\frac{\min_{i\in[n-1]}E_i^Y}{\sum_{i=1}^{n-1}E_i^Y}\cdot\frac{\sum_{i=1}^{n-1}E_i^X}{\min_{i\in[n-1]}E_i^X}.
\]
\noindent For $s\ge1$, consider the two events
\[
 A_s=\bigl\{E_1^X\le 1/s,\ E_2^X\ge1\bigr\},
 \qquad
 B=\bigl\{1\le E_i^Y\le 2 \text{ for all } i\in[n-1]\bigr\},
\]
i.e., the first internal exponential variable of $X$ is at most $1/s$ and the second is at least $1$, while all internal exponential variables of $Y$ lie in $[1,2]$. Whenever both events hold, we have $\min_i E_i^X\le1/s$, $\sum_i E_i^X\ge1$, $\min_i E_i^Y\ge1$ and $\sum_i E_i^Y\le2n$, and hence $c(X,Y)\ge s/(2n)$.

By independence,
\[
 \PP(A_s\cap B)=\bigl(1-e^{-1/s}\bigr)e^{-1}\bigl(e^{-1}-e^{-2}\bigr)^{n-1}\ge\frac{a_n}{s}
\]
for some constant $a_n>0$. Thus $\PP(c(X,Y)\ge u)\ge a_n/(2nu)$ for all $u\ge1$, and
\[
 \mathbb{E}[c(X,Y)]\ge\int_1^\infty\PP\bigl(c(X,Y)\ge u\bigr)\,du=\infty.
\]
\end{proof}
\section{Oscillating permutations}
\label{sec:overview-oscilating-verification}

We now describe the oscillation construction appearing in the technical overview in more detail, and give some more intuition on how to define it for general $k$. Let $S_{i \rightarrow j}$ be the set of permutations sending $i$ to $j$.

\begin{proposition}
    Let $k > 0$ be odd and let $k\leq i < j \leq n-k$ satisfy $j - i = (k-1)^2$. Then, there exists $\pi \in S_{i \rightarrow j}$ satisfying
    $$\operatorname{dist}(\pi) = (2k-1)^2,\quad\quad\operatorname{disp}(\pi) = (k-1)^2.$$
\end{proposition}

\begin{figure}[h]
    \centering
    \begin{tikzpicture}[x=0.48cm,y=1cm]
        \foreach \r in {0,...,24} {
            \coordinate (T\r) at (\r,3);
            \coordinate (B\r) at (\r,0);
        }

        \foreach \r/\q in {0/0,1/9,2/10,3/19,4/20,
            5/21,6/18,7/11,8/8,9/1,
            10/2,11/7,12/12,13/17,14/22,
            15/23,16/16,17/13,18/6,19/3,
            20/4,21/5,22/14,23/15,24/24} {

            \pgfmathsetmacro{\redcomponent}{1-\r/24}
            \pgfmathsetmacro{\bluecomponent}{\r/24}
            \definecolor{arrowcolor}{rgb}
                {\redcomponent,0,\bluecomponent}

            \draw[
                arrowcolor,
                line width=0.7pt,
                ->,
                shorten <=2.5pt,
                shorten >=2.5pt
            ] (T\r) -- (B\q);
        }

        \foreach \r in {0,...,24} {
            \fill (T\r) circle (1.45pt);
            \fill (B\r) circle (1.45pt);
            \node[above=2pt,font=\scriptsize] at (T\r) {$\r$};
            \node[below=2pt,font=\scriptsize] at (B\r) {$\r$};
        }
    \end{tikzpicture}

    \caption{The permutation induced by matching the LR and UD snakes according to their common grid positions. Its distortion is at most $(2k-1)^2$. Colors interpolate from red to blue.}
    \label{fig:sigma_from_snake}
    \label{fig:snake-permutation}
\end{figure}
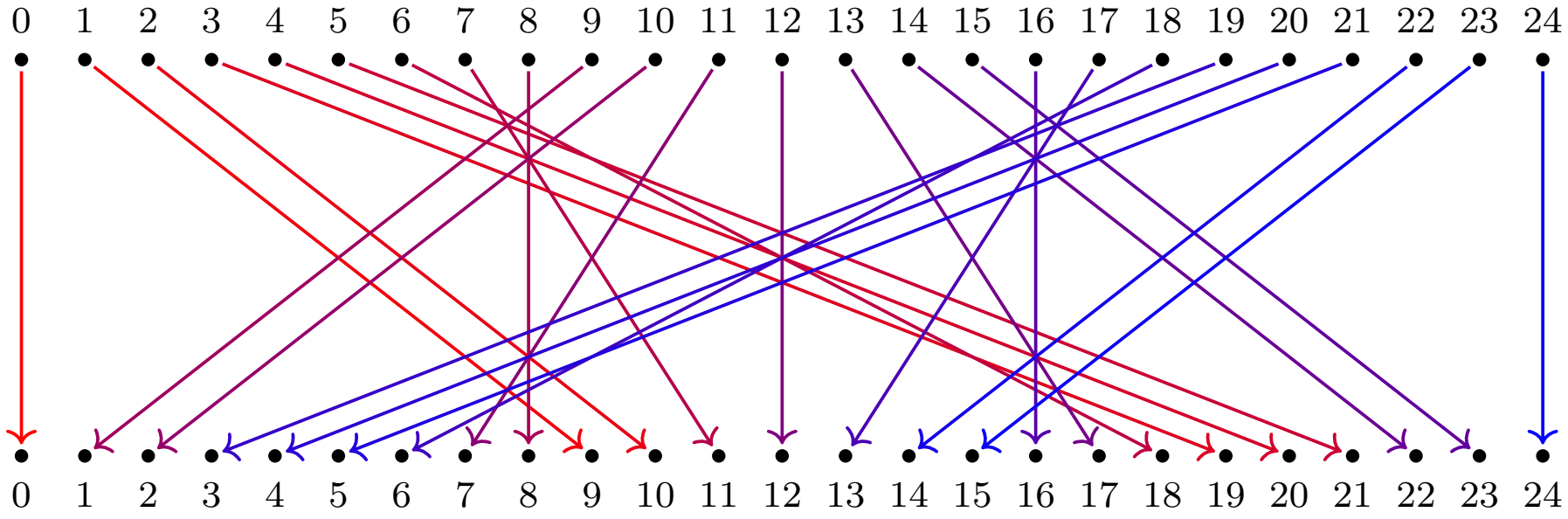

\begin{proof}[Sketch of the proof]
    
    Let us construct $\pi \in S_{i \rightarrow j}$. For all $r \in [n]$ outside of $[i-k+1, j+k-1]$ we let $\pi$ to be the identity, i.e. $\pi(r) = r$. Now we need to define the remaining permutation $\pi|_{[i-k+1, j+k-1]}$. In order to do that, consider the following specific permutation $\sigma$ of the set $[0, k^2-1]$ that satisfies
    \begin{itemize}
        \item[P1.] $\operatorname{dist}(\sigma) = (2k-1)^2$;
        \item[P2.] $\sigma(k-1) = k^2 - k$; 
        \item[P3.] $\sigma(0) = 0$ and $\sigma(k^2-1) = k^2-1$;
        \item[P4.] $\operatorname{disp}(\sigma) = (k-1)^2$.
    \end{itemize}
    If $\sigma$ exists then it is easy to verify that the following definition of $\pi$ suffices: for every $r \in [0, k^2-1]$, set 
    $$\pi\biggl(\bigl(i-k+1\bigr) + r\biggr) = \bigl(i-k+1\bigr) + \sigma(r).$$
    It remains to verify the existence of such $\sigma$ and this will yield the idea of oscillations.

    \begin{claim}
        There exists a permutation $\sigma$ of the set $[0, k^2-1]$ satisfying P1, P2, P3, and P4.
    \end{claim}

    \begin{proof}[Sketch of the proof]
        Let $M = [0, k-1]^2$ denote the set of pairs of size $k^2$. Let the ``snake'' bijection $s:M \to [0, k^2-1]$ be defined as follows: for $a, b \in [0, k-1]$ set $s(a, b) = ak+b$ for $a$ being even and set $s(a, b) = (a+1)k-b-1$ for $a$ being odd. It is easy to verify that $s$ is indeed a bijection. Let $\operatorname{inv}:M \to M$ be the symmetry bijection, i.e. $\operatorname{inv}(a, b) = (b, a)$. The left-hand side of \Cref{fig:snake} represents the horizontal ``snake'' bijection, i.e. the set $M$ is represented as a square and every cell $(a,b)$ of the square contains the value of $s(a,b)$ inside. At the same time the right-hand side of \Cref{fig:snake} represents the vertical snake, i.e. that is the depiction of map $s \cdot \operatorname{inv}$. Now, every cell of the square has two values written inside, given by $s$ and $s \cdot \operatorname{inv}$ respectively. Consider the map $\sigma$ so that for every $(a, b) \in M$ it maps $s(a, b)$ to $(s \circ \operatorname{inv})(a,b)$ (see, \Cref{fig:sigma_from_snake}), i.e. for all $r \in [0, k^2-1]$
        $$\sigma(r) = (s \circ \operatorname{inv}\circ s^{-1}) (r).$$
        
        It is easy to verify that $\sigma$ possess properties P1-P4 by picture, see \Cref{fig:snake}.

        We say that $\sigma$ oscillates (see, \Cref{fig:snake-permutation}) because if we follow the motion of $\sigma(r)$ as $r$ increases then it makes oscillating movements --- it first monotonically increases then decreases then increases again and so on. In addition, all phases of subsequent increments/decrements take about the same amount of time of length $\Theta(k)$ and, also, every step $|\sigma(r+1)-\sigma(r)|$ is at most $\Theta(k)$.
    \end{proof}

\end{proof}

\begin{figure}[h]
    \centering
    \begin{tikzpicture}[
        x=1.15cm,
        y=1.15cm,
        grid point/.style={
            circle,
            draw=black,
            fill=white,
            minimum size=8mm,
            inner sep=0pt,
            font=\scriptsize
        }
    ]
        \begin{scope}
            \draw[blue!65!black, line width=1.5pt, rounded corners=3pt]
                (0,0) -- (4,0) -- (4,1) -- (0,1) -- (0,2)
                -- (4,2) -- (4,3) -- (0,3) -- (0,4) -- (4,4);

            \foreach \a in {0,...,4} {
                \foreach \b in {0,...,4} {
                    \ifodd\a
                        \pgfmathtruncatemacro{\r}{5*(\a+1)-\b-1}
                    \else
                        \pgfmathtruncatemacro{\r}{5*\a+\b}
                    \fi
                    \node[grid point] at (\b,\a) {$\r$};
                }
            }
        \end{scope}

        \begin{scope}[xshift=7cm]
            \draw[blue!65!black, line width=1.5pt, rounded corners=3pt]
                (0,0) -- (0,4) -- (1,4) -- (1,0) -- (2,0)
                -- (2,4) -- (3,4) -- (3,0) -- (4,0) -- (4,4);

            \foreach \a in {0,...,4} {
                \foreach \b in {0,...,4} {
                    \ifodd\b
                        \pgfmathtruncatemacro{\r}{5*(\b+1)-\a-1}
                    \else
                        \pgfmathtruncatemacro{\r}{5*\b+\a}
                    \fi
                    \node[grid point] at (\b,\a) {$\r$};
                }
            }
        \end{scope}
    \end{tikzpicture}
    \caption{The row-wise and column-wise snake orderings for $k=5$.}
    \label{fig:snake}

\end{figure}
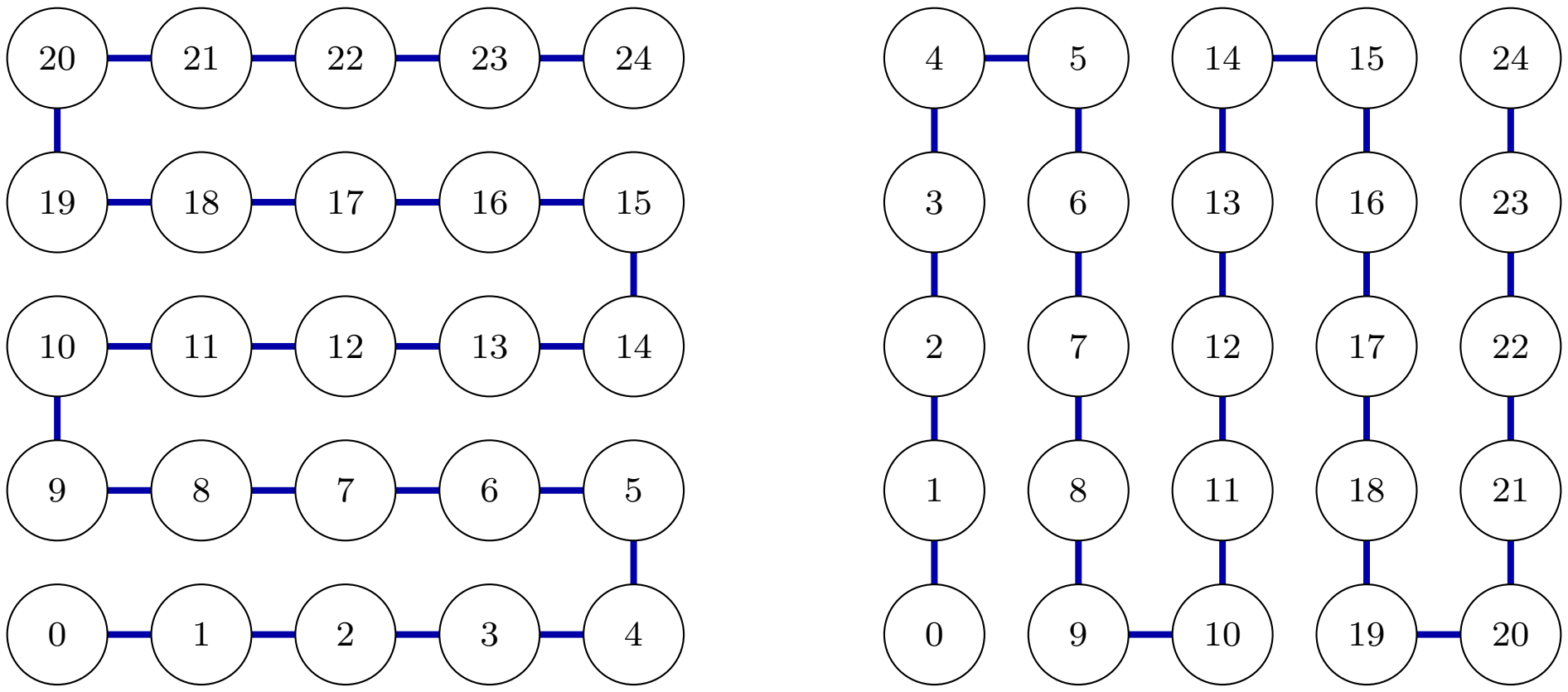

\end{document}